\documentclass[a4paper,10p]{article}
\usepackage[utf8]{inputenc}
\usepackage[english]{babel} 
\usepackage{graphicx}
\usepackage{amssymb}
\usepackage{amsthm}
\usepackage{amsmath}
\usepackage{listings}

\usepackage{makeidx}
\usepackage{xcolor}
\usepackage{amsfonts}
\usepackage{subfig}
\usepackage{float}
\usepackage{appendix}
\usepackage{cite}
\usepackage{multirow}

\usepackage{pgf,tikz,pgfplots}
\pgfplotsset{compat=1.15}
\usepackage{mathrsfs}
\usepackage{comment}
\usetikzlibrary{arrows}

\newtheorem{theo}{Theorem}
\newtheorem{defi}{Definition}

\newtheorem{alg}{Algorithm}

\usepackage{titlesec}

\titleformat{\paragraph}
{\normalfont\normalsize\bfseries}{\theparagraph}{1em}{}
\titlespacing*{\paragraph}
{0pt}{3.25ex plus 1ex minus .2ex}{1.5ex plus .2ex}
\begin{document}

%%%%%%%%% TITLE of abstract %%%%%%%%%

\begin{center}
\textbf{\large High-order well-balanced methods for systems of balance laws: a control-based approach}
\end{center}

%\medskip

%%%%%%%%%  SPEAKER %%%%%%%%%

\begin{center}
I. G\'omez, M.J. Castro, C. Par\'es\\
\emph{University of M\'alaga.}\\
\end{center}

\begin{abstract}
In some previous works, two of the authors have introduced a strategy to develop high-order numerical methods for systems of balance laws that preserve all the stationary solutions of the system. The key ingredient of these methods is a well-balanced reconstruction operator. A strategy has been also introduced to modify any standard reconstruction operator like MUSCL, ENO, CWENO, etc. in order to be well-balanced. This strategy involves a non-linear problem at every cell at every time step that consists in finding the stationary solution whose average is the given cell value. So far this strategy has been only applied to systems whose stationary solution are known either in explicit or implicit form. The goal of this paper is to present a general implementation of this technique that can be applied to any system of balance laws. To do this, the nonlinear problems to be solved in the reconstruction procedure are interpreted as control problems: they consist in finding a solution of an ODE system whose average at the computation interval is given. These problems are written in functional form and the gradient of the functional is computed on the basis of the adjoint problem. Newton's method is applied then to solve the problems. Special care is put to analyze the effects of computing the averages and the source terms using quadrature formulas. To test their efficiency and well-balancedness, the methods are applied to a number of systems of balance laws, ranging from easy academic systems consisting of Burgers equation with some nonlinear source terms to the shallow water equations or Euler equations of gas dynamics with gravity effects. 
\end{abstract}

\section{Introduction}
Let us consider a PDE system of the form:
\begin{equation} \label{sle}
U_t(x,t)+ f(U(x,t))_x= S(U(x,t)) H_x(x), \quad x\in\mathbb{R}, \, t>0,
\end{equation}
where $U(x,t)$ takes values on an open convex set $\Omega \subset \mathbb{R}^N$, $f:\Omega \longrightarrow \mathbb{R}^N$ is the flux function, $S:\Omega \longrightarrow \mathbb{R}^N$, and $H$ is a continuous known function from $\mathbb{R}$ to $\mathbb{R}$ (possibly the identity function $H(x)=x$). It is supposed that system (\ref{sle}) is strictly hyperbolic, that is,
$D_f(U)= \displaystyle \frac{\partial f}{\partial U} (U)$ has $N$ real different eigenvalues and eigenvectors. Moreover, we suppose that  $\lambda_i (U) \neq 0 , \, i=1, \ldots, N$. 

Systems of the form \eqref{sle} have non trivial stationary solutions that satisfy the ODE system:
\begin{equation}\label{sblst}
f(U)_x = S(U)H_x.
\end{equation}
A numerical method is said to be well-balanced if  it solves exactly or with enhanced accuracy all the stationary solutions
of the system or, at least, a relevant family of them. The use of methods with this property is of major importance when the waves generated for
small perturbations of a steady state are going to be simulated: this is the case, for instance, for tsunami waves in the Ocean. 
Well-balanced methods have been studied by many authors: see, for instance,
\cite{Audusse},    \cite{Bouchut},
\cite{bicaper},  \cite{EFN-03}, \cite{EFN-CRAS}, 
\cite{Chandrashekar15}, \cite{Chandrashekar17}, \cite{Desveaux16}, \cite{desveauxRipa2016},
\cite{Gosse00}, \cite{Gosse01}, \cite{Gosse02}, \cite{Kappeli14}, \cite{LR96},
\cite{LR97}, \cite{LeVeque98}, \cite{Lukacova}, \cite{Noelle06},
\cite{Noelle07}, \cite{Pelanti}, \cite{PerSi1}, \cite{PerSi2},
  \cite{Russo}, \cite{Tang04},  
 \cite{touma2016}, \cite{Xing06},... See \cite{handbook} and its references for a recent review on this topic.

Recently, in \cite{CastroPares2019} the following family of semidiscrete high-order well-balanced finite-volume methods for \eqref{sle} has been discussed:
\begin{equation} \label{met_num}
\frac{dU_i}{dt}= -\frac{1}{\Delta x} \left( F_{i+\frac{1}{2}}(t) - F_{i-\frac{1}{2}}(t) \right) + \frac{1}{\Delta x} \int_{x_{i-\frac{1}{2}}}^{x_{i+\frac{1}{2}}} S(P_i^t (x)) H_x(x) \, dx ,
\end{equation}
where: 
\begin{itemize}
\item  $I_i=\left[x_{i-\frac{1}{2}},x_{i+\frac{1}{2}}\right]$ are the computational cells, whose length $\Delta x$ is supposed to be constant for simplicity;
\item $U_i(t)$ is the approximation of the average of the exact solution at the $i$th cell at time $t$, that is, 
\begin{equation*}
U_i(t) \cong \frac{1}{\Delta x} \int_{x_{i-\frac{1}{2}}}^{x_{i+\frac{1}{2}}}  U(x,t) \, dx;
\end{equation*}
\item $P_i^t (x)$ is the approximation of the solution at the $i$th cell given by a high-order reconstruction operator from the sequence of cell averages $\{ U_i(t) \}$:
\begin{equation*}
P_i^t(x)= P_i(x; \{U_j(t)\}_{j \in \mathcal{S}_i});
\end{equation*}
where $\mathcal{S}_i$ denotes the set of indexes of the cells belonging to the stencil of the $i$th cell.
\item $ F_{i+\frac{1}{2}} = \mathbb{F}(U_{i+\frac{1}{2}}^{t, -}, U_{i+\frac{1}{2}}^{t, +})$, where $U_{i+\frac{1}{2}}^{t, \pm}$ are the reconstructed states at the intercells, i.e.
\begin{equation*}
U_{i+\frac{1}{2}}^{t, -}=P_i^t(x_{i+\frac{1}{2}}), \quad U_{i+\frac{1}{2}}^{t, +}=P_{i+1}^t(x_{i+\frac{1}{2}}),
\end{equation*}
and $\mathbb{F}$ is a consistent first order numerical flux.
\end{itemize}

It can be then easily shown that, if the reconstruction operator is well-balanced for a stationary solution $U$ of \eqref{sle} then the numerical method is also well-balanced for $U$
according to the following definitions: 
\begin{defi}
Given a stationary solution $U$ of (\ref{sle}):
\begin{itemize}
\item The numerical method (\ref{met_num}) is said to be well-balanced for $U$ if the vector of cell averages of $U$ is an equilibrium of the ODE system (\ref{met_num}).
\item The reconstruction operator is said to be well-balanced for $U$ if
\begin{equation}
P_i(x)=U(x), \quad  \forall x \in [x_{i-\frac{1}{2}}, x_{i+\frac{1}{2}}], \, \forall i,
\end{equation}
where $P_i$ is the approximation of $U$ obtained by applying the reconstruction operator to the vector of cell averages of $U$.
\end{itemize}
\end{defi}

The following strategy to design a well-balanced reconstruction operator $P_i$ on the basis of a standard operator $Q_i$ was introduced in \cite{sinum2008}: given a family of cell values $\{U_i\}$, at every cell $I_i=[x_{i-\frac{1}{2}}, x_{i+\frac{1}{2}}]$:

\begin{enumerate}
\item Look for the stationary solution $U_i^*(x)$ such that: 
\begin{equation} \label{step1}
\frac{1}{\Delta x} \int_{x_{i-\frac{1}{2}}}^{x_{i+\frac{1}{2}}} U_i^* (x) \, dx = U_i.
\end{equation}
\item Apply the reconstruction operator to the cell values $\{V_j\}_{j \in S_i}$ given by
\begin{equation*}
V_j= U_j - \frac{1}{\Delta x} \int_{x_{j-\frac{1}{2}}}^{x_{j+\frac{1}{2}}} U_i^* (x) \, dx ,
\end{equation*}
to obtain:
\begin{equation*}
Q_i(x)=Q_i(x;\{V_j\}_{j \in \mathcal{S}_i}).
\end{equation*}
\item Define
\begin{equation} \label{step3}
P_i(x)=U_i^*(x)+Q_i(x).
\end{equation}
\end{enumerate}

It can be then easily shown that the reconstruction operator $P_i$ in (\ref{step3}) is well-balanced for every stationary solution provided that the reconstruction operator $Q_i$ is exact for the null function. Moreover, if $Q_i$ is conservative, then $P_i$ is conservative, that is,
\begin{equation*}
\frac{1}{\Delta x} \int_{x_{i-\frac{1}{2}}}^{x_{i+\frac{1}{2}}} P_i (x) \, dx = U_i, \, \forall i,
\end{equation*} 
and $P_i$ is high-order accurate provided that the stationary solutions are smooth.

The well-balanced property of the method can be lost if a quadrature formula is used to compute the integral appearing at right-hand side of \eqref{met_num}. In order to circumvent this difficulty,  the authors of \cite{CastroPares2019} proposed to rewrite the methods as follows
\begin{equation}\label{met_num2}
\begin{split}
\frac{dU_i}{dt} &=-\frac{1}{\Delta x} \left( F_{i+\frac{1}{2}}(t)-f\left( U_i^{t,*}(x_{i+\frac{1}{2}})\right) - F_{i-\frac{1}{2}}(t) +f\left( U_i^{t,*}(x_{i-\frac{1}{2}})\right) \right) \\
& + \frac{1}{\Delta x} \int_{x_{i-\frac{1}{2}}}^{x_{i+\frac{1}{2}}} \left(\left[ S(P_i^t (x)) -S(U_i^{t,*}(x)) \right] H_x(x) \right) \, dx .
\end{split}
\end{equation}
where $U_i^{t,*}$ is the stationary solution found in (\ref{step1}) at the $i$th cell and time $t$. Once the method has been rewritten in this form,  a quadrature formula can be applied  to the integral without losing the well-balanced property.

The main difficulty when this strategy is applied comes from the first step of the well-balanced reconstruction operator: a nonlinear problem of the form \eqref{step1} has to be solved at every time step. Since the stationary solutions of \eqref{sle} are the solutions of the ODE system
\begin{equation}\label{state_equation}
f(U)_x= S(U) \, H_x,
\end{equation}
then, problem \eqref{step1} is equivalent to find the solution of an ODE system with prescribed average in the integration domain. In some cases, the explicit form of the general solution of the ODE is known and \eqref{step1} can be solved by hand or by using standard iterative methods for nonlinear problems: this is the case, for instance, of the shallow water equations (see \cite{CastroPares2019} and its references).

The goal of this paper is to describe a general methodology to solve numerically  problems of the form \eqref{step1} and to apply it to the implementation of well-balanced reconstruction operators for general systems of balance laws whether or not the analytical expression of the stationary solutions is known. 

The organization of the article is as follows: in Section 2 problem \eqref{step1} is interpreted as a control problem: it is first written in functional form; then, the gradient of the functional is computed using the adjoint equation. Once the expression of the gradient is available, Newton's method can be applied to solve numerically the problem: this is done in Section 3. In practice, the state and the adjoint equations, the cell-averages, and the integral appearing at the source terms are computed numerically: Section 4 is devoted to describe the well-balanced reconstruction operator and the numerical method taking into account this fact. The well-balanced property satisfied by the numerical methods is precisely stated. A number of numerical tests are presented to check the accuracy and the well-balancedness of the methods and to analyze their performance: both scalar problems and systems are considered, ranging from easy problems in which \eqref{step1} can be easily solved by hand (what will allow us to measure the computational cost of solving it by using control techniques) to systems that appear in real applications, like the shallow water model or the Euler equations of gas dynamics including gravity effects. In particular, we show that the numerical methods studied here are able to preserve subcritical and supercritical moving stationary solution: as far as we know this is the first time that such a method is obtained. Some conclusions are drawn in Section 6 and further developments are also discussed.

%
%\begin{obs}
%In the first step of the well-balanced reconstruction procedure, it has been assumed that (\ref{step1}) has always only one solution, what may not happen in practice. If (\ref{step1}) has no solution, then 
%$u_i^*\equiv 0$ is chosen in the first step and the reconstruction operator obtained is the standard one. Notice that the reconstruction operator is still well-balanced, since (\ref{step1}) has always at least one solution when the reconstruction operator is applied to the cell averages of a stationary solution. When the equation (\ref{step1}) has more than one solution, a criterion to select one of them is needed (see, for intance, \cite{diaz2013high}). 
%\end{obs}
%\begin{obs}
%The well-balanced reconstruction procedure described above can be easily adapted to obtain an operator which is well-balanced for a fixed set $\mathit{C}$ of stationary solutions: in the first step we have to look for the stationary solution of $\mathit{C}$ which satisfies (\ref{step1}) and apply the previous remark if it doesn't exist. If there is only a stationary solution $u^*$ to preserve, then the fisrt step is omitted and $u_i^*=u^*$ in steps 2 and 3.
%\end{obs}

\
\section{Control problem}

As it has been seen in the previous section, the well-balanced reconstruction procedure described in the previous section leads to find, at every cell, the solution of the ODE
system \eqref{state_equation} whose average at the cell $[x_{i-1/2}, x_{i+1/2}]$ is $U_i$. 
These problems  may have no solution or to have more than one. Observe that saying that \eqref{step1} has no solution at the $i$th cell  is equivalent to say that $U_i$ cannot be the average of any stationary solution. Therefore, at this cell there isn't any stationary solution to preserve and thus the standard reconstruction operator is applied, i.e. $U_i^* \equiv 0$ is chosen in the first step.
On the other hand, if it  has more than one solution, a criterion to select one of them is needed: see, for instance, \cite{lopez2013} where a well-balanced reconstruction operator for the shallow water equation has been introduced. We will assume in this section that 
\eqref{step1} has a unique solution.

As solving Cauchy problems is easier, we can state the problem as a control one, in which the control variable is the initial condition and the state equation is \eqref{state_equation}:

 Find the
initial condition $U_{i-1/2}$ of the Cauchy problem
\begin{equation}\label{cauchy_problem_sistema}
\begin{cases}
\displaystyle f(U)_x= S(U) \, H_x,\\
U(x_{i-\frac{1}{2}})=U_{i-1/2},
\end{cases}
\end{equation}
such that  the solution $U_i^*$ satisfies  \eqref{step1}.  

Let us write this problem in functional form.  In order to simplify the notation, let us assume that $x_{i-1/2} = 0$ and let us  denote $U_i$ by $W \in \mathbb{R}^N$. 
 The problem to be solved is then: 
 
 Find $U_0 \in \Omega$ such that
\begin{equation}\label{problem}
\mathcal{F}(U_0) = W,
\end{equation}
where $\mathcal{F}: \Omega \mapsto \mathbb{R}^N$ is given by
\begin{equation}\label{F}
\mathcal{F}(U_0) = \frac{1}{\Delta x} \int_0^{\Delta x} U(x, U_0) \, dx,
\end{equation}
where  $U(x, U_0)$  denotes the solution of the Cauchy problem 
\begin{equation}\label{cauchy_system}
\begin{cases}
\displaystyle U_x= G(U,x),\\
U(0)=U_0.
\end{cases}
\end{equation}
Here,
$G$ is the function $G:\Omega \times \mathbb{R}  \longrightarrow \mathbb{R}^N$ defined by 
\begin{equation}\label{functionG}
G(U,x)= D_f(U)^{-1}S(U) H_x. 
\end{equation} 
Remember that we assume that the eigenvalues of $D_f$ do not vanish: situations in which one of them vanishes are called resonants and they are  out of the scope of the present article.

Let us compute the gradient of $\mathcal{F}$ using the adjoint problem. To do this, given a variation $\delta \in \mathbb{R}^N $ of the initial
condition $U_0$, let us derive with respect to $s$ the function
$$
g(s)=\mathcal{F}(U_0+s \delta)  = \frac{1}{\Delta x} \int_0^{\Delta x} U(x,U_0 + s\delta ) \,  dx.
$$
For $j = 1, \dots, N$ one has:
$$
g'_j(s) = \frac{1}{\Delta x} \frac{d}{ds}\left( \int_0^{\Delta x} u_j(x, U_0 + s\delta) \, dx \right),
$$
where $U = [u_1, \dots, u_N]^T$.
In what follows, the dependency of $U$ or $u_j$ with respect to $x$ and $U_0 + s\delta$ will not be written to simplify the notation. 
Let us rewrite the integral appearing in this expression as follows:
\begin{equation}\label{ju0_system}
\begin{split}
& \int_0^{\Delta x}  u_j  \, dx = \int_0^{\Delta x} \left( u_j  + \vec\lambda_j \cdot \left( G(U,x) - U_x \right) \right) \, dx \\
&=  \int_0^{\Delta x} \left( u_j + \vec\lambda_j \cdot G(U,x) + \frac{d \vec\lambda_j}{d x} \cdot U \right)  \, dx  - \vec\lambda_j (\Delta x) \cdot U (\Delta x) + \vec\lambda_j(0) \cdot U(0),
\end{split}
\end{equation}
where $\vec\lambda_j : \mathbb{R} \longrightarrow \mathbb{R}^N$ is an arbitrary function to be selected: the so-called adjoint variables. If we denote by $\vec e_j$ the $jth$ vector of the canonical basis we get:
\begin{equation}\label{ju0'_system}
\begin{split}
 & \frac{d}{ds} \left(\int_0^{\Delta x} u_j  \, dx \right) \\
 &\qquad = \int_0^{\Delta x} \left( \frac{d u_j}{ds}  + \vec\lambda_j \cdot \frac{\partial G(U,x)}{\partial s} + 
 \frac{d \vec\lambda_j}{d x} \cdot \frac{dU}{ds}\right)  \, dx 
-\vec\lambda_j (\Delta x) \cdot \frac{dU}{ds}  (\Delta x) + \vec\lambda_j(0) \cdot \delta \\
& \qquad = \int_0^{\Delta x} \left( \vec e_j \cdot \frac{d U}{ds}  + \vec\lambda_j \cdot \left(\nabla_U G(U,x) \cdot \frac{dU}{ds} \right) + \frac{d \vec\lambda_j}{d x} \cdot \frac{dU}{ds}\right)   dx 
-\vec\lambda_j (\Delta x) \cdot \frac{dU}{ds}  (\Delta x) + \vec\lambda_j(0) \cdot \delta \\
&\qquad = \int_0^{\Delta x} \left( \left( \vec e_j  + \nabla_U G(U,x)^{T} \cdot \vec\lambda_j  + \frac{d \vec\lambda_j}{d x} \right) \cdot \frac{dU}{ds}\right)  \, dx 
-\vec\lambda_j (\Delta x) \cdot \frac{dU}{ds}  (\Delta x) + \vec \lambda_j(0) \cdot \delta, \\
\end{split}
\end{equation}
where
\begin{equation}\label{gradG}
\nabla_U G(U,x) = 
\left[
\begin{array}{ccc}
\displaystyle \frac{\partial G_1}{\partial u_1}(U,x) & \dots & \displaystyle \frac{\partial G_1}{\partial u_N}(U,x)\\
\vdots & \ddots & \vdots\\
\displaystyle \frac{\partial G_N}{\partial u_1}(U,x) & \dots & \displaystyle \frac{\partial G_N}{\partial u_N}(U,x)
\end{array}
\right].
\end{equation}
Since it is difficult to obtain an exact expression for $\displaystyle \frac{dU}{ds}$, we will choose the $j$th adjoint variable $\vec\lambda_j$ satisfying the adjoint problem:
\begin{equation}
\begin{cases} \label{adjoint_system}
\displaystyle \frac{d \vec\lambda_j}{d x}(x)=-\vec e_j  - \nabla_U G(U,x)^{T} \cdot \vec\lambda_j ,\\
\\
\vec\lambda_j(\Delta x)=0,
\end{cases}
\end{equation}
\medskip
so that (\ref{ju0'_system})  reduces to
\begin{equation} \label{ju0'_system2}
\begin{split}
g_j'(0) = \frac{1}{\Delta x}\left. \frac{d}{ds} \left(\int_0^{\Delta x} u_j  \, dx \right)\right|_{s = 0} &=  \frac{1}{\Delta x} \vec\lambda_j(0) \cdot \delta.
\end{split}
\end{equation}
\medskip
Therefore:
\begin{equation}\label{g'0_system}
g'(0)=  \frac{1}{\Delta x} \Lambda(0)^T \cdot \delta,
\end{equation}
where  $\Lambda $ denotes the matrix whose columns are $\lambda_1(x), \ldots, \lambda_N(x)$, that is,
\begin{equation}\label{Lambda}
\Lambda (x) = \begin{bmatrix}
\lambda_1(x) \, | \quad \ldots \quad | \, \lambda_N(x)
\end{bmatrix}.
\end{equation}

Therefore
$$
D\mathcal{F}(U_0) =  \frac{1}{\Delta x} \Lambda(0)^T.
$$

\section{Numerical algorithm}

\subsection{Newton's method}

Since problem \eqref{step1} has to be solved at every intercell at every time step, it is crucial to  choose an efficient numerical method. Since the gradient of $\mathcal{F}$ is available, Newton's method can be applied. Observe that a sensible choice for the initial guess $U_0^0$ is $W$: if $\Delta x$ is small, the average of the solution of the Cauchy problem is expected to be close to the initial condition.
 The algorithm is then as follows: 
 
\begin{alg}{Newton's method}
  \begin{itemize}
\item $U^0_0 = W$;

\item For k = 0,1,2\dots

\begin{itemize}

\item Compute the solution $U_k$ of \eqref{cauchy_system} with initial condition $U^k_0$ in the interval $[0, \Delta x]$.

\item For $j = 1, \dots, N$ compute the solution $\vec\lambda_j$ of \eqref{adjoint_system} with $U = U_k$ in the interval $[0, \Delta x]$.

\item Compute $V_k$ by solving the linear system:
$$
\Lambda(0)^T V_k = \Delta x (\mathcal{F}(U^k_0) - W),
$$
where $\Lambda(x)$ is given by \eqref{Lambda}. 

\item Update $U^k_0$:
$$
U^{k+1}_0 = U^k_0 - V_k.
$$

\end{itemize}

\end{itemize}

\end{alg}

At every iteration of the method $N+1$ Cauchy problems and a $N \times N$ linear system have to be solved.

The computational cost can be reduced by using the modified Newton method in which the matrix $\Lambda(0)$ is only updated every $K$ iterations, where $K$ is a fixed integer.

\subsection{Numerical integration}
In practice, the integral in the definition $\mathcal{F}$ given by \eqref{F} is computed using a quadrature rule in $[0, \Delta x]$
$$
\int_0^{\Delta x} g(x) \, dx \cong \Delta x \sum_{l= 0}^M \alpha_l g(x_l),
$$
 and the initial and final value problems to compute $U_k$ and $\vec\lambda_j$ at the iterations of the algorithms are approximated with a numerical method to solve ODE problems using a mesh of the interval $[0, \Delta x]$
whose maximum step will be denoted by $h$. This mesh will be chosen so that all the quadrature points $x_l$ are nodes. The
order of the method and the size of $h$ will be chosen so that  errors are close to machine precision.   

Therefore, in practice the algorithms solves the numerical problem:

Find $U_0$ such that
$$
\mathcal{F}_h(U_0):= \sum_{l=0}^M \alpha_l U_{h,l} = W,
$$
where $U_{h,l}$ represents the numerical approximation of $U(x, U_0)$ at the quadrature point $x_l$ given by the numerical method chosen to solve the ODE.

\section{Discrete well-balanced reconstruction operator}
To implement the well-balanced reconstruction operator the following ingredients have to be chosen first:
\begin{itemize}
\item Quadrature rules at the cells
$$
\int_{x_{i-1/2}}^{x_{i+1/2}} g(x) \, dx \cong \Delta x \sum_{l= 0}^M \alpha^i_l g(x^i_l).
$$

\item A numerical method for solving Cauchy problems.

\item  Meshes of maximum step $h$ at the cells $[x_{i-1/2}, x_{i+1/2}]$ whose set of nodes include the quadrature points $x^i_l$ and 
$x_{i\pm 1/2}$.

\end{itemize}

Once these ingredients have been chosen, the reconstruction procedure is as follows: 

Given a family of cell values $\{U_i\}$, at every cell $I_i=[x_{i-\frac{1}{2}}, x_{i+\frac{1}{2}}]$:

\begin{enumerate}
\item Look for $U^*_{i-1/2}$ such that: 
\begin{equation} \label{step1h}
\sum_{l=0}^M \alpha^i_l U^{*,i}_{h,l} = U_i,
\end{equation}
where $U^{*,i}_{h,l}$ represents the numerical approximation of 
$$
U^*_i(x) = U(x, U^*_{i-1/2})
$$
at the quadrature point $x^i_l$ given by the numerical method chosen to solve the ODE.

\item  Obtain approximations of $U^*_i$
$$
U^{*,i}_{h,j,l}, \quad l=0,\dots, M, \quad j \in \mathcal{S}_i,
$$
 at the quadrature point $x^j_l$ of the cells of the stencil  using the chosen numerical method. Notice that $U^{*,i}_{h,i,l} = U^{*,i}_{h,l}$ have been already computed at step 1.

\item Apply the standard reconstruction operator $Q_i$ to the cell values $\{V_j\}_{j \in \mathcal{S}_i}$ given by
\begin{equation*}
V_j= U_j - \sum_{l = 0}^M \alpha_l^j U^{*,i}_{h,j,l}
\end{equation*}
to obtain:
\begin{equation*}
Q_i(x)=Q_i(x;\{V_j\}_{j \in \mathcal{S}_i}).
\end{equation*}
\item Compute:
\begin{eqnarray*}
& & U^+_{i-1/2}   =   U^*_{i-1/2} + Q_i (x_{i-1/2}),\\
& & U^-_{i+1/2}  =   U^*_{i+1/2} + Q_i (x_{i+1/2}),\\
& & P^i_l   =  U^{*,i}_{h,l} + Q_i (x^i_l), \quad l = 0, \dots, M,
\end{eqnarray*}
where $  U^*_{i+1/2} $ is the approximation to $U^*_i$ provided by the chosen numerical method at $x_{i+1/2}$.
\end{enumerate}

The semidiscrete numerical method to solve \eqref{sle} writes then as follows:
\begin{equation}\label{met_num_disc}
\begin{split}
\frac{dU_i}{dt} &=-\frac{1}{\Delta x} \left( F_{i+\frac{1}{2}}(t)-f\left( U_{i+1/2}^{*}(t) \right) - F_{i-\frac{1}{2}}(t) +f\left( U_{i-1/2}^{*}(t)\right) \right) \\
& + \sum_{l = 0}^M \alpha^i_l \left(S(P^{i}_l(t)) -S(U_{h,l}^{*,i}(t)) \right) H_x(x^i_l),
\end{split}
\end{equation}
where 
 $$ F_{i+\frac{1}{2}}(t) = \mathbb{F}(U_{i+\frac{1}{2}}^{-}(t), U_{i+\frac{1}{2}}^{+}(t)). $$
The values of $U^\pm_{i \pm 1/2}(t)$,  $P^{i}_l(t)$, $U_{i\pm 1/2}^{*}(t)$, $U_{h,l}^{*,i}(t)$ are given by the well-balanced reconstruction operator applied to the cell averages $\{ U_i(t)\}$. 

The proof of the following result is straightforward:
\begin{theo}\label{th:wbc}
The numerical method \eqref{met_num_disc} is well-balanced for every stationary solution $u$ in the sense that the vector of the cell-averages $U_{h,i}$ given by
$$
U_{h,i} = \sum_{l = 0}^M \alpha^i_l U^i_{h,l}
$$
is an equilibrium of the ODE system \eqref{met_num_disc}. Here, $U^i_{h,l}$ are approximations of $U$ at the quadrature points obtained with the numerical method selected for solving 
\eqref{state_equation} in the well-balanced reconstructions  using the same meshes at the cells.
\end{theo}

\subsection{First and second order methods}
First and second order methods can be implemented in an easier way if the mid-point rule is used to approach the cell averages:
$$\frac{1}{\Delta x}\int_{x_{i-1/2}}^{x_{i+1/2}} U(x)\, dx \cong U(x_i).$$
In effect, in this case the first step in the reconstruction procedure reduces to:
\begin{enumerate}
\item Look for the stationary solution $U^*_{i}$ such that: 
\begin{equation} \label{step1h_1and2order}
 U^*_i(x_i) = U_i.
\end{equation}
\end{enumerate}
Therefore, there is no need to solve a nonlinear problem: it is enough to solve  the Cauchy problem
\begin{equation}\label{Cauchyfsom}
\left\{
\begin{array}{l}
U_x = G(U,x), \\
U(x_i) = U_i.
\end{array}
\right.
\end{equation}
A first order reconstruction operator is then given by
\begin{eqnarray*}
& & U^+_{i-1/2}   =   U^*_{i-1/2} ,\\
& & U^-_{i+1/2}  =   U^*_{i+1/2} ,\\
& & P^i_i   =  U_i, 
\end{eqnarray*}
where $U^*_{i \pm 1/2}$ represent the numerical approximations of the solution of the Cauchy problem \eqref{Cauchyfsom} at $x_{i \pm 1/2}$ computed with the chosen numerical method. It can be easily checked that, if the mid-point formula is used again to approximate the integral term in \eqref{met_num_disc}, the expression of the numerical method reduces to
\begin{equation}\label{met_num_disc_fo}
\frac{dU_i}{dt} =-\frac{1}{\Delta x} \left( F_{i+\frac{1}{2}}(t)-f\left( U_{i+1/2}^{*}(t) \right) - F_{i-\frac{1}{2}}(t) +f\left( U_{i-1/2}^{*}(t)\right) \right).
\end{equation}
For second order methods, the solution of \eqref{Cauchyfsom} has to be numerically approximated at the center of the cells of the
stencil, $x_j$, and then steps 3-4 are performed. In the particular case of the MUSCL reconstruction, the numerical method writes in the form \eqref{met_num_disc_fo}: see \cite{CastroPares2019}. 

\section{Numerical experiments} \label{numerical_experiments}

%\subsection{Implementation} \label{subsection_implementation}
In order to implement the well-balanced procedure described in the previous paragraphs, the following choices have been made:
\begin{itemize}
\item The fourth order Runge-Kutta method is selected to solve the state and the adjoint ODE problems.
\item The following quadrature rules are chosen:
\begin{itemize}
\item The midpoint rule is used in first and second order schemes. 
\item The Gauss two points quadrature rule is used in third order schemes.
\end{itemize}
\end{itemize}
Concerning the mesh at the cell $[x_{i-1/2}, x_{i+1/2}]$ to solve the ODE systems, once the quadrature points $x^i_l$ have been selected, we consider  uniform partitions of the 
intervals
$$[x_{i-1/2}, x^i_0], [x^i_0, x^i_1], \dots, [x^i_{M-1},  x^i_M],  [x^i_M, x_{i+1/2}]$$
with $N_p$ subintervals, so that the total mesh has $N_p(M+2) + 1$ points.

To implement the numerical method \eqref{met_num_disc} we consider:
\begin{itemize} 
\item Rusanov  numerical flux;

\item the second order MUSCL (see \cite{MUSCL}) and the third order CWENO reconstructions (see \cite{CWENO}, \cite{CWENO2});  

\item the third order TVD Runge-Kutta for solving the ODE system \eqref{met_num_disc}: see \cite{Gottlieb98}.

\end{itemize}

When the initial condition is a stationary solution $U^*$ in an interval $[a,b]$, we approximate its cell averages either by applying the quadrature formula to the
exact solution (when it is avalable) or by
$$
U^*_{h,i} = \sum_{l = 0}^M \alpha^i_l U^{*,i}_{h,l}
$$
where $U^{*,i}_{h,j}$ are the approximations at the quadrature points obtained using RK4 to approximate \eqref{state_equation} with initial condition
$$
U(a) = U^*(a).
$$

Observe that the only information about the particular problem required by the numerical method is $f$, $S$, $H$, $G$,  $\nabla G$ (see \eqref{sle}, \eqref{functionG},  \eqref{gradG}) what leads to very general algorithms.

The following symbols will be used in this section to denote the different methods considered:

\begin{itemize}

\item SM$i$, $i=1, 2, 3$: numerical method of order $i$ based on the Rusanov flux and the standard reconstruction operators. 

\item WBM$i$, $i=1, 2, 3$: numerical method of order $i$ based on the Rusanov flux and the well-balanced reconstruction operators in which problems \eqref{step1} are exactly solved.

\item DWBM$i$, $i=1, 2, 3$: numerical method of order $i$ based on the Rusanov flux and the well-balanced reconstruction operators in which problems \eqref{step1} are solved numerically 
(by solving the Cauchy problems \eqref{Cauchyfsom} for first and second order methods or by solving
\eqref{step1h}
using Newton's method for third order methods).

%\item DWBM$i$, $i=1, 2, 3$: numerical method of order $i$ based on the Rusanov flux and the well-balanced reconstruction operators in which problems \eqref{step1} are solved using descent methods. 

\end{itemize}

\subsection{Problem 1: Burgers equation with a nonlinear source term I}

Let us consider the Burgers equation with a non-linear source term:
\begin{equation}\label{burgers}
\begin{cases}
\displaystyle u_t + \left( \frac{u^2}{2} \right)_x= u^2 , \quad x\in \mathbb{R}, \, t>0,\\
u(x,0)=u_0(x),
\end{cases}
\end{equation}
This problem is the particular case of \eqref{sle} corresponding to:
$$U = u, \quad f(U)= \displaystyle \frac{u^2}{2}, \quad S(U)=u^2, \quad H(x)=x.$$
The ODE satisfied by the stationary solutions is
\begin{equation}\label{ODEBurgers1}
\frac{du}{dx} = u,
\end{equation}
whose solutions  are 
$$u(x)=C e^x, \quad C \in \mathbb{R}.$$
Therefore:
$$
G(u,x) = u, \quad \partial_u G(u,x) = 1.
$$

Since the expression of the stationary solutions is known, the first step of the well-balanced reconstruction procedure can be easily solved: given a family of cell values $\{u_i\}$, and given a quadrature formula, the stationary solution $u^*_i$ which solves the non-linear problem 
\[
\sum_{l=0}^M \alpha_l^i u_i^*(x_l^i)=u_i,
\]
reduces to
\begin{equation}\label{expresionpaso1}
u^*_i(x)= \frac{u_i}{\sum_{l=0}^M e^{x_l^i}} e^x
\end{equation}

Although WBM$i$,$i=1,2,3$ can be easily implemented using this explicit expression, DWB$i$, $i=1,2,3$ will be also applied to check their efficiencies and their sensitivity to the numerical discretization of the ODE \eqref{ODEBurgers1}. Since \eqref{ODEBurgers1} is a linear equation, Newton's method converges in only one iteration.

%Moreover, the descent method taking $u_i$ as the initial guess converges as well in only one iteration if the step 
%$$
%\rho_i = \frac{\Delta x^2}{2} \frac{e^{x_{i-1/2}}(e^{\Delta x} - 1) - \Delta x}{e^{x_{i-1/2}}(e^{\Delta x} - 1)^2(e^{\Delta x} - 1 - \Delta x)}
%$$
%is chosen, as it can be easily checked. Since the computational cost of both approaches is similar,  only Newton's method will be considered. 

\subsubsection{Test 1.1}

We consider $x \in [-1, 1]$ and $t \in [0,5]$. The $CFL$ parameter is set to 0.9, and the function $u_0(x)=e^x$ is chosen as initial condition. The boundary condition
$$u(-1,t) = e^{-1}$$
is imposed at $x=-1$ and open boundary conditions are set at $x = 1$. The tolerance considered to stop the Newton's method  is $\varepsilon=10^{-8}$. 

\begin{figure}[H]
\begin{center}
  \subfloat[Errors (logaritmic scale) ]{
   \includegraphics[width=0.5\textwidth]{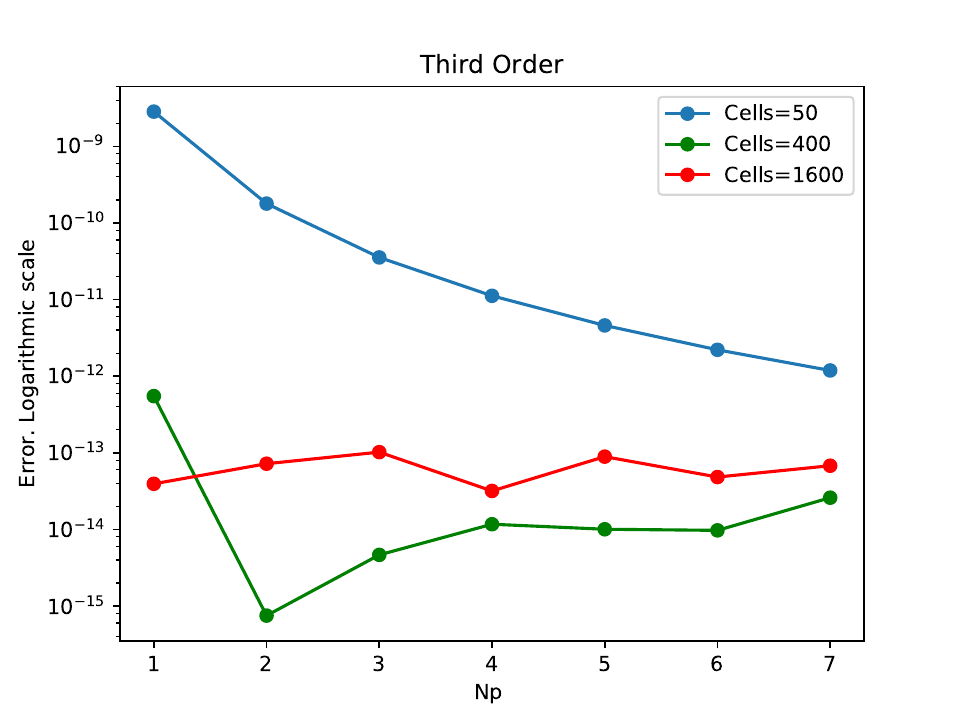}}
 \subfloat[CPU time]{
   \includegraphics[width=0.5\textwidth]{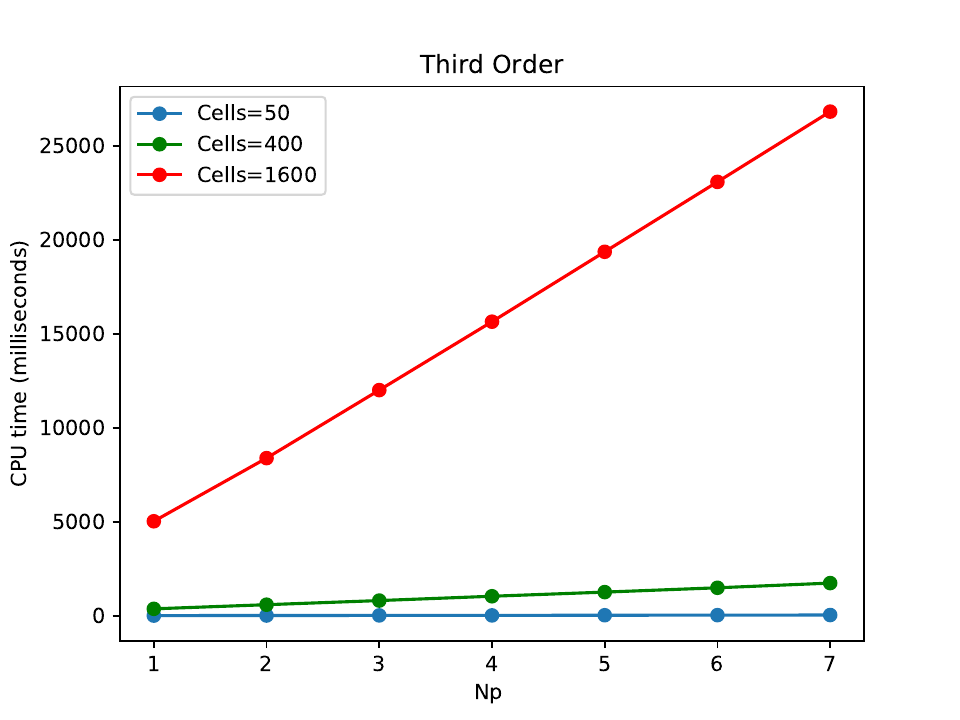}}
    \caption{Test 1.1. Errors and CPU times corresponding to DWBM3 with different number of cells and different values of $N_p$.} 
    \label{np}
 \end{center}
\end{figure}

The initial condition considered for the numerical methods is the sequence of cell averages of the exact solution computed with the quadrature formulas.  Different values for $N_p$ have been compared. Figure \ref{np} shows the errors at logarithmic scale and the CPU times corresponding to different values of $N_p$ for the third order method. As it can be seen, all the errors are below $ 10^{-12}$ except for the coarsest mesh (50 cells) and $N_p =1$ or 2.  Therefore in this case $N_p=1$ is enough for fine meshes and $N_p=2$ or 3 is a good choice in all cases. The results and conclusions are similar for the first and second order methods. 

\begin{figure}[H]
\begin{center}
  \subfloat[ SM$i$, $i =1,2,3$ ]{
   \includegraphics[width=0.5\textwidth]{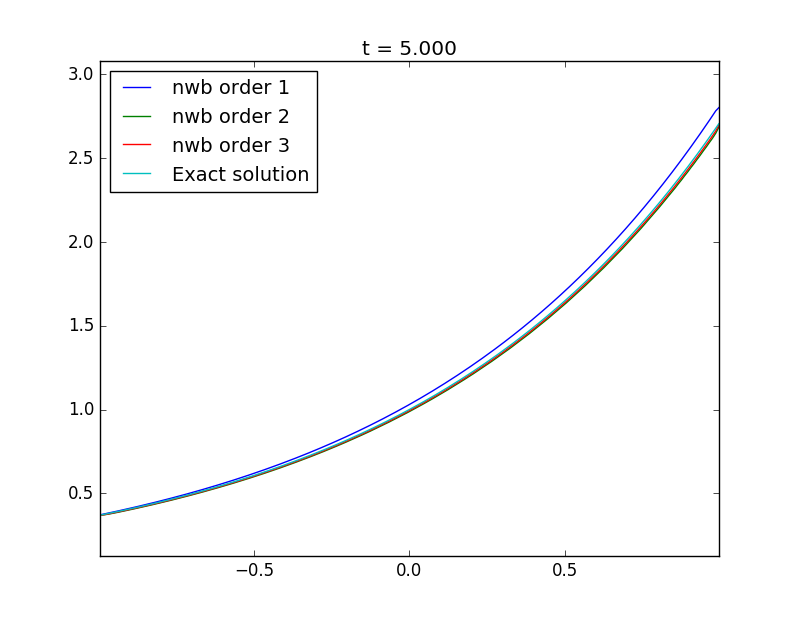}}
 \subfloat[WBM$i$, $i = 1,2,3$)]{
   \includegraphics[width=0.5\textwidth]{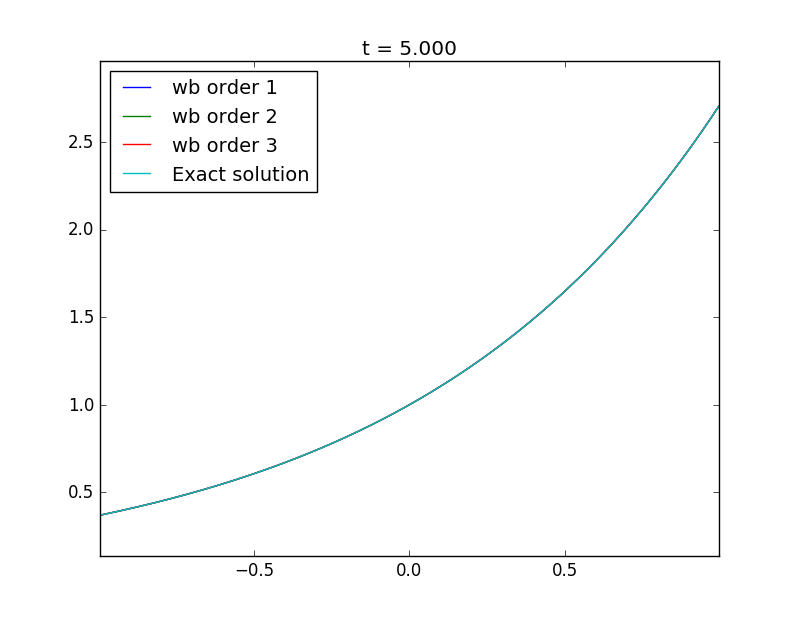}}
    \caption{Test 1.1. Numerical solutions at $t = 5s$.
Number of cells: 200.} \label{test1_1}
 \end{center}
\end{figure}

Figure \ref{test1_1} shows the numerical solutions obtained with SM$i=1,2,3$ and WBM$i$, $i=1,2,3$ (the graphs corresponding to DWBM$i$, $i=1,2,3$ are similar).
Tables \ref{ex1_error_nwb}, \ref{ex1_error_wb_exacta}, and \ref{ex1_error_wb_newton} show the errors corresponding to SM$i$, WBM$i$, DWBM$i$, $i = 1,2,3$. 
\begin{table}[H]
\centering
\begin{tabular}{|c|cc|cc|cc|} \hline
Cells & SM1: Error &Order& SM2: Error & Order&SM3: Error&Order\\\hline
100& 7.53E-2 & - & 2.44E-3  & - & 7.66E-6 & -\\
200 & 3.78E-2 & 0.995 & 8.09E-4 & 1.591 & 9.62E-7 & 2.993\\
400 & 1.89E-2 & 1.002 & 2.16E-4 & 1.905 & 1.21E-7 & 2.995\\
800 & 9.43E-3 & 1.0001 & 5.54E-5 & 1.963 & 1.51E-8 & 2.998\\ \hline
\end{tabular}
\caption{Test 1.1. Errors in $L^1$ norm and convergence rates for SM$i$, $Ii= 1,2,3$.} \label{ex1_error_nwb}

\end{table}

\begin{table}[H]
\centering
\begin{tabular}{|c|c|c|c|} \hline
Cells & WBM1: Error &WBM2: Error &WBM3: Error \\\hline
100& 4.21E-15 &  8.87E-16  &  3.20E-16 \\
200 & 2.90E-15  & 4.42E-16 &   2.54E-16 \\
400 & 1.84E-14 & 1.82E-15 &   7.40E-14 \\
800 & 4.45E-16 &  1.83E-16 &  2.61E-15 \\  \hline
\end{tabular}
\caption{Test 1.1. Errors in $L^1$ norm  for WBM$i$, $i= 1,2,3$.} \label{ex1_error_wb_exacta}

\end{table}

\begin{table}[H]
\centering
\begin{tabular}{|c|cc|cc|cc|} \hline
Cells &\multicolumn{2}{|c|}{DWBM1: Error }& \multicolumn{2}{|c|}{DWBM2: Error)}&\multicolumn{2}{|c|}{DWBM3: Error  }\\
 & $N_p=1$ &$N_p=3$&$N_p=1$&$N_p=3$&$N_p=1$&$N_p=3$\\\hline
100 & 1.70E-10 & 2.10E-12 & 1.74E-10 & 2.14E-12  & 1.79E-10 & 2.20E-12 \\
200 & 1.07E-11 & 1.29E-13 & 1.08E-11  &  1.09E-13 & 1.11E-11 & 5.03E-14 \\
400 & 6.72E-13 & 1.04E-14 & 6.72E-13 &  4.86E-16 & 5.50E-13 & 4.65E-15 \\
800 & 1.77E-14 & 2.85E-15 & 5.13E-16 & 4.61E-16 & 1.98E-15 & 1.59E-14 \\  \hline
\end{tabular}
\caption{Test 1.1. Errors in $L^1$ norm  for DWBM$i$, $i=1,2,3$} \label{ex1_error_wb_newton}

\end{table}

Notice that the errors for SM$i$, $i=1,2,3$ decrease  with the number of cells at the expected rate. While WBM$i$, $i=1,2,3$ capture the exact solution with machine precision, the errors for DWBM$i$, $i=1,2,3$ depends on the tolerance used in Newton's method  and on the discretization error corresponding to RK4, whose order is $O(h^4)$. The computational costs are shown in Table \ref{ex1_times}. It can be seen that the well-balanced  modification of the reconstruction operator based on the exact solution of \eqref{step1}  multiplies the computational cost by a factor ranging from 1.5 to 7.5. On the other hand, the numerical resolution of \eqref{step1} increases the computational cost of the well-balanced methods by a factor of 1--1.5 if $N_p = 1$. This extra cost increases linearly with $N_p$.

\begin{table}[H]
\centering
\begin{tabular}{|c|c|c|c|cc|}
\hline 
     Cells              & $i$ & SM$i$ & WBM$i$ & \multicolumn{2}{|c|}{ DWBM$i$} \\  
   & & & & $N_p=1$ & $N_p=3$             \\   \hline
\multirow{3}{*}{100} & 1 & 20 & 30 & 30 &40\\ %\cline{2-5} 
                  & 2 & 30 & 60 & 70 & 140 \\ %\cline{2-5} 
                  & 3 & 40 & 190 & 200 & 380 \\ \hline
\multirow{3}{*}{200} & 1 & 20 & 60 & 80 & 100 \\ %\cline{2-5} 
                  & 2 & 40 & 190 & 200 & 430 \\ %\cline{2-5} 
                  & 3 & 110 & 480 & 580 & 1170 \\ \hline         
\multirow{3}{*}{400} & 1& 50 & 180 & 220 & 330 \\ %\cline{2-5} 
                  & 2 & 100 & 530 & 610 & 1250 \\ %\cline{2-5} 
                  & 3 & 350 & 1680 & 1950 & 3820 \\ \hline
\multirow{3}{*}{800} & 1 & 140 & 570 & 650 & 810 \\ %\cline{2-5} 
                  & 2 & 270 & 2040 & 2080 & 4190 \\ %\cline{2-5} 
                  & 3& 1080 & 5540 & 6360 & 14970 \\ \hline                  
\end{tabular}
\caption{Test 1.1. Computational time (milliseconds).}\label{ex1_times}
\end{table}

\subsubsection{Test 1.2}
\begin{figure}[H]
\begin{center}
  \subfloat{
   \includegraphics[width=0.5\textwidth]{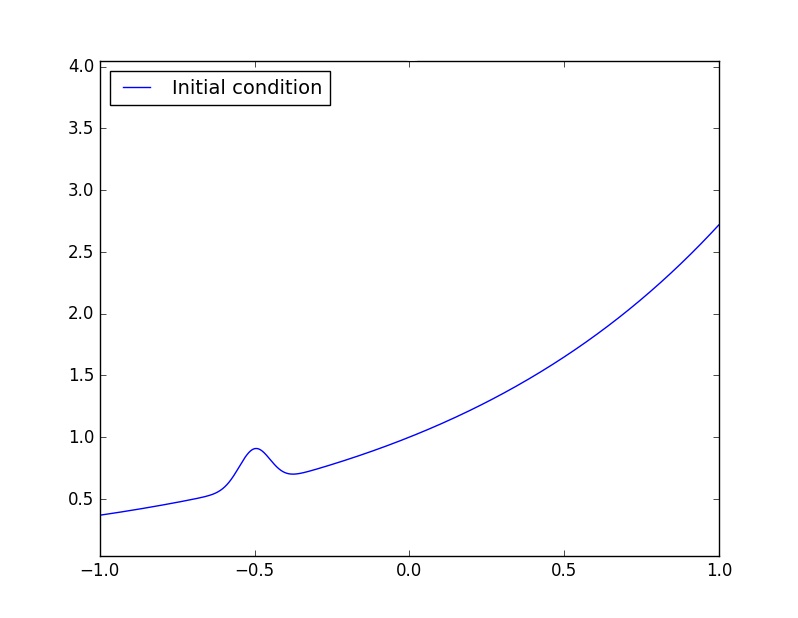}}
 \caption{Test 1.2. Initial condition.} \label{test2_cini}
 \end{center}
\end{figure}

The evolution of a perturbation of the stationary solution considered in the previous test is now simulated. The only difference with Test 1.1. is that, here, the initial condition is given by:  
$$u_0(x)=e^x+0.3e^{-200(x+0.5)^2},$$
see Figure \ref{test2_cini}. $N_p=1$ has been considered for DWBM$i$, $i = 1, 2, 3$.

\begin{figure}[H]
\begin{center}
  \subfloat[SM$i$, $i=1,2,3$. $t = 0.5s$. ]{
   \includegraphics[width=0.5\textwidth]{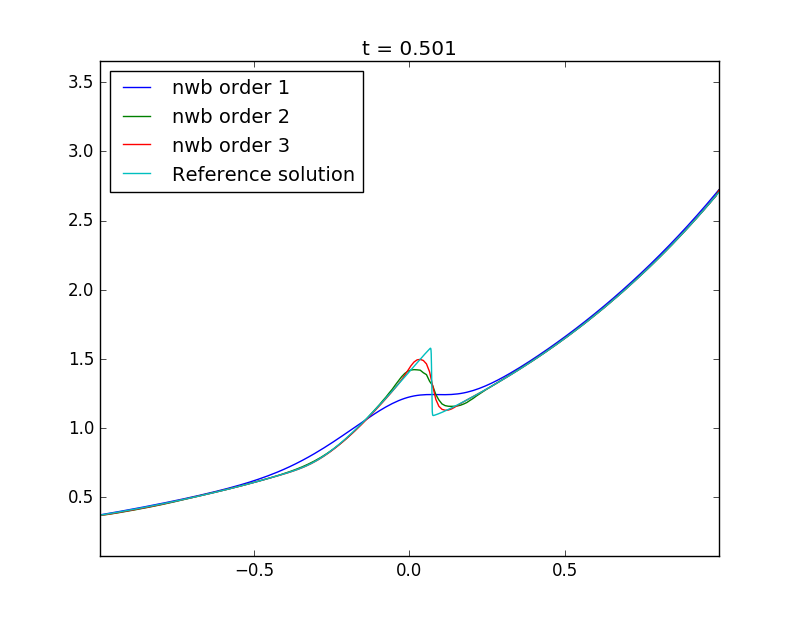}}
 \subfloat[WBM$i$, $i=1,2,3$. $t = 0.5s$.]{
   \includegraphics[width=0.5\textwidth]{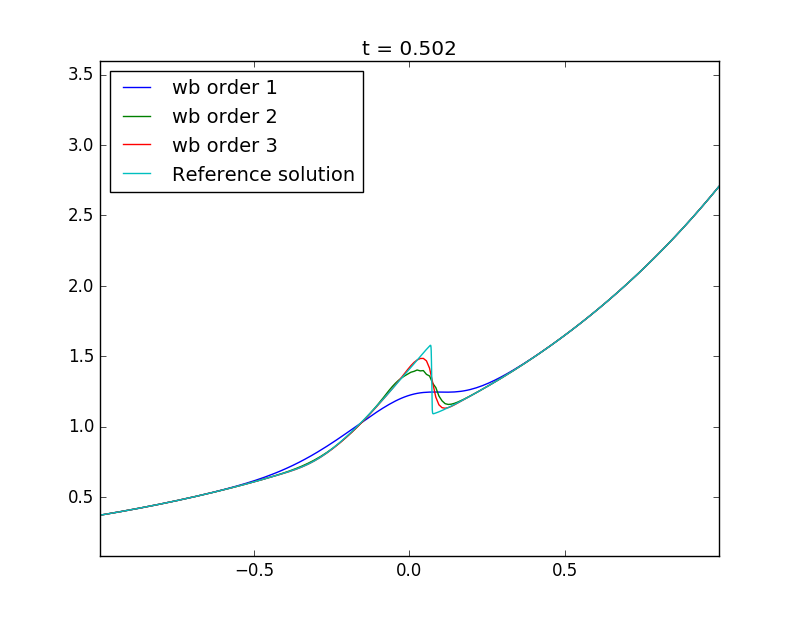}}\vspace{0.0001mm}
   \subfloat[SM$i$, $i$=1,2,3. $t = 1s$.]{
   \includegraphics[width=0.5\textwidth]{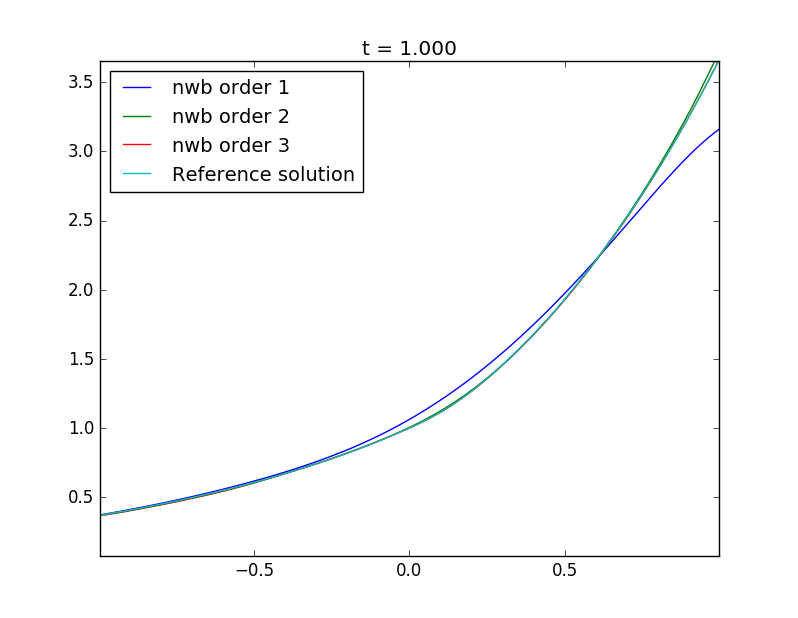}}
   \subfloat[WBM$i$, $i=1,2,3$.  $t = 1s$.]{
   \includegraphics[width=0.5\textwidth]{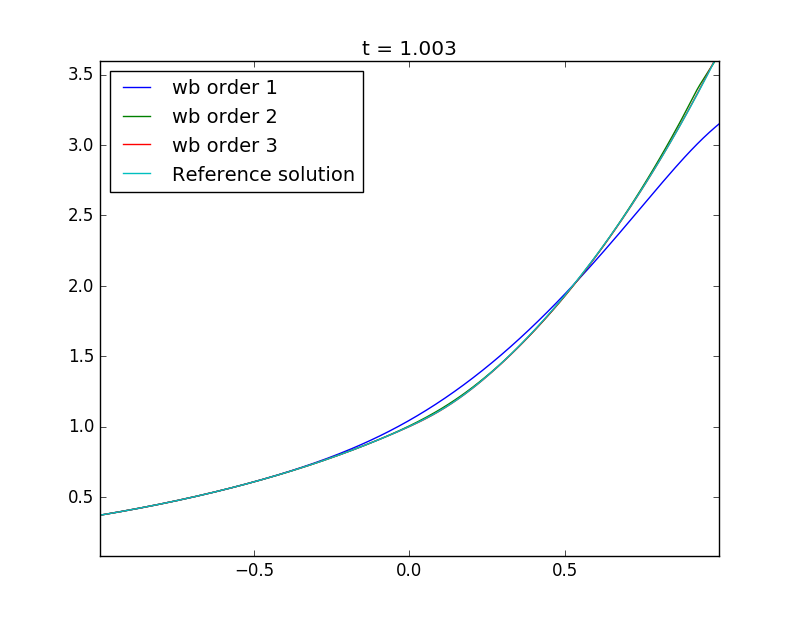}}\vspace{0.0001mm}
   \subfloat[SM$i$, $i$=1,2,3. $t = 5s$. ]{
   \includegraphics[width=0.5\textwidth]{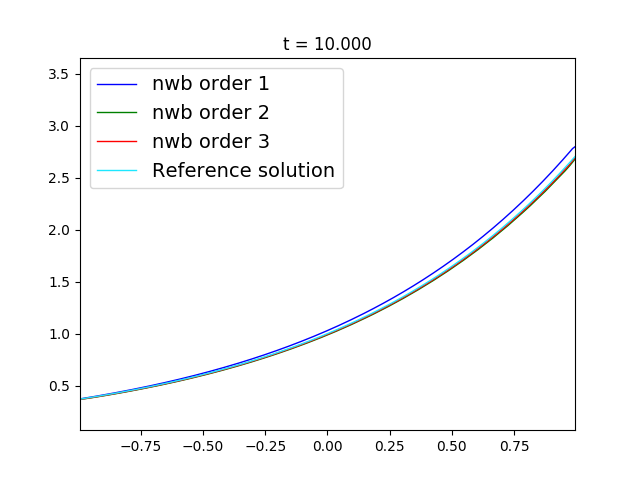}}
 \subfloat[WBM$i$, $i=1,2,3$.  $t = 5s$.]{
   \includegraphics[width=0.5\textwidth]{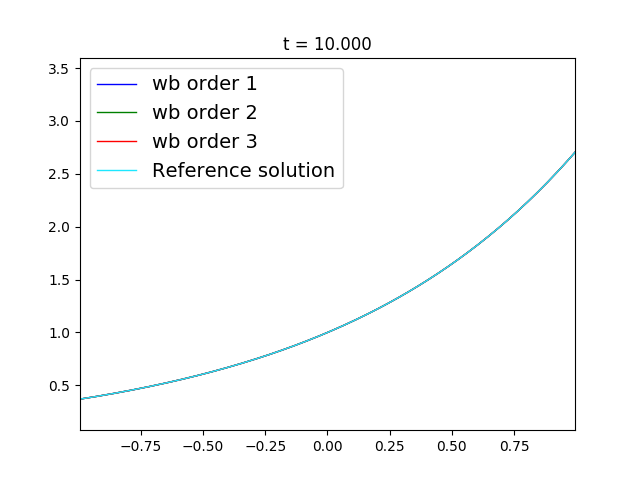}}
    \caption{Test 1.1. Reference and numerical solutions at $t = 0.5, 1, 10s$. Number of cells: 200.} \label{test2}
 \end{center}
\end{figure}

Figure \ref{test2} shows the propagation of the perturbation at times $t=0.5,1, 10s$ given by SM$i$ and WBM$i$, $i=1,2,3$ (the graphs corresponding to DWBM$i$, $i=1,2,3$ are similar). A reference solution has been computed with WBM1 using a fine mesh.

Although during the propagation of the perturbation there are not important differences between the well-balanced and non well-balanced schemes, once the propagation has left the domain  the well-balanced methods recover the stationary solution, as expected, while the non-well balanced methods perturb it. Tables \ref{ex2_error_nwb}, \ref{ex2_error_wb_exacta}, and \ref{ex2_error_wb_newton} show the errors corresponding to the different methods at time $t = 10s.$  
%Therefore, conclusions concerning to computational cost are similar to the ones included in table \ref{ex1_times}.

\begin{table}[H]
\centering
\begin{tabular}{|c|c|c|c|} \hline
Cells & Error ($i=1$)&Error ($i=2$)& Error ($i=3$)\\\hline
100& 1.58E-1 &  5.28E-2  &  3.12E-2\\
200 & 7.51E-2 &  2.83E-2 &  1.61E-2\\
400 & 3.66E-2 &  1.47E-2 &  8.16E-3\\
800 & 1.81E-2 &  7.51E-3 &  4.11E-3\\ \hline
\end{tabular}
\caption{Test 1.2. Errors in $L^1$ norm for SM$i$, $i=1,2,3$. $t = 10s$.} \label{ex2_error_nwb}

\end{table}

\begin{table}[H]
\centering
\begin{tabular}{|c|c|c|c|} \hline
Cells & Error ($i = 1$)&Error ($i = 2$)& Error ($i=3$)\\\hline
100& 5.03E-11 &  1.37E-14  &  1.44E-14\\
200 & 3.54E-14 &  1.75E-14 &  4.44E-14\\
400 & 8.49E-14 &  9.09E-14 &  2.51E-13\\
800 & 3.59E-15 &  6.72E-14 &  3.80E-13\\ \hline
\end{tabular}
\caption{Test 1.2. Errors in $L^1$ norm for WBM$i$, $i=1,2,3$. $t = 10s$.} \label{ex2_error_wb_exacta}

\end{table}

\begin{table}[H]
\centering
\begin{tabular}{|c|c|c|c|} \hline
Cells & Error ($i=1$)&Error ($i=2$)& Error ($i=3$)\\\hline
100& 2.57E-10 &  5.57E-10  &  2.65E-10\\
200 & 1.61E-11 &  1.60E-11 &  1.65E-11\\
400 &9.73E-13 &  6.73E-13 &  8.32E-13\\
800 & 7.31E-14 &  4.86E-15 &  1.69E-13\\ \hline
\end{tabular}
\caption{Test 1.2. Errors in $L^1$ norm for DWBM$i$, $i=1,2,3$. $t = 10s$.} \label{ex2_error_wb_newton}

\end{table}

\subsection{Problem 2: Burgers equation with a nonlinear source term II}
We now consider Burgers equation with a different non-linear source term:
\begin{equation} \label{burgers2}
\begin{cases}
u_t +  \left( \displaystyle \frac{u^2}{2} \right)_x= \sin(u), \quad x\in \mathbb{R}, \, t>0,\\
u(x,0)=u_0(x).
\end{cases}
\end{equation}
This problem is the particular case of \eqref{sle} corresponding to:
$$U = u, \quad f(U)= \displaystyle \frac{u^2}{2}, \quad S(U)=\sin(u), \quad H(x)=x.$$
The ODE satisfied by the stationary solutions is
\begin{equation}\label{ODEBurgers2}
\frac{du}{dx} = \frac{ \sin(u)}{u}.
\end{equation}
Therefore:
$$
G(x,U) = \frac{\sin(u)}{u}, \quad \partial_U G(x,U) = \frac{u\cos(u) - \sin(u)}{u^2}.
$$
In this case, the stationary solutions cannot be expressed in terms of elementary functions so that \eqref{step1} has to be numerically solved.

\subsubsection{Test 2.1}

\begin{figure}[H]
\begin{center}
  \subfloat{
   \includegraphics[width=0.5\textwidth]{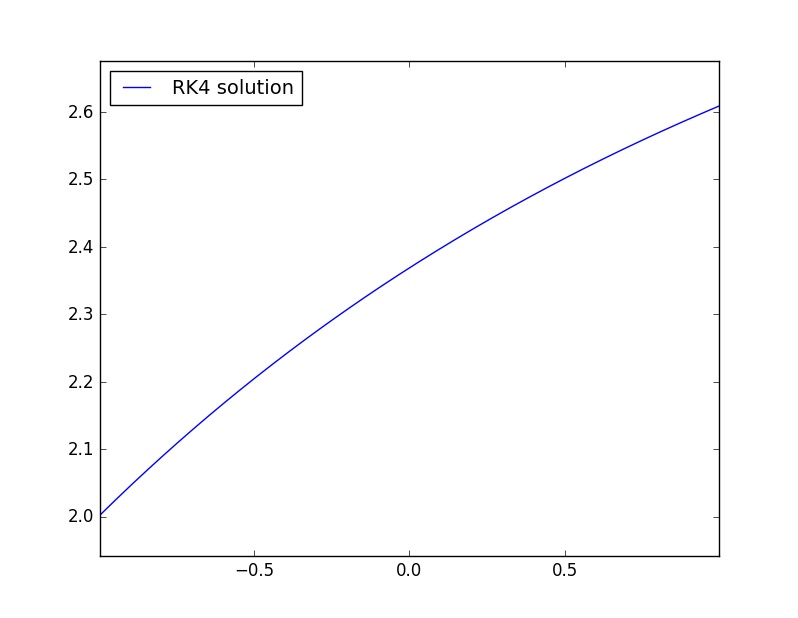}}
    \caption{Test 2.1. Initial condition: a stationary solution approximated with the RK4 method.}
    \label{Burgers2-ic}
 \end{center}
\end{figure}

We consider $x \in [-1, 1]$, $t \in [0,5]$ and $CFL=0.9$. The initial condition is the solution of the Cauchy problem
consisting  of \eqref{ODEBurgers2} with initial condition
$$
u(-1) = 2,
$$
which is a stationary solution of the problem. This solution is approximated using the RK4 method: see Figure \ref{Burgers2-ic}. $N_p=1$ is considered.

%$$\begin{cases}
%\displaystyle \frac{du}{dx} \left( x \right)= \frac{\sin(u(x))}{u(x)},\\
%u(-1)=2,
%\end{cases}$$

$u(-1,t)=2$ is imposed at $x=-1$ and free boundary conditions are considered at $x = 1$.

\begin{figure}[H]
\begin{center}
  \subfloat[SM$i$, $i=1,2,3$]{
   \includegraphics[width=0.5\textwidth]{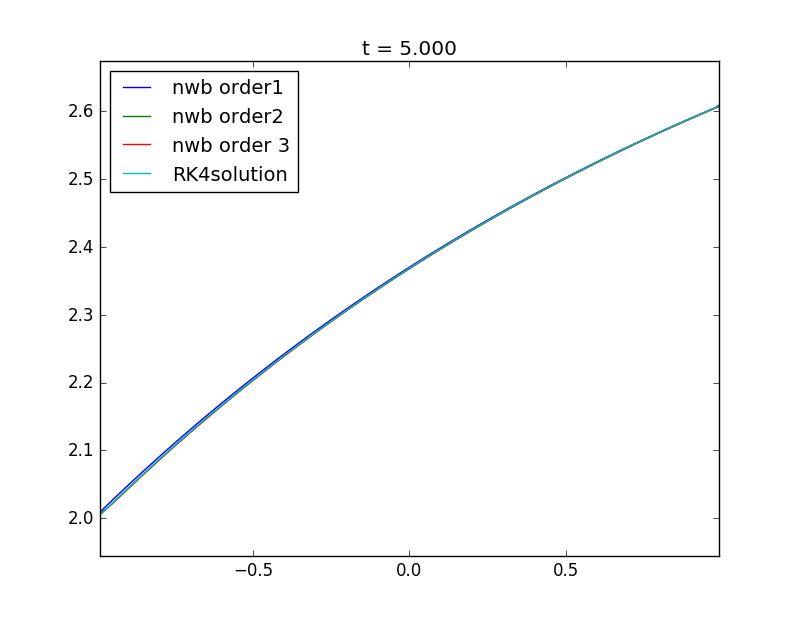}}
 \subfloat[DWBM$i$, $i=1,2,3$ ]{
   \includegraphics[width=0.5\textwidth]{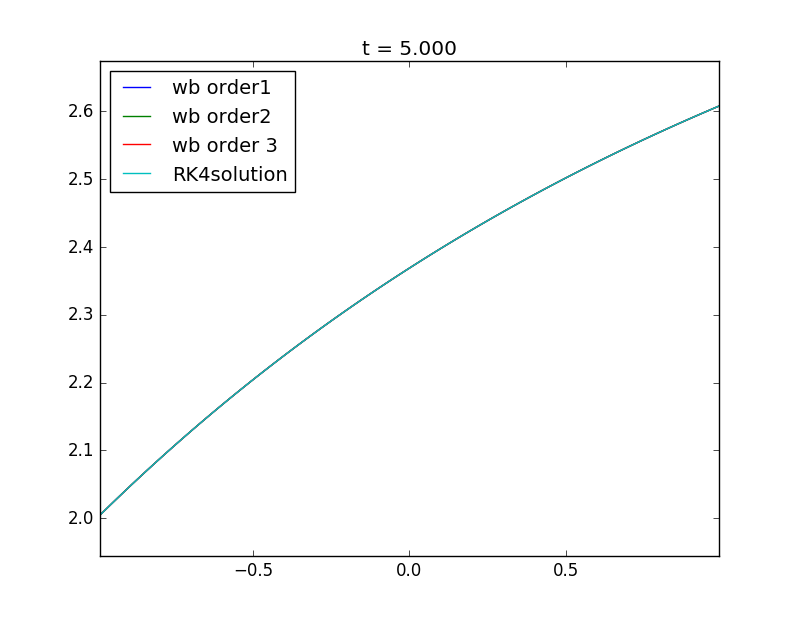}}
    \caption{Test 2.1. Numerical solutions at $t = 5s$. Number of cells: 100.} \label{test3_1}
 \end{center}
\end{figure}

\begin{figure}[H]
\begin{center}
  \subfloat{
   \includegraphics[width=0.5\textwidth]{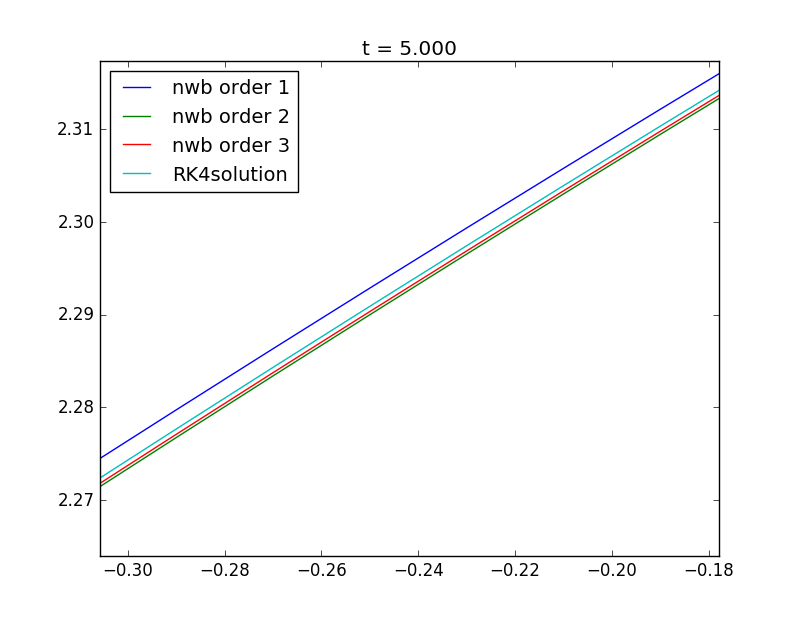}}
    \caption{Test 2.1. Zoom of the numerical solutions obtained with SM$i$, $i=1,2,3$ at $t = 5s$. Number of cells: 100.}
\label{zoomtest2.1}
 \end{center}
\end{figure}

Figure \ref{test3_1}  shows the numerical results obtained with SM$i$, $i=1,2,3$ (left) and DWBM$i$, $i=1,2,3$ (right). Notice that the non well-balanced methods perturb the stationary solution,  specially in a neighborhood of the left extreme: see the zoom in Figure \ref{zoomtest2.1}.

The maximum number of iterations required to solve the nonlinear problem \eqref{step1} applying Newton's method is two and it converges in only
one iteration for meshes with 200 cells or more. 
Tables \ref{ex3_error_nwb} and \ref{ex3_error_wb_newton} show the errors corresponding to SM$i$ and DWBM$i$, $i=1,2,3$ respectively. Computational costs are shown in  Table \ref{ex3_times}:
the well-balanced traitment increases the cost with a factor ranging from 5 to 15.6 in this case. 

\begin{table}[H]
\centering
\begin{tabular}{|c|cc|cc|cc|} \hline
Cells & Error ($i=1$)&Order&Error ($i=2$)& Order&Error ($i=3$)&Order\\\hline
100& 2.72E-3 & - & 1.43E-4 & - & 7.66E-5 & -\\
200 & 1.34E-3 & 1.021 & 2.43E-6 & 5.879 & 9.62E-8 & 10.506\\
400 & 6.58E-4 & 1.026 & 8.19E-7 & 1.569 & 1.21E-10 & 7.254\\
800 & 3.24E-4 & 1.022 & 2.47E-7 & 1.729 & 1.51E-11 & 2.922\\ \hline
\end{tabular}
\caption{Test 2.1. Errors in $L^1$ norm and convergence rates for SM$i$, $i=1,2,3$.} \label{ex3_error_nwb}

\end{table}

\begin{table}[H]
\centering
\begin{tabular}{|c|c|c|c|} \hline
Cells & Error ($i =1$)&Error ($i=2$)& Error ($i=3$)\\\hline
100& 1.71E-13 &  1.76E-13  &  2.54E-13\\
200 & 1.08E-14 &  3.33E-15 &  3.60E-14\\
400 & 1.98E-14 &  7.99E-16 &  2.12E-14\\ 
800 & 5.27E-15 &  9.19E-16 &  9.11E-14\\ \hline
\end{tabular}
\caption{Test 2.1. Errors in $L^1$ norm for DWBM$i$, $i=1,2,3$} \label{ex3_error_wb_newton}

\end{table}

%\begin{table}[H]
%\centering
%\begin{tabular}{|c|c|c|c|} \hline
%Cells & Error ($i=1$)&Error ($i=2$)& Error ($i=3$)\\\hline
%100& 5.19E-10 &  3.98E-13  &  5.64E-14\\
%200 & 3.15E-11 &  4.08E-15 &  2.60E-14\\
%400 & 8.92E-10 &  1.60E-14 &  5.47E-14\\ 
%800 & 1.09E-10 &  1.21E-14 &  8.02E-14\\ \hline
%\end{tabular}
%\caption{Errors in $L^1$ norm for DWBM$i$, $i=1,2,3$.} \label{ex3_error_wb}
%
%\end{table}

%\begin{table}[H]
%\centering
%\begin{tabular}{|c|c|c|c|c|}
%\hline
%          Cells        & $i$ & SM$i$ & DWBM$i$ & DWBM$i$\\
% \hline
%\multirow{3}{*}{100} & 1& 10 & 360 & 390 \\ %\cline{2-5} 
%                  & 2 & 20 & 950 & 1040 \\ %\cline{2-5} 
%                  & 3& 50 & 1590 & 2260 \\ \hline
%\multirow{3}{*}{200} & 1& 30 & 810 & 1420 \\ %\cline{2-5} 
%                  & 2& 60 & 2760 & 4000 \\ %\cline{2-5} 
%                  & 3& 190 & 6200 & 8880 \\ \hline
%%\multirow{3}{*}{400} & $1^{st}$ O & 100 & 3060 &  4310\\ %\cline{2-5} 
%%                  & $2^{nd}$ O & 210 & 10850 & 13350 \\ %\cline{2-5} 
%%                  & $3^{rd}$ O & 660 & 24510 & 30020 \\ \hline
%%\multirow{3}{*}{800} & $1^{st}$ O & 350 & 12060 & 16980 \\ %\cline{2-5} 
%%                  & $2^{nd}$ O & 740 & 43410 & 53230 \\ %\cline{2-5} 
%%                  & $3^{rd}$ O & 2470 & 97850 & 119790 \\ \hline
%\end{tabular}
%\caption{Computational times (milliseconds). $t=5s$.}\label{ex3_times}
%\end{table}

\begin{table}[H]
\centering
\begin{tabular}{|c|c|c|c|}
\hline
          Cells        & $i$ & SM$i$ & DWBM$i$ \\
 \hline
\multirow{3}{*}{100} & 1& 10 & 50 \\ %\cline{2-5} 
                  & 2 & 20 & 370 \\ %\cline{2-5} 
                  & 3& 50 &  760\\ \hline
\multirow{3}{*}{200} & 1& 30 & 200  \\ %\cline{2-5} 
                  & 2& 60 & 940 \\ %\cline{2-5} 
                  & 3& 190 & 2220 \\ \hline
%\multirow{3}{*}{400} & $1^{st}$ O & 100 & 3060 &  4310\\ %\cline{2-5} 
%                  & $2^{nd}$ O & 210 & 10850 & 13350 \\ %\cline{2-5} 
%                  & $3^{rd}$ O & 660 & 24510 & 30020 \\ \hline
%\multirow{3}{*}{800} & $1^{st}$ O & 350 & 12060 & 16980 \\ %\cline{2-5} 
%                  & $2^{nd}$ O & 740 & 43410 & 53230 \\ %\cline{2-5} 
%                  & $3^{rd}$ O & 2470 & 97850 & 119790 \\ \hline
\end{tabular}
\caption{Computational times (milliseconds). $t=5s$.}\label{ex3_times}
\end{table}

%\begin{table}[H]
%\centering
%\begin{tabular}{|c|c|c|c|c|c|c|}
%\hline
%\multicolumn{7}{|c|}{Newton's methods} \\ \cline{1-7} 
%\multirow{2}{*}{Cells} & \multicolumn{2}{c|}{$1^{st}$ Order} & \multicolumn{2}{c|}{$2^{nd}$ Order} & \multicolumn{2}{c|}{$3^{rd}$ Order} \\ \cline{2-7} 
%   &     Min      &     Max      &     Min      &     Max      &      Min     &     Max      \\ \hline
% 100        &     2      &     2      &      2     &     2      &      1     &      1     \\ %\hline
%200   &      1     &    1       &        1   &     1      &     1      &     1      \\ %\hline
% 400  &     1      &     1      &       1    &      1     &      1     &     1      \\ %\hline
%800  &      1     &       1    &     1      &     1      &       1    &     1      \\ \hline
%\multicolumn{7}{|c|}{Descent methods} \\ \cline{1-7} 
%\multirow{2}{*}{Cells} & \multicolumn{2}{c|}{$1^{st}$ Order} & \multicolumn{2}{c|}{$2^{nd}$ Order} & \multicolumn{2}{c|}{$3^{rd}$ Order} \\ \cline{2-7} 
%   &     Min      &     Max      &     Min      &     Max      &      Min     &     Max      \\ \hline
% 100        &     3      &     3      &      3    &     3      &      3     &      3     \\ %\hline
%200   &      3     &    3       &        3   &     3      &     3      &     3      \\ %\hline
% 400  &     2      &     2      &       2    &      2     &      2     &     2     \\ %\hline
%800  &      2     &       2    &     2      &     2      &       2    &     2      \\ \hline
%\end{tabular}
%\caption{Maximum and minimum number of iterations required to solve the nonlinear problem.}\label{ex3_it} 
%\end{table}

\subsubsection{Test 2.2}
\begin{figure}[H]
\begin{center}
  \subfloat{
   \includegraphics[width=0.5\textwidth]{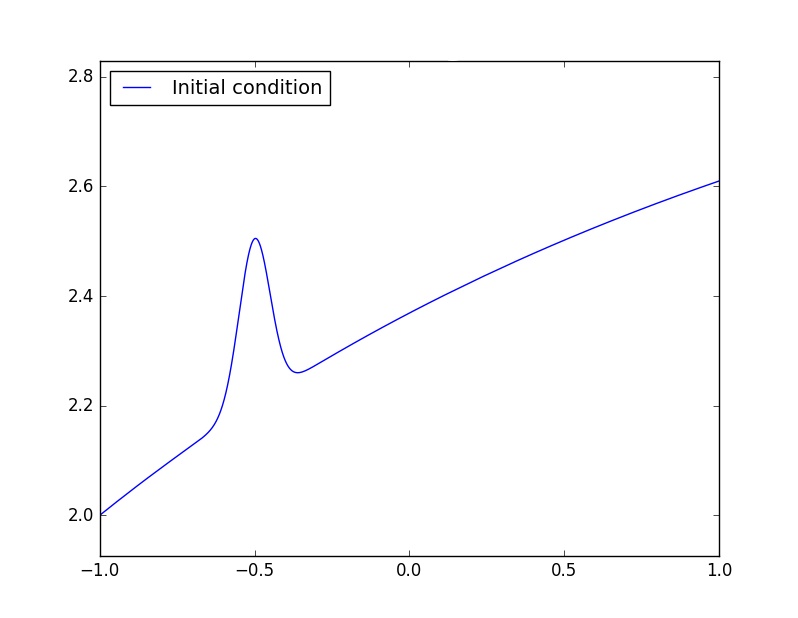}}
 \caption{Test 2.2. Initial condition.}
 \label{test22ic}
 \end{center}
\end{figure}

\begin{figure}[H]
\begin{center}
  \subfloat[SM$i$, $i =1,2,3$. $t = 0.3s$. ]{
   \includegraphics[width=0.5\textwidth]{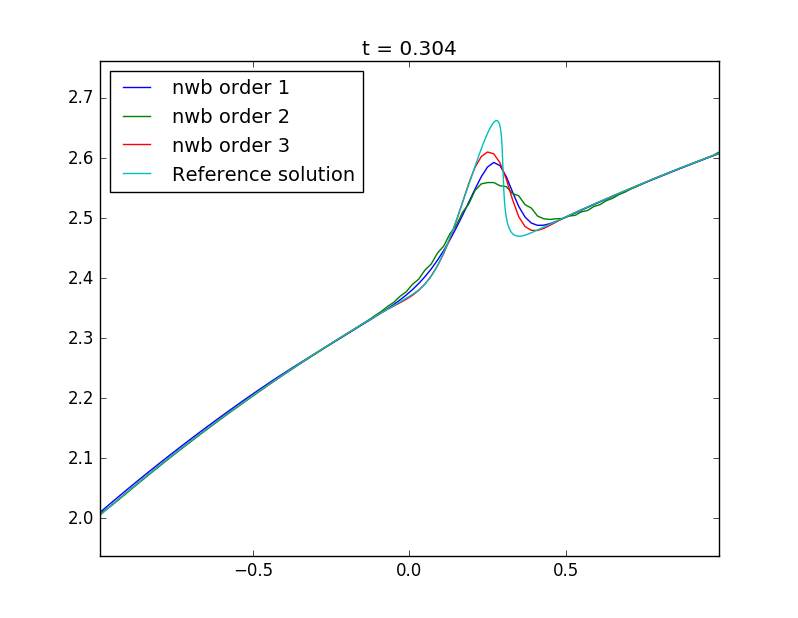}}
 \subfloat[DWBM$i$, $i=1,2,3$. $t = 0.3s$.]{
   \includegraphics[width=0.5\textwidth]{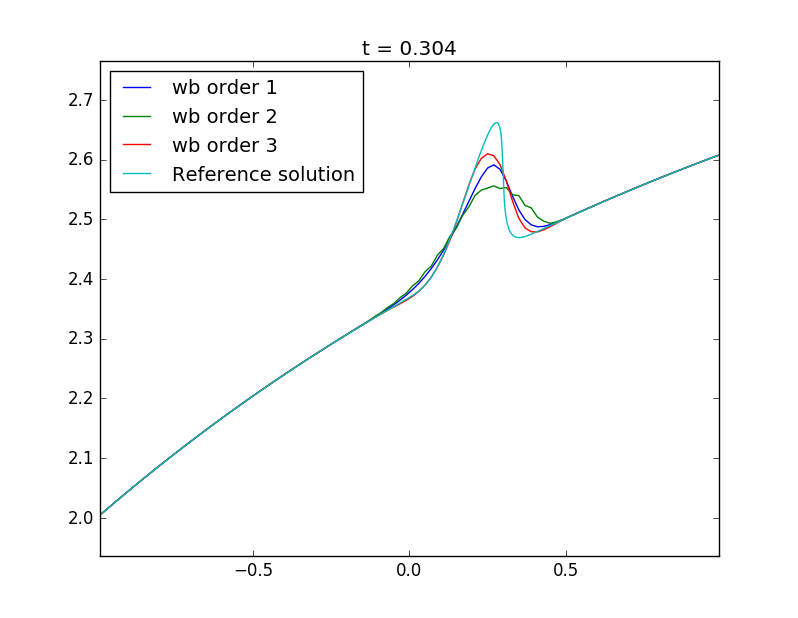}}\vspace{0.0001mm}
   \subfloat[SM$i$, $i =1,2,3$.. $t = 5s$.]{
   \includegraphics[width=0.5\textwidth]{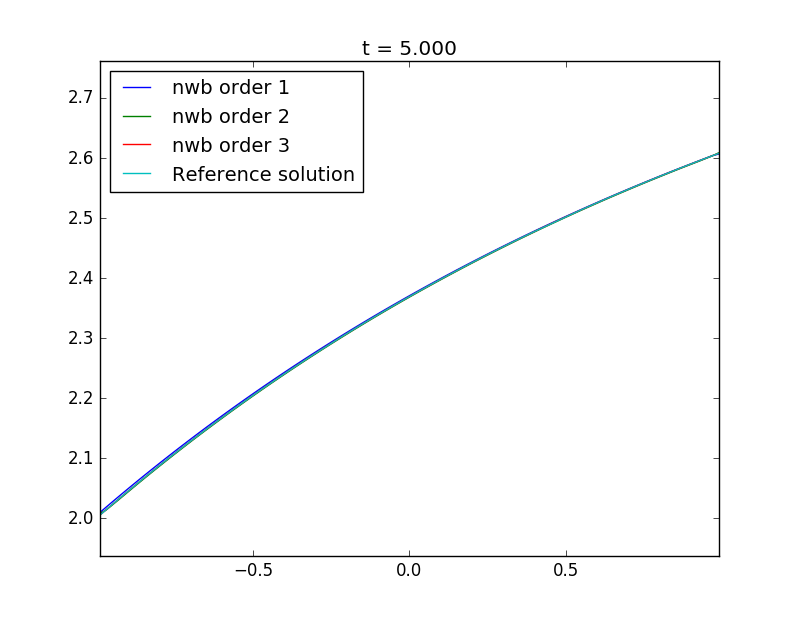}}
   \subfloat[DWBM$i$, $i=1,2,3$.. $t = 5s$.]{
   \includegraphics[width=0.5\textwidth]{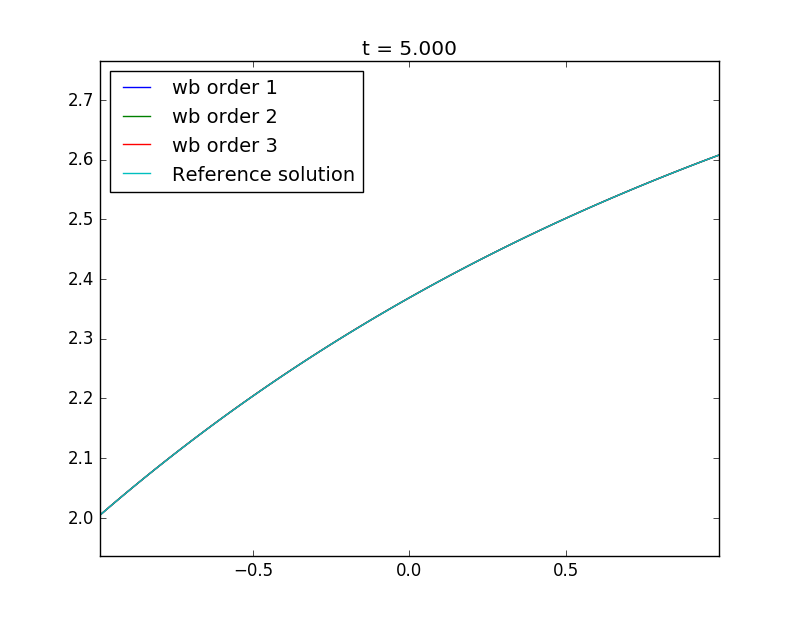}}
    \caption{Test 2.2. Reference and numerical solutions at time $t = 0.3, 5s$. Number of cells: 100.} \label{test4}
 \end{center}
\end{figure}

The evolution of a perturbation of the stationary solution considered in the previous test is now simulated. The only difference with Test 2.1. is that now the initial condition is given by:  
$$u_0(x)=u^*(x)+0.3e^{-200(x+0.5)^2},$$
where $u^*(x)$ is again the stationary solution satisfying $u^*(-1) = 2$: see Figure \ref{test22ic}.
Figure \ref{test4} shows the evolution of the perturbation at times $t=0.3, 5$ obtained with SM$i$, $i=1,2,3$ and DWBM$i$, $i=1,2,3$.
A reference solution has been computed with a first order well-balanced scheme on a fine mesh (12800 cells). 

Again, the main differences between the non-well balanced and the well-balanced methods are found once the perturbation has left the domain: only the well-balanced methods preserve the stationary solutions. This is very clear in Tables \ref{ex4_error_nwb} and \ref{ex4_error_wb_newton}, where the errors at time  $t = 5s$ are shown.  %Therefore, conclusions concerning to computational effort are similar to the ones included in table \ref{ex3_times}.

\begin{table}[H]
\centering
\begin{tabular}{|c|c|c|c|} \hline
Cells & Error ($i=1$)&Error ($i = 2$)& Error ($i = 3$)\\\hline
100& 3.43E-3 &  1.71E-3  &  1.06E-3\\
200 & 1.72E-3 &  8.48E-4 &  5.27E-4\\
400 & 8.59E-4 &  4.25E-4 &  2.61E-4\\ 
800 & 4.30E-4 &  2.11E-4 &  1.30E-4\\ \hline
\end{tabular}
\caption{Test 2.2. Errors in $L^1$ norm for SM$i$, $i=1,2,3$.} \label{ex4_error_nwb}

\end{table}

\begin{table}[H]
\centering
\begin{tabular}{|c|c|c|c|} \hline
Cells & Error ($i=1$)&Error ($i = 2$)& Error ($i = 3$)\\\hline
100& 3.15E-15 &  4.21E-15  &  2.46E-12\\
200 & 2.93E-15 &  1.30E-15 &  1.79E-13\\
400 & 3.81E-15 &  1.48E-15 &  5.77E-14\\ 
800 & 4.02E-15 &  2.42E-15 &  1.16E-13\\ \hline
\end{tabular}
\caption{Errors in $L^1$ norm for DWBM$i$, $i=1,2,3$.} \label{ex4_error_wb_newton}

\end{table}

%\begin{table}[H]
%\centering
%\begin{tabular}{|c|c|c|c|} \hline
%Cells & Error ($i=1$)&Error ($i = 2$)& Error ($i = 3$)\\\hline
%100& 5.32E-10 &  4.89E-12  &  2.84E-12\\
%200 & 3.37E-11 &  1.42E-13 &  5.50E-14\\
%400 & 9.83E-10 &  2.16E-12 &  1.23E-12\\ 
%800 & 1.23E-10 &  1.44E-13 &  6.04E-14\\ \hline
%\end{tabular}
%\caption{Errors in $L^1$ norm for DWBM$i$, $i=1,2,3$.} \label{ex4_error_wb}
%
%\end{table}
%

\subsection{Problem 3: coupled Burgers equations with nonlinear source terms}
Let us consider the system of balance laws
\begin{equation}\label{burgers_system}
\begin{cases}
\displaystyle \frac{\partial u_1}{\partial t} + \displaystyle \frac{\partial}{\partial x}\left( \frac{u_1^2}{2} \right)= 2 u_1^2 + u_1 u_2 ,\vspace{2mm} \\
\displaystyle \frac{\partial u_2}{\partial t} + \displaystyle \frac{\partial}{\partial x}\left( \frac{u_2^2}{2} \right)=  - u_1 u_2 + 3u_2^2,
\end{cases}
\end{equation}
which is the particular case of \eqref{sle} corresponding to the choices $N=2$,
$$ U =\begin{pmatrix}
u_1 \\
u_2 \\
\end{pmatrix} , \quad f(U) =\begin{pmatrix}
\displaystyle \frac{u_1^2}{2} \vspace{2mm}  \\
\displaystyle \frac{u_2^2}{2}\\
\end{pmatrix}, \quad S(U) =\begin{pmatrix}
 2 u_1^2 + u_1 u_2 \vspace{2mm} \\
- u_1 u_2 + 3u_2^2\\
\end{pmatrix}, \quad H(x)=x.$$
The system of ODE satisfied by the stationary solutions is the linear system:
\begin{equation} \label{state_equation_test5}
\left\{ 
\begin{array}{l}
\displaystyle \frac{du_1}{dx} = 2u_1 + u_2, \vspace{2mm} \\
\displaystyle  \frac{du_2}{dx} = - u_1  +3 u_2.  \\
\end{array}
\right.
\end{equation}
Therefore:
$$
G(x,U) = \left[ \begin{array}{c}  2u_1 + u_2 \\   - u_1  +3 u_2. \end{array} \right], \quad \nabla G(x, U) = \left[ \begin{array}{cc} 2 & 1 \\ -1 & 3 \end{array} \right]. 
$$

%The stationary solutions of the equation are given by the general solution of the ODE system: 
%\begin{equation}\label{ex5_stationary}
%\begin{cases}
%\begin{pmatrix} \displaystyle \frac{d u_1}{d x} \vspace{1mm}\\ \displaystyle  \frac{d u_2}{d x} \end{pmatrix}= \begin{pmatrix} 2 & 1 \\  -1 & 3 \end{pmatrix} \begin{pmatrix} u_1\\  u_2 \end{pmatrix},
%\end{cases}
%\end{equation}
The stationary solutions are given by the general solution of the ODE system:
\begin{equation} \label{test5_estac}
\left\{
\begin{array}{l}
u_1 (x) = \displaystyle c_1 e^{5x/2} \cos\left( \frac{\sqrt{3}}{2}x \right) +  c_2 e^{5x/2} \sin\left( \frac{\sqrt{3}}{2} x\right) , \\
u_2 (x) = \displaystyle \left( \frac{c_1}{2} + \frac{\sqrt{3}}{2}c_2 \right) e^{5x/2} \cos\left( \frac{\sqrt{3}}{2} x\right) +
 \left( - \frac{\sqrt{3}}{2} c_1 + \frac{c_2}{2} \right) e^{5x/2} \sin\left( \frac{\sqrt{3}}{2} x \right) , 
\end{array}
\right.
\end{equation}
%\begin{equation} \label{test5_estac}
%\begin{cases}
%u_1(x)=  \, \displaystyle \frac{2 c_1 e^{\frac{5x}{2}} \sin{ \left( \displaystyle \frac{\sqrt{3}}{2}x \right)} }{\sqrt{3}} + \frac{1}{3} c_2 e^{\frac{5x}{2}} \left( 3 \cos{ \left( \frac{\sqrt{3}}{2}x \right)} - \sqrt{3} \sin{\left( \frac{\sqrt{3}}{2}x \right)} \right),  \vspace{2mm} \\
%u_2(x)=  \, - \displaystyle \frac{2 c_2 e^{\frac{5x}{2}} \sin{ \left( \displaystyle \frac{\sqrt{3}}{2}x \right)} }{\sqrt{3}} + \frac{1}{3} c_1 e^{\frac{5x}{2}} \left( 3 \cos{ \left( \frac{\sqrt{3}}{2}x \right)} + \sqrt{3} \sin{\left( \frac{\sqrt{3}}{2}x \right)} \right).  \\
%\end{cases}
%\end{equation}

Since the expression of the stationary solutions is known, the first step of the well-balanced reconstruction procedure can be easily solved: given a family of cell values $ \left\lbrace U_i = \begin{pmatrix}
\displaystyle u_i^1 \vspace{2mm}  \\
u_i^2\\
\end{pmatrix} \right\rbrace$, the stationary solution $U^*_i=\begin{pmatrix}
\displaystyle u^*_1 \vspace{2mm}  \\
u^*_2\\
\end{pmatrix}$ which solves the non-linear problem \eqref{step1} is 
\begin{equation}
U^*_i(x)=\begin{pmatrix}
 \displaystyle a_i e^{5x/2} \cos\left( \frac{\sqrt{3}}{2}x \right) +  b_i e^{5x/2} \sin\left( \frac{\sqrt{3}}{2} x\right) \vspace{2mm}  \\
\displaystyle \left( \frac{a_i}{2} + \frac{\sqrt{3}}{2}b_i \right) e^{5x/2} \cos\left( \frac{\sqrt{3}}{2} x\right) +
 \left( - \frac{\sqrt{3}}{2} a_i + \frac{b_i}{2} \right) e^{5x/2} \sin\left( \frac{\sqrt{3}}{2} x \right)\\
\end{pmatrix},\label{step1cB}
\end{equation}
where
\begin{equation*}
\begin{split}
a_i & = \displaystyle \frac{\Delta x}{3} \displaystyle \frac{\sqrt{3} (4 u_i^1 - 5 u_i^2) \left[ e^{ \frac{5 \Delta x}{2}} \sin \left( \frac{\sqrt{3}}{2} x_{i+ \frac{1}{2}} \right) - \sin \left(  \frac{\sqrt{3}}{2} x_{i- \frac{1}{2}} \right) \right]}{e^{ \frac{5 }{2}x_{i- \frac{1}{2}}} \left[ e^{5 \Delta x} +1 - 2 e^{ \frac{5 \Delta x}{2}} \cos \left(  \frac{\sqrt{3}}{2} \Delta x \right) \right] } \\
& + \displaystyle \frac{\Delta x}{3} \displaystyle \frac{ 3 (2 u_i^1 + u_i^2) \left[ e^{ \frac{5 \Delta x}{2}} \cos \left(  \frac{\sqrt{3}}{2} x_{i+ \frac{1}{2}} \right) - \cos \left(  \frac{\sqrt{3}}{2} x_{i- \frac{1}{2}} \right) \right]}{e^{ \frac{5 }{2}x_{i- \frac{1}{2}}} \left[ e^{5 \Delta x} +1 - 2 e^{ \frac{5 \Delta x}{2}} \cos \left(  \frac{\sqrt{3}}{2} \Delta x \right) \right] }, 
\end{split}
\end{equation*}
\begin{equation*}
\begin{split}
b_i & = \displaystyle \frac{\Delta x}{3} \displaystyle \frac{3 (2 u_i^1 + u_i^2) \left[ e^{ \frac{5 \Delta x}{2}} \sin \left( \frac{\sqrt{3}}{2} x_{i+ \frac{1}{2}} \right) - \sin \left(  \frac{\sqrt{3}}{2} x_{i- \frac{1}{2}} \right) \right]}{e^{ \frac{5 }{2}x_{i- \frac{1}{2}}} \left[ e^{5 \Delta x} +1 - 2 e^{ \frac{5 \Delta x}{2}} \cos \left(  \frac{\sqrt{3}}{2} \Delta x \right) \right] } \\
& + \displaystyle \frac{\Delta x}{3} \displaystyle \frac{ \sqrt{3} (5 u_i^1 - 4 u_i^2) \left[ e^{ \frac{5 \Delta x}{2}} \cos \left(  \frac{\sqrt{3}}{2} x_{i+ \frac{1}{2}} \right) - \cos \left(  \frac{\sqrt{3}}{2} x_{i- \frac{1}{2}} \right) \right]}{e^{ \frac{5 }{2}x_{i- \frac{1}{2}}} \left[ e^{5 \Delta x} +1 - 2 e^{ \frac{5 \Delta x}{2}} \cos \left(  \frac{\sqrt{3}}{2} \Delta x \right) \right] }.
\end{split}
\end{equation*}

Therefore the well-balanced reconstruction can be easily implemented using this explicit expression. The techniques described in the previous sections to compute numerically the solution of 
\eqref{step1} will be used to measure their efficiencies and its sensitivity to the numerical discretization of the ODE \eqref{state_equation_test5}. 

\subsubsection{Test 3.1}

\begin{figure}[H]
\begin{center}
  \subfloat{
   \includegraphics[width=0.5\textwidth]{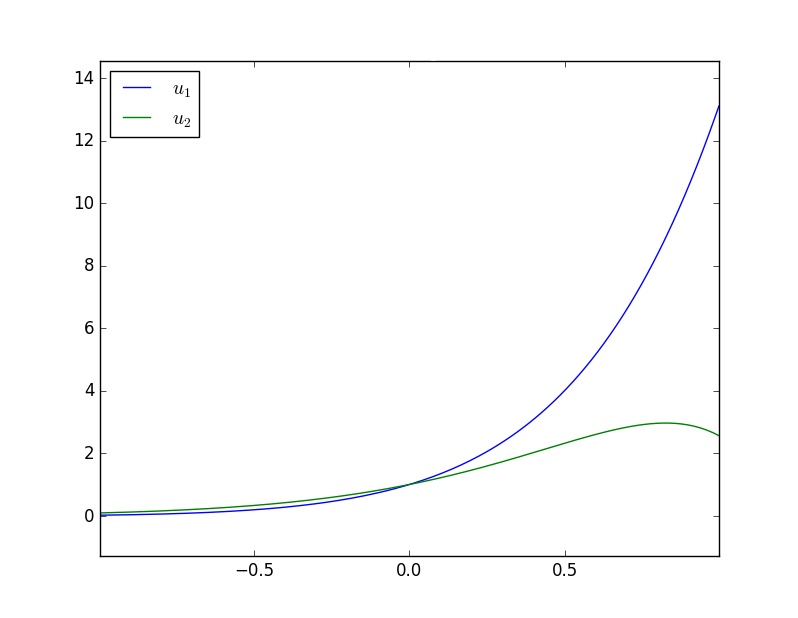}}
    \caption{Initial condition for Test 3.1.} \label{test5_cini}
 \end{center}
\end{figure}

We consider $x \in [-1, 1]$, $t \in [0,5]$, and $CFL=0.9$. The initial condition, shown in Figure \ref{test5_cini}, is the stationary solution
\begin{equation} \label{test5_ic}
\left\{
\begin{array}{l}
u^*_1 (x) = \displaystyle e^{5x/2} \cos\left( \frac{\sqrt{3}}{2}x \right) +  \frac{ \sqrt{3}}{3}  e^{5x/2} \sin\left( \frac{\sqrt{3}}{2} x\right) , \\
u^*_2 (x) = \displaystyle e^{5x/2} \cos\left( \frac{\sqrt{3}}{2} x\right)  - \frac{\sqrt{3}}{3}  e^{5x/2} \sin\left( \frac{\sqrt{3}}{2} x \right) .  
\end{array}
\right.
\end{equation}
In this problem, since \eqref{state_equation_test5} is linear, Newton's method converges again in only one iteration. 

\begin{figure}[H]
\begin{center}
  \subfloat[Errors (logaritmic scale) ]{
   \includegraphics[width=0.5\textwidth]{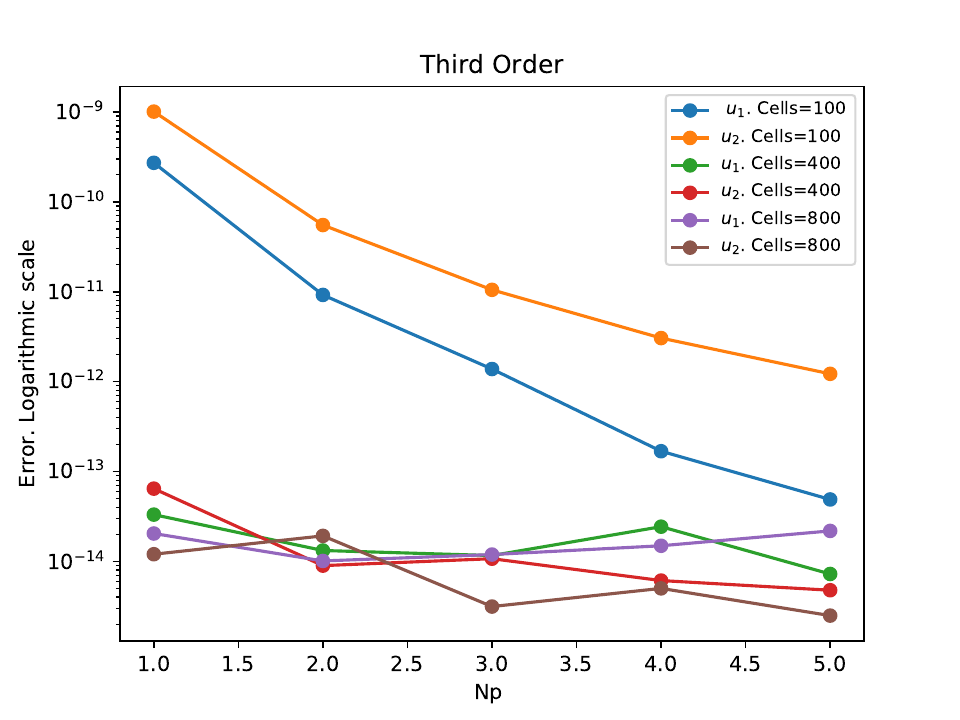}}
 \subfloat[CPU time]{
   \includegraphics[width=0.5\textwidth]{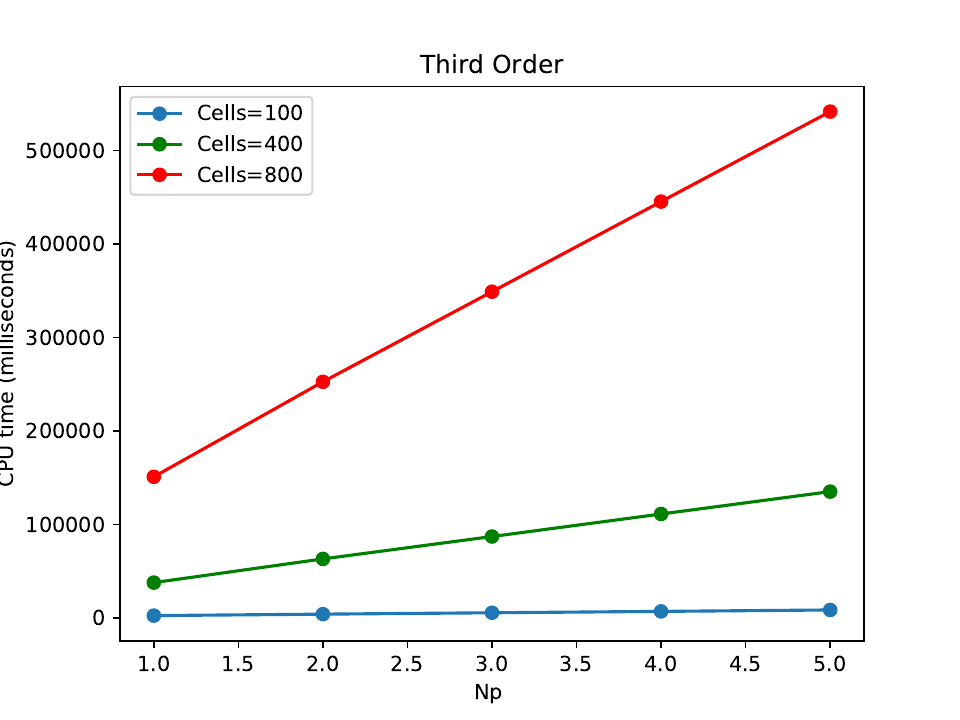}}
    \caption{Test 3.1. Errors and CPU times corresponding to DWBM3 with different number of cells and different values of $N_p$.} 
    \label{np_test3.1}
 \end{center}
\end{figure}

Boundary conditions
$$
u_j(-1, t) = u^*_j(-1), \quad j=1,2
$$
are imposed at the left extreme of the interval and free boundary conditions 
 at $x = 1$. Different values for $N_p$ have been compared. Figure \ref{np_test3.1} shows the errors at logarithmic scale and the CPU times corresponding to different values of $N_p$ for the third order method. All the errors are below $10^{-13}$ for fine enough meshes with $N_p = 1$. For $N_p = 3$ the errors for all the meshes are below that threshold.
 The results and conclusions are similar for the first and second order methods. 

\begin{figure}[H]
\begin{center}
  \subfloat[SM$i$, $i=1,2,3$]{
   \includegraphics[width=0.5\textwidth]{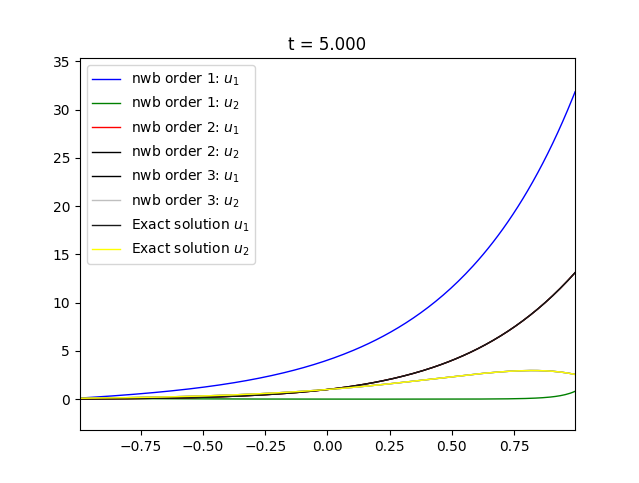}}
 \subfloat[WBM$i$, $i=1,2,3$]{
   \includegraphics[width=0.5\textwidth]{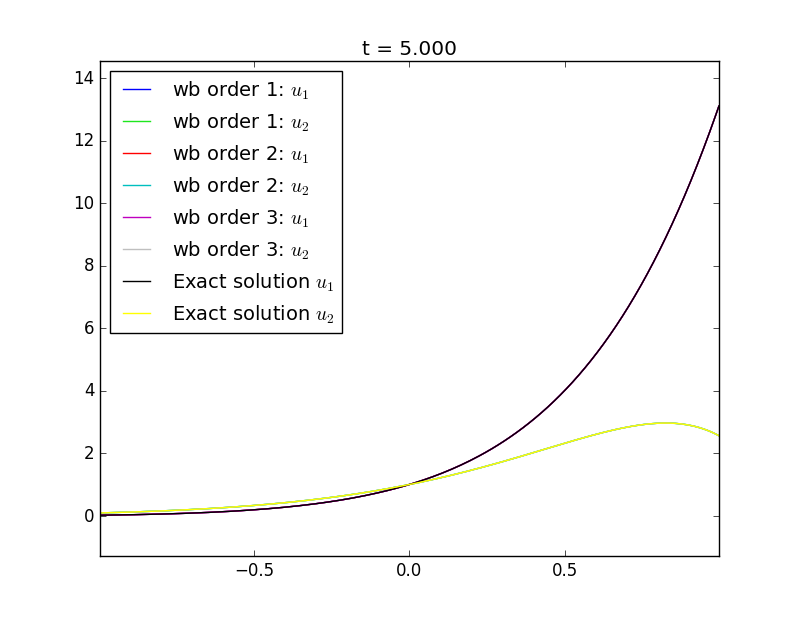}}
    \caption{Test 3.1. Exact and numerical solutions at $t = 5s$. Number of cells: 200.} \label{test5_1}
 \end{center}
\end{figure}

Figure \ref{test5_1} shows the numerical solutions at time $t= 5s$ with SM$i$, $i=1,2,3$ and WBM$i$, $i=1,2,3$ and
Tables \ref{ex5_error_nwb}, \ref{ex5_error_wb_exacta} and \ref{ex5_error_wb_newton} show the errors. The conclusions are similar to the previous test cases, although in this case the difference between the exact solution and  the numerical solutions obtained with the non-well balanced in a coarse mesh is much bigger.

\begin{table}[H]
\centering
\begin{tabular}{|c|cc|cc|cc|} \hline
Cells & Error ($i = 1$)  &Order&Error ($i = 2$)  & Order&Error ($i=3$) &Order\\
 & \multicolumn{2}{|c|}{$u_1$}  & \multicolumn{2}{|c|}{$u_1$} & \multicolumn{2}{|c|}{$u_1$} \\\hline
100& 7.81 & - & 1.06E-1  & - & 4.43E-3 & -\\
200 & 3.08 & 1.342 & 2.69E-2 & 1.978 & 5.59E-4 & 2.986\\
400 & 1.50 & 1.038 & 7.01E-3 & 1.940 & 7.02E-5 & 2.978\\
800 & 7.42E-1 & 1.015 & 1.88E-3 & 1.899 & 8.79E-6 & 2.997\\  \hline
Cells & Error ($i = 1$)  &Order&Error ($i = 2$)  & Order&Error ($i=3$) &Order\\
 & \multicolumn{2}{|c|}{$u_2$}  & \multicolumn{2}{|c|}{$u_2$} & \multicolumn{2}{|c|}{$u_2$} \\\hline
100& 1.98 & - & 9.12E-2  & - & 3.57E-3 & -\\
200 & 1.33 & 0.574 & 2.28E-2 & 2.000 & 4.53E-4 & 2.978\\
400 & 7.08E-1 & 0.910 & 5.79E-3 & 1.977 & 5.69E-5 & 2.993\\
800 & 3.31E-1 & 1.0971 & 1.48E-3 & 1.968 & 7.13E-6 & 2.997\\  \hline
\end{tabular}
\caption{Test 3.1. Errors in $L^1$ norm and convergence rates for SM$i$, $i=1,2,3$.} \label{ex5_error_nwb}

\end{table}

\begin{table}[H]
\centering
\begin{tabular}{|c|cc|cc|cc|} \hline
Cells & Error ($i = 1$)  &&Error ($i = 2$)  & &Error ($i=3$) &\\
  & $u_1$&$u_2$ &$u_1$& $u_2$ &$u_1$&$u_2$\\\hline
100& 2.01E-14 & 1.93E-14 & 2.15E-15  & 5.69E-15 & 1.80E-13 & 9.92E-14\\
200 & 1.29E-14 & 2.22E-14 & 1.80E-15 & 5.43E-15 & 4.98E-13 & 2.71E-13\\
400 & 5.28E-14 & 6.64E-14 & 3.11E-14 & 3.11E-14 & 6.96E-13 & 4.75E-13\\
800 & 8.29E-14 & 7.65E-14 & 1.65E-14 & 1.61E-15 & 1.21E-13 & 8.58E-13\\ \hline
\end{tabular}
\caption{Errors in $L^1$ norm for the WBM$i$, $i=1,2,3$.} \label{ex5_error_wb_exacta}

\end{table}

\begin{table}[H]
\centering
\begin{tabular}{|c|cc|cc|cc|} \hline
Cells & Error ($i = 1$)  &&Error ($i = 2$)  & &Error ($i=3$) &\\
  & $u_1$&$u_2$ &$u_1$& $u_2$ &$u_1$&$u_2$\\\hline
100& 7.03E-10 & 9.57E-9 & 6.14E-10  & 1.40E-9 & 2.72E-10 & 1.01E-9\\
200 & 2.10E-11 & 2.13E-10 & 2.16E-11 & 3.96E-11 & 5.38E-12 & 1.13E-11\\
400 & 6.89E-13 & 4.41E-12 & 2.24E-13 & 4.31E-13 & 3.31E-14 & 6.44E-14\\
800 & 1.53E-14 & 6.29E-14 & 2.16E-15 & 8.46E-16 & 2.04E-14 & 1.20E-14\\ \hline
\end{tabular}
\caption{Errors in $L^1$ norm for  DWBM$i$, $i=1,2,3$.} \label{ex5_error_wb_newton}

\end{table}

Again, the behavior of the errors is as expected for both the non-well-balanced and the well-balanced methods. Computational costs are shown in Table \ref{ex5_times}: 
observe that, in this case, DWBM$i$  are less costly than WBM$i$, due to the large number of operations required to compute the exact solution of \eqref{step1}: see \eqref{step1cB}. The unexpected computational cost corresponding to the first order method using the mesh of 100 cells is due to the fact that the numerical solution is very far from the stationary solution and takes very large values, what implies an important reduction of the time step.

\begin{table}[H]
\centering
\begin{tabular}{|c|c|c|c|c|}
\hline
          Cells        & $i$ & SM$i$ & WBMS$i$ & DWBM$i$ \\\hline
\multirow{3}{*}{100} & 1& 1720 & 340 & 180 \\ %\cline{2-5} 
                  & 2& 180 & 920 & 850 \\ %\cline{2-5} 
                  & 3 & 460 & 2670 & 2270 \\ \hline
\multirow{3}{*}{200} & 1& 470 & 1230 & 610 \\ %\cline{2-5} 
                  & 2  & 640 & 3560 & 3170 \\ %\cline{2-5} 
                  & 3 & 1650 & 10790 & 9450 \\ \hline
%\multirow{3}{*}{400} & $1^{st}$ O & 2350 & 9480 &  29550\\ %\cline{2-5} 
%                  & $2^{nd}$ O & 5330 & 27190 & 83860 \\ %\cline{2-5} 
%                  & $3^{rd}$ O & 14550 & 55540 & 161940 \\ \hline
%\multirow{3}{*}{800} & $1^{st}$ O & 9100 & 37050 & 108850 \\ %\cline{2-5} 
%                  & $2^{nd}$ O & 21090 & 109230 & 316670 \\ %\cline{2-5} 
%                  & $3^{rd}$ O & 58190 & 208980 & 648630 \\ \hline
\end{tabular}
\caption{Test 3.1. Computational times (milliseconds). $t=5s$.}\label{ex5_times}
\end{table}

\subsection{Problem 4: shallow water equations}
Let us consider the shallow water model, which is the particular case of \eqref{sle} corresponding to the choices $N=2$,
$$ U =\begin{pmatrix}
h \\
q \\
\end{pmatrix} , \quad f(U) =\begin{pmatrix}
q \\
\displaystyle \frac{q^2}{h}+\displaystyle \frac{g}{2}h^2\\
\end{pmatrix}, \quad S(U) =\begin{pmatrix}
0\\
gh\\
\end{pmatrix}.$$
The variable $x$ makes reference to the axis of the channel and $t$ is the time; $q(x,t)$ and $h(x,t)$ are the discharge and the thickness, respectively; $g$ is the gravity and $H(x)$ is the depth function measured from a fixed reference level. We denote by $u=q/h$ the depth-averaged velocity and $c=\sqrt{gh}$. 

The eigenvalues of the Jacobian matrix $D_f(U)$ of the flux function $f(U)$ are the following:
$$r_1=u-\sqrt{c}, \quad r_2=u+\sqrt{c}.$$
The Froude number, given by 
\begin{equation}
Fr(U)= \displaystyle \frac{|u|}{c},
\end{equation}
indicates the flow regime: subcritical ($Fr<1$), critical ($Fr=1$) or supercritical ($Fr>1$). 

The system of ODE satisfied by the stationary solutions is:
\begin{equation}
\begin{cases}
q_x=0,\\
\left( \displaystyle \frac{q^2}{h} + \frac{1}{2} g h^2 \right)_x= ghH_x.\\
\end{cases}
\end{equation}
It can be easily checked that, while $Fr(U) \not= 1$, this system can be written as follows: 
%\begin{equation}\label{u'=Gsw}
%\begin{cases}
%h_x= \displaystyle \frac{gh^3 H_x}{- {q^2}+ gh^3 },\\
%q_x=0,
%\end{cases}
%\end{equation}
\begin{equation}\label{u'=Gsw}
\begin{cases}
h_x= \displaystyle \frac{ghH_x}{- {u^2}+ gh } ,\\
q_x=0,
\end{cases}
\end{equation}
that is
\begin{equation}\label{Gsw}
G(U,x) = \left[ \begin{array}{c}
\displaystyle
\frac{ghH_x}{- {u^2}+ gh}\\
\\
0 \\
\end{array}
\right]
\end{equation}
%\begin{equation}\label{Gsw}
%G(U,x) = \left[ \begin{array}{c}
%\displaystyle
%\frac{gh^3 H_x}{- {q^2}+ gh^3 }\\
%\\
%0 \\
%\end{array}
%\right]
%\end{equation}
and thus
\begin{equation}\label{nablaGsw}
\nabla_U G  = \left[ \begin{array}{cc}
\displaystyle
-\frac{3gu^2 H_x}{(- {u^2}+ gh)^2 } & \displaystyle \frac{2gu H_x}{(- {u^2}+ gh)^2 } \\
& \\
0  & 0\\
\end{array}
\right].
\end{equation}
The stationary solutions are given in implicit form by:
\begin{equation}\label{incondimpl}
q = C_1, \quad \frac{q^2}{2 h^2} + gh -g H = C_2, \quad C_1, C_2 \in \mathbb{R},
\end{equation}
In \cite{lopez2013} a family of high-order well-balanced methods numerical methods was presented in which \eqref{step1} was solved on the basis of
this implicit form. 
%\begin{equation}\label{nablaGsw}
%\nabla_u G  = \left[ \begin{array}{cc}
%\displaystyle
%-\frac{3gH_xh^2q^2}{(- {q^2}+ gh^3)^2 } & \displaystyle \frac{2gH_xh^3q}{(- {q^2}+ gh^3)^2 } \\
%& \\
%0  & 0\\
%\end{array}
%\right].
%\end{equation}
%

In this case, the expression of Newton's method is particularly simple. In effect, notice first that the equation for $\lambda_2 = [\lambda_{2,1}, \lambda_{2,2}]^T$ is
(see \eqref{adjoint_system}):
$$
\left\{
\begin{array}{l}
\displaystyle \frac{d \lambda_{2,1}}{dx} = \frac{3gu^2 H_x}{(- {u^2}+ gh)^2 }\lambda_{2,1},\smallskip\\
\displaystyle \frac{d \lambda_{2,2}}{dx} = - 1  -\frac{2gu H_x}{(- {u^2}+ gh)^2 }\lambda_{2,1} ,\smallskip\\
\lambda_{2,1}(\Delta x) = \lambda_{2,2}(\Delta x)  = 0,
\end{array}
\right.
$$
whose solution is
$$
\lambda_{2,1}(x) = 0, \quad \lambda_{2,2}(x) =  \Delta x - x.
$$

As a consequence:
$$
\Lambda(0)^T  = \left[ \begin{array}{cc}
\lambda_{1,1}(0) & \lambda_{1,2}(0) \\
0  & \Delta x
\end{array}
\right].
$$
Let us suppose that Newton's method is used to solve \eqref{problem}, with $W = [ \bar h, \bar q]^T$ and $U^0_0 = W$. Then, to compute $\mathcal{F}(U^0_0)$ (where $\mathcal{F}$ is given by \eqref{F}), first the solution $U(x,W ) = [h(x, W), q(x, W)]^T$ of  system \eqref{u'=Gsw} has to be solved with initial condition:
$$
 h(0) = \bar h, \quad q(0) = \bar q, 
$$
Clearly, the solution for $q$ is
$$
q(x, W) = \bar q, \quad \forall x,
$$
and thus
$$
\mathcal{F}(U^0_0) = \left[ \begin{array}{c} \displaystyle  \frac{1}{\Delta x} \int_0^{\Delta x} h(x, W) \, dx  \\  \bar q \end{array} \right].
$$
In order to update $U$, the following linear system has then to be solved:
$$
\left[ \begin{array}{cc}
\lambda_{1,1}(0) & \lambda_{1,2}(0) \\
0  & \Delta x
\end{array}
\right]. \left[ \begin{array}{c} v_{0,1} \\ v_{0,2} \end{array} \right] =  \left[ \begin{array}{c} \displaystyle  \int_0^{\Delta x} h(x, W) \, dx  - \Delta x \bar h \\  0 \end{array} \right],
$$
whose solution is
$$
v_{0,1} =  \frac{\Delta x}{\lambda_{1,1}(0)}\left( \frac{1}{\Delta x} \int_0^{\Delta x} h_0(x) \, dx - \bar h\right)  , \quad v_{0,2} = 0,
$$
and then
$$
U^1_0  = \left[
\begin{array}{c}
h^1_0 
\smallskip\\
q^1_0
\end{array}
\right] =
\left[
\begin{array}{c}
\displaystyle h^0_0 -    \frac{\Delta x}{\lambda_{1,1}(0)}\left( \frac{1}{\Delta x} \int_0^{\Delta x} h_0(x) \, dx - \bar h\right)    \\ \bar q
\end{array}
\right].
$$

Reasoning by induction, it can be easily checked that Newton's method writes in this case as follows:
\begin{alg}{Newton's method}
  \begin{itemize}
\item $h^0_0 = \bar h$;

\item For k = 0,1,2\dots

\begin{itemize}

\item Compute the solution $h_k$ of 
$$
\left\{
\begin{array}{l}
\displaystyle \frac{d h}{dx}= \frac{gh^3 H_x}{- {\bar q^2}+ gh^3 }, \smallskip\\
h(0) = h^k_0
\end{array}
\right.
$$
 in the interval $[0, \Delta x]$.

\item Compute the solution $\lambda_k$ of
$$
\left\{
\begin{array}{l}
\displaystyle \frac{d \lambda}{dx}= -1  + \frac{3gH_xh_k^2\bar q^2}{(- {\bar q^2}+ gh_k^3)^2 }\lambda, \smallskip\\
\lambda(\Delta x) = 0
\end{array}
\right.
$$
 in the interval $[0, \Delta x]$.

\item Update $h^k_0$:
\begin{equation}\label{NewtonSW}
h^{k+1}_0 = h^k_0 -  \frac{\Delta x}{\lambda_k(0)}\left( \frac{1}{\Delta x} \int_0^{\Delta x} h_k(x) \, dx - \bar h\right).
\end{equation}

\end{itemize}

\end{itemize}

\end{alg}

%
%\begin{alg}{Gradient method}
%  \begin{itemize}
%\item $h^0_0 = \bar h$;
%
%\item For k = 0,1,2\dots
%
%\begin{itemize}
%
%\item Compute the solution $h_k$ of 
%$$
%\left\{
%\begin{array}{l}
%\displaystyle \frac{d h}{dx}= \frac{gh^3 H_x}{- {\bar q^2}+ gh^3 }, \smallskip\\
%h(0) = h^k_0
%\end{array}
%\right.
%$$
% in the interval $[0, \Delta x]$.
%
%\item Compute the solution $\lambda_k$ of
%$$
%\left\{
%\begin{array}{l}
%\displaystyle \frac{d \lambda}{dx}= -1  + \frac{3gH_xh_k^2\bar q^2}{(- {\bar q^2}+ gh_k^3)^2 }\lambda, \smallskip\\
%\lambda(\Delta x) = 0
%\end{array}
%\right.
%$$
% in the interval $[0, \Delta x]$.
%
%\item Update $h^k_0$:
%$$
%h^{k+1}_0 = h^k_0 -  \frac{\rho_k}{\Delta x} 2 \lambda_k(0) \left( \frac{1}{\Delta x} \int_0^{\Delta x} h_k(x) \, dx - \bar h\right),
%$$
%where $\rho_k$ is the descent step.
%\end{itemize}
%
%\end{itemize}
%
%\end{alg}

\subsubsection{Test 4.1}

\begin{figure}[H]
\begin{center}
  \subfloat[Initial condition. Free surface and bottom. ]{
   \includegraphics[width=0.5\textwidth]{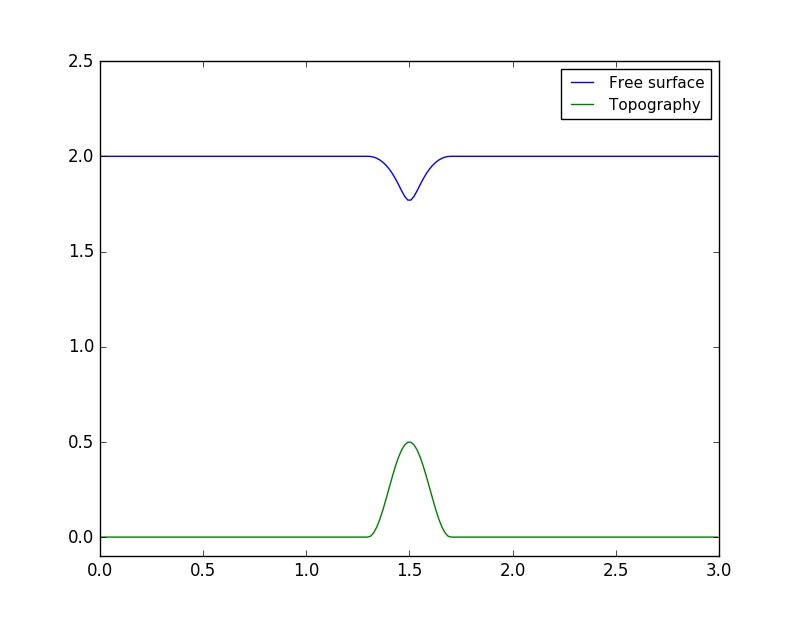}}
 \subfloat[Initial condition. Velocity.]{
   \includegraphics[width=0.5\textwidth]{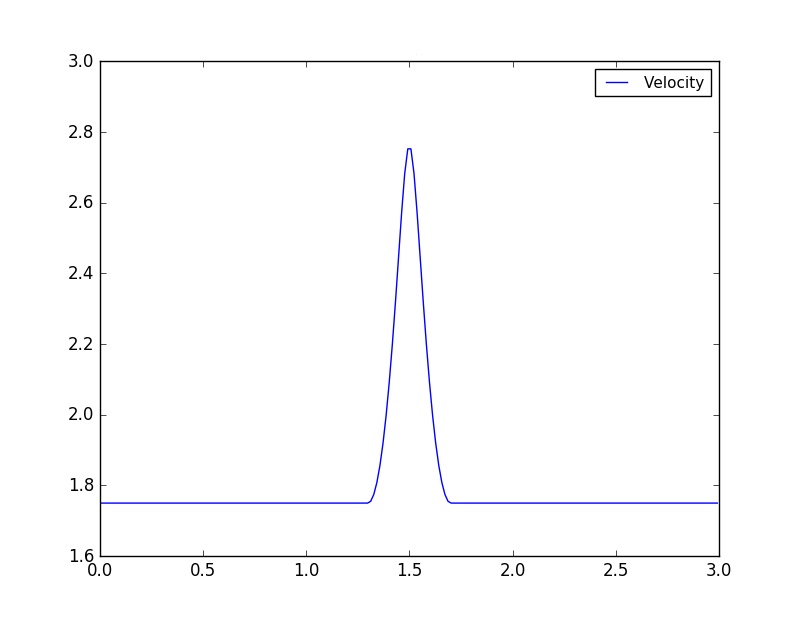}}
    \caption{Test 4.1. Initial condition: a subcritical stationary solution computed with RK4.} \label{test6_1}
 \end{center}
\end{figure}

\begin{figure}[H]
\begin{center}
  \subfloat[SM$i$, $i=1,2,3$. Free surface and bottom. ]{
   \includegraphics[width=0.5\textwidth]{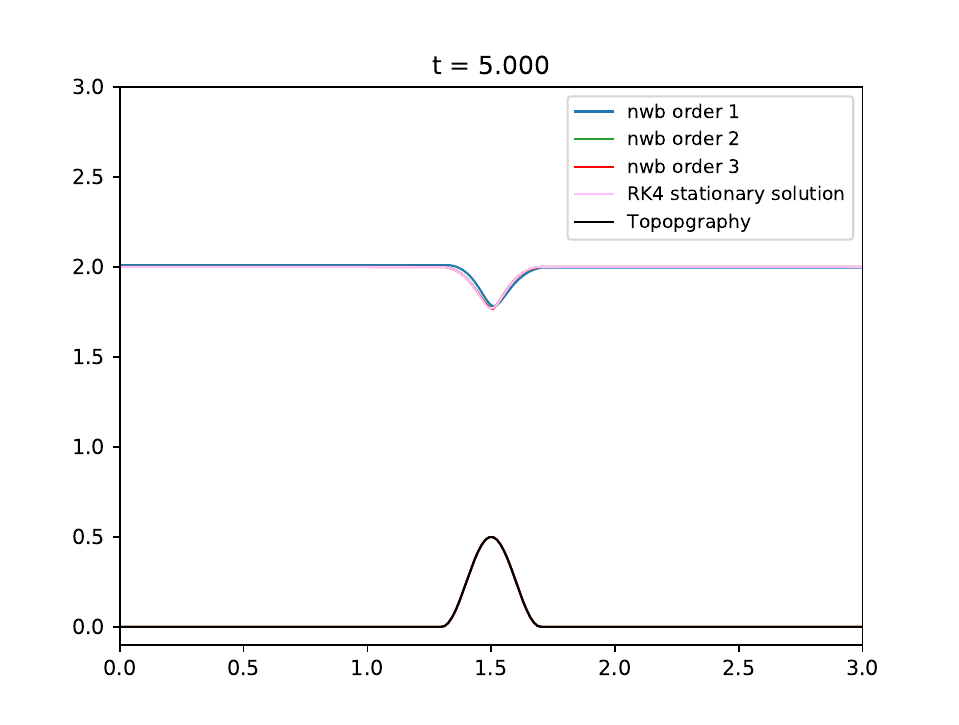}}
 \subfloat[SM$i$, $i=1,2,3$.  Velocity.]{
   \includegraphics[width=0.5\textwidth]{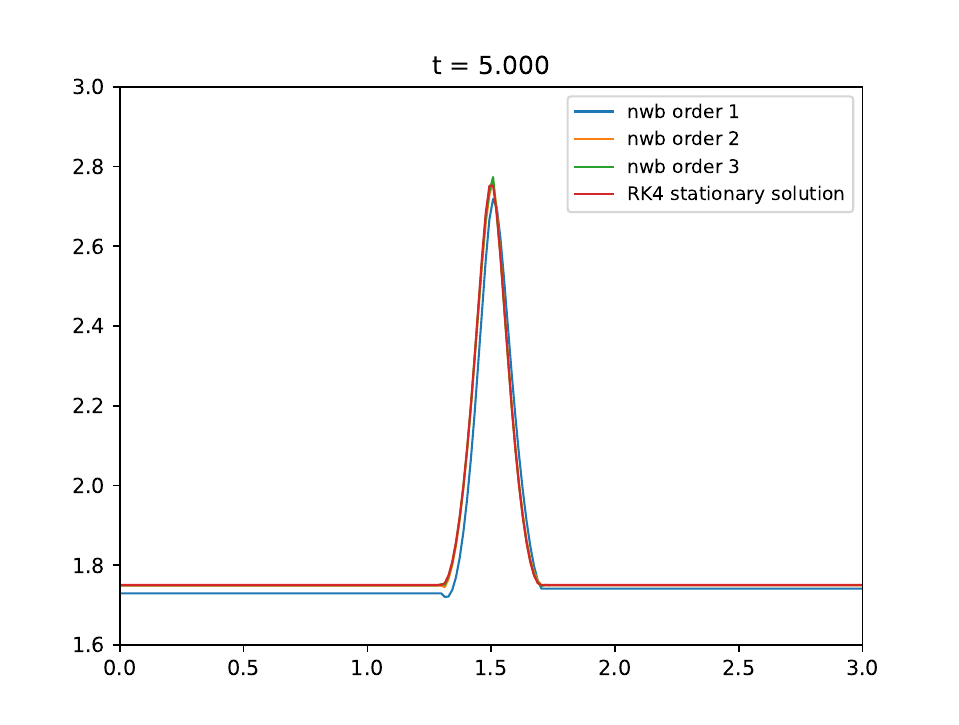}}\vspace{0.0001mm}
   \subfloat[DWBM$i$, $i=1,2,3$. Free surface and bottom.]{
   \includegraphics[width=0.5\textwidth]{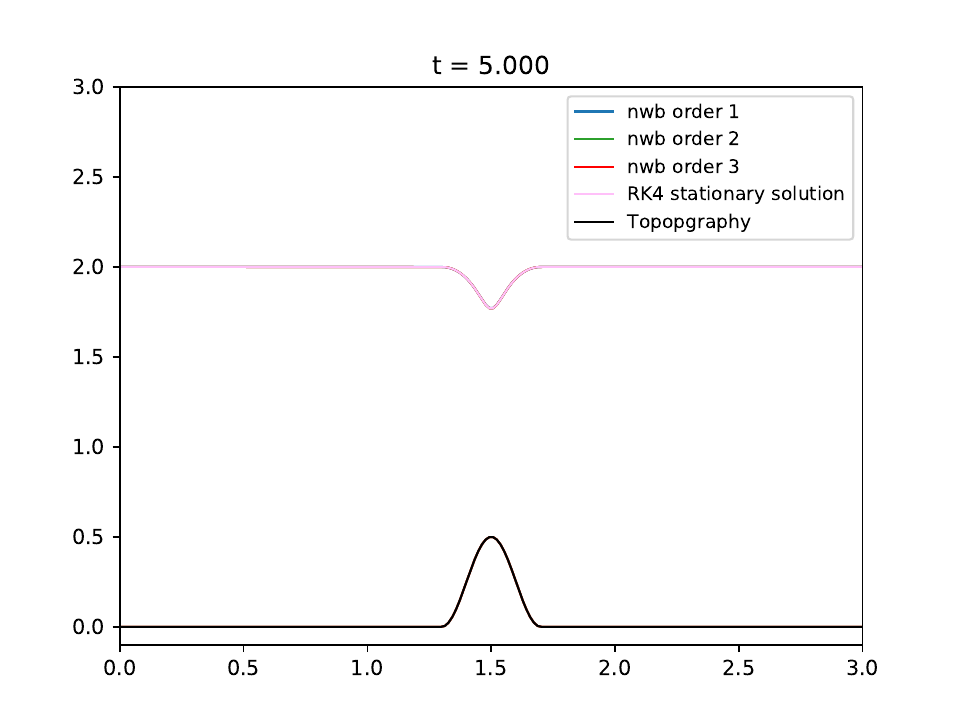}}
   \subfloat[DWBM$i$, $i=1,2,3$.  Velocity.]{
   \includegraphics[width=0.5\textwidth]{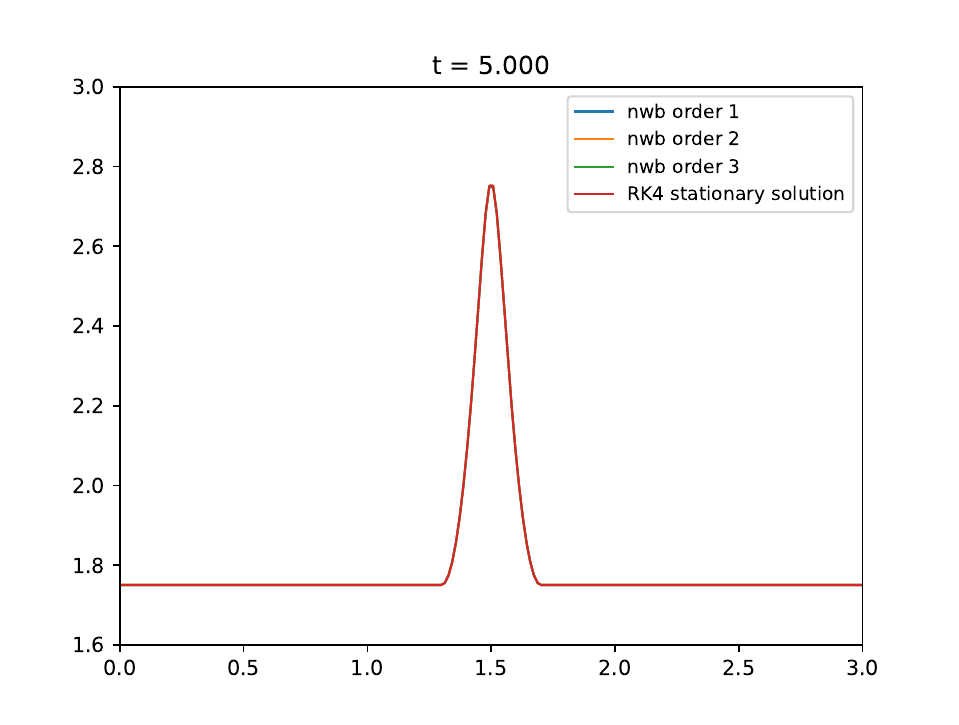}}\vspace{0.0001mm}
    \caption{Test 4.1. Numerical solution at $t = 5s$. Number of cells: 200.} \label{test6}

\end{center}
\end{figure}

Let us consider a test case taken from \cite{lopez2013}: $x \in [0,3]$, $t \in [0,5]$, and the depth function is given by:
\begin{equation}
H(x)= \left\{ \begin{array}{lcc}
             -0.25(1+\cos(5 \pi (x+0.5))) &   \text{if}  & 1.3 \leq x \leq 1.7, \\
             \\ \hspace{20mm} 0 &  &\text{otherwise}. 
             
             \end{array}
   \right.
\end{equation}
As initial condition, we consider the subcritical stationary solution that solves the Cauchy problem:
\begin{equation}
\begin{cases}
q_x=0,\\
h_x= \displaystyle \frac{ghH_x}{- \displaystyle \frac{q^2}{h^2} + gh },\\
h(0)=2,\, q(0)=3.5,
\end{cases}
\end{equation}
see Figure \ref{test6_1}. 
Conditions
$$
h(-1, t) = 2, \quad q(-1,t) = 3.5
$$
are imposed at the left extreme of the interval and free boundary conditions at $x = 1$. The CFL parameter is set again to 0.9. The conclusions are similar to the previous
test cases: Figure \ref{test6} shows the numerical solutions at  $t = 5s$ obtained with SM$i$, $i=1,2,3$ and DWBM$i$, $i=1,2,3$;
Tables \ref{ex6_error_nwb} and \ref{ex6_error_newton_wb} show the errors corresponding to the different methods. 
The results are also similar to the obtained in \cite{lopez2013}. 

\begin{figure}[H]
\begin{center}

 \includegraphics[width=0.99\textwidth]{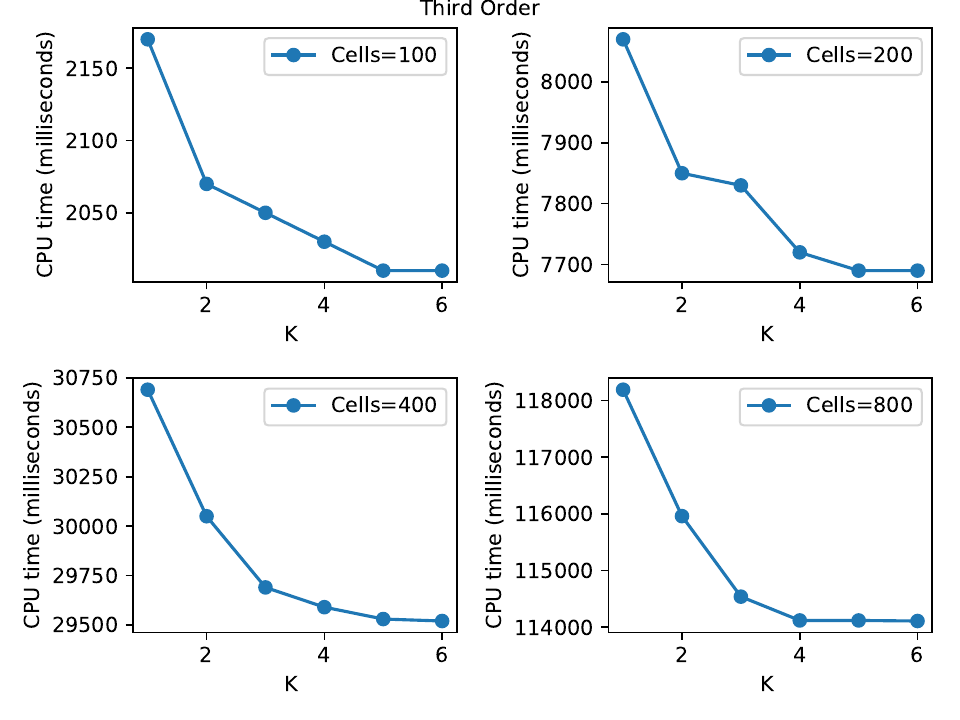}

  \caption{Test 4.1. CPU times corresponding to DWBM3 with different number of cells and different values of $K$} \label{k-tiempo}
\end{center}
\end{figure}

Concerning the computational times, we have checked the effect of using Newton's method or its modification in which $\lambda_k(0)$ in \eqref{NewtonSW} is  recomputed every $K$ iterations. 
Since, in this case, the maximum number of iterations of Newton's method throughout the computations is 6, we have compared the computational effort for values
of $K$ ranging from 1 ($\lambda_k(0)$ is recomputed at every iteration) to 6 (it is only computed once at the beginning in all cases): Figure \ref{k-tiempo} 
shows the CPU times for the third order method. As it can be seen, in this case the best option is to solve the adjoint problem only once at the beginning. 

% The results and conclusions are similar for first and second order methods.

%\begin{figure}[H]
%\begin{center}
%  \subfloat[Free surface and topography ]{
%   \includegraphics[width=0.5\textwidth]{images/test5_nwb_ok.png}}
% \subfloat[Velocity]{
%   \includegraphics[width=0.5\textwidth]{images/test5_nwb_vel.png}}
%    \caption{Solution at $t = 5s$: first, second and third order non well-balanced schemes.
%Number of cells: 200.} \label{test6_nwb}
% \end{center}
%\end{figure}
%
%\begin{figure}[H]
%\begin{center}
%  \subfloat[Free surface and topography ]{
%   \includegraphics[width=0.5\textwidth]{images/test5_wb_ok.png}}
% \subfloat[Velocity]{
%   \includegraphics[width=0.5\textwidth]{images/test5_wb_vel.png}}
%    \caption{Solution at $t = 5s$: first, second and third order well-balanced schemes.
%Number of cells: 200.} \label{test6_wb}
% \end{center}
%\end{figure}

%As expected, the approximate well-balanced schemes preserve the subcritical solution, whereas the non well-balanced methods perturb it. 

\begin{table}[H]
\centering
\begin{tabular}{|c|cc|cc|cc|} \hline
Cells & Error ($i=1$)  &Order&Error ($i=2$)  & Order&Error ($i=3$) &Order\\
 & \multicolumn{2}{|c|}{$h$}  & \multicolumn{2}{|c|}{$h$} & \multicolumn{2}{|c|}{$h$} \\\hline
100& 5.16E-2 & - & 9.39E-3  & - & 5.98E-3 & -\\
200 & 2.58E-2 & 1.000 & 2.23E-3 & 2.074 & 9.16E-4 & 2.707\\
400 & 1.28E-2 & 1.011 & 5.36E-4 & 2.057 & 1.21E-4 & 2.920\\
800 & 6.34E-3 & 1.014 & 1.30E-4 & 2.044 & 1.60E-5 & 2.919\\  \hline
Cells & Error ($i=1$)  &Order&Error ($i=2$)  & Order&Error ($i=3$) &Order\\
 & \multicolumn{2}{|c|}{$q$}  & \multicolumn{2}{|c|}{$q$} & \multicolumn{2}{|c|}{$q$} \\\hline
100& 1.94E-1 & - & 3.51E-2  & - & 2.12E-2 & -\\
200 & 9.74E-2 & 0.994 & 8.43E-3 & 2.058 & 3.23E-3 & 2.714\\
400 & 4.83E-2 & 1.012 & 2.01E-3 & 2.068 & 4.26E-4 & 2.923\\
800 & 2.40E-2 & 1.009 & 4.89E-4 & 2.039 & 5.47E-5 & 2.961\\  \hline
\end{tabular}
\caption{Test 4.1. Errors in $L^1$ norm and convergence rates for SM$i$, $i=1,2,3$.} \label{ex6_error_nwb}

\end{table}

%\begin{table}[H]
%\centering
%\begin{tabular}{|c|cc|cc|cc|} \hline
%Cells & \multicolumn{2}{|c|}{Error ($i=1$)} &\multicolumn{2}{|c|}{Error ($i=2$)}&\multicolumn{2}{|c|}{Error ($i=3$)}\\

%  & $h$&$q$ &$h$& $q$ &$h$&$q$\\\hline
%100& 1.86E-8 & 3.86E-8 & 1.92E-8  & 4.61E-8 & 3.75E-8 & 9.16E-8\\
%200 & 3.83E-9 & 1.17E-8 & 6.04E-10 & 1.45E-9 & %1.40E-9 & 3.39E-9\\
%400 & 2.13E-9 & 7.75E-9 & 1.92E-11 & 4.51E-11 & 6.50E-11 & 1.58E-10\\
%800 & 4.76E-9 & 1.76E-8 & 8.07E-13 & 1.15E-12 & 3.24E-12 & 8.53E-12\\ \hline
%\end{tabular}
%\caption{Test 4.1. Errors in $L^1$ norm for DWBM$i$, $i=1,2,3$. $N_p=3$} \label{ex6_error_newton_wb}

%\end{table}

\begin{table}[H]
\centering
\begin{tabular}{|c|cc|cc|cc|} \hline
Cells & \multicolumn{2}{|c|}{Error ($i=1$)} &\multicolumn{2}{|c|}{Error ($i=2$)}&\multicolumn{2}{|c|}{Error ($i=3$)}\\

  & $h$&$q$ &$h$& $q$ &$h$&$q$\\\hline
100& 6.15E-6 & 4.19E-6 & 6.01E-6  & 4.26E-6 & 7.45E-6 & 1.83E-5\\
200 & 2.00E-7 & 1.35E-7 & 1.97E-7 & 1.33E-7 & 1.93E-7 & 4.70E-7\\
400 & 7.46E-9 & 4.64E-9 & 7.41E-9 & 4.61E-9 & 6.63E-9 & 1.62E-8\\
800 & 2.67E-10 & 1.76E-10 & 2.66E-10 & 1.75E-10 & 2.42E-10 & 5.94E-10\\ \hline
\end{tabular}
\caption{Test 4.1. Errors in $L^1$ norm for DWBM$i$, $i=1,2,3$. } \label{ex6_error_newton_wb}

\end{table}

%Notice that the errors for the non well-balanced schemes decrease when increasing the order, and the convergence rates in Table \ref{ex6_error_nwb} are the ones expected. Moreover, the errors for the  well-balanced schemes are of the order of the machine precision. Some results concerning to computational times are included in Table \ref{ex6_times}.

\begin{table}[H]
\centering
\begin{tabular}{|c|c|c|c|}
\hline
          Cells        & $i$ & SM$i$ & DWBM$i$\\  \hline
\multirow{3}{*}{100} & 1 & 60 & 150 \\ %\cline{2-5} 
                  & 2& 160 & 750  \\ %\cline{2-5} 
                  & 3& 300 & 2010  \\ \hline
\multirow{3}{*}{200} & 1  & 220 & 500 \\ %\cline{2-5} 
                  & 2  & 490 & 2790 \\ %\cline{2-5} 
                  & 3 & 1020 & 7690 \\ \hline
%\multirow{3}{*}{400} & $1^{st}$ O & 770 & 12520 &  13170\\ %\cline{2-5} 
%                  & $2^{nd}$ O & 1790 & 39870 & 41340 \\ %\cline{2-5} 
%                  & $3^{rd}$ O & 3920 & 86550 & 90020 \\ \hline
%\multirow{3}{*}{800} & $1^{st}$ O & 2880 & 48790 & 50310 \\ %\cline{2-5} 
%                  & $2^{nd}$ O & 6870 & 155160 & 159380 \\ %\cline{2-5} 
%                  & $3^{rd}$ O & 15580 & 340420 & 349680 \\ \hline
\end{tabular}
\caption{Test 4.1. Computational times (milliseconds).}\label{ex6_times}
\end{table}

%Although the non well-balanced methods are obviously less expensive, the computational cost for well-balanced schemes is affordable. Moreover, the Newton's method is faster than the descent method, since the maximum number of iterations required to solve the nonlinear problem at every cell is smaller.

\medskip

In order to check the sensibility of the well-balanced property to the way in which the initial condition is computed, instead of using RK4 we compute now the initial condition on the basis of the implicit representation \eqref{incondimpl} of the initial condition: once the constants $C_1$ and $C_2$ have been selected, for any given $x$ the value of $h(x)$ is obtained by solving a third order polynomial equation. Table \ref{ex6_error_wb_newton_ec} shows the errors in $L^1$ norm: as it can be checked they are similar to those obtained by approximating the initial condition using RK4.

\begin{comment}
\begin{figure}[H]
\begin{center}
   \includegraphics[width=0.5\textwidth]{images/fondoec.png}
    \caption{Depth function} \label{fondoec}
 \end{center}
\end{figure}

If we use the notation $\bar{h}$ to denote the left state (in Figure \ref{fondoec}, $\bar{h}=h(a)$), the stationary solutions satisfy:
\begin{equation}
\frac{q^2}{2 h^2} + gh= \frac{q^2}{2 \bar{h}^2} + g \bar{h} +\Delta H,
\end{equation}
where $\Delta H$ is the difference between the right and left states of $H$ (in Figure \ref{fondoec}, $\Delta H = H(x_0) - H(a)$) .

Let $\Phi$ be the following function:
\begin{equation}
\varphi(h)= \frac{q^2}{2 h^2} + gh.
\end{equation}
Let us consider, at every cell, a mesh where the quadrature nodes used to approximate the objective function integral in the control problem are included. Then, each initial cell average can be obtained by solving the equation
\begin{equation}
\varphi(h) = \varphi(\bar{h}) + \Delta H, 
\end{equation}
at every node of the mesh, where $\bar{h}$ is the left state and $\Delta H$ is the difference of states of $H$, and by applying then the same quadrature rule used to approximate the objective function.

\end{comment}
\medskip

\begin{table}[H]
\centering
\begin{tabular}{|c|cc|cc|cc|} \hline
Cells & \multicolumn{2}{|c|}{Error ($i=1$)} &\multicolumn{2}{|c|}{Error ($i=2$)}&\multicolumn{2}{|c|}{Error ($i=3$)}\\

  & $h$&$q$ &$h$& $q$ &$h$&$q$\\\hline
100& 4.68E-6 & 1.20E-6 & 4.47E-6  & 8.47E-6 & 9.36E-6 & 2.56E-5\\
200 & 2.17E-7 & 5.25E-7 & 2.21E-7 & 5.50E-7 & 2.42E-7 & 6.17E-7\\
400 & 7.25E-9 & 1.70E-9 & 7.92E-9 & 1.97E-8 & 8.01E-9 & 2.00E-8\\
800 & 2.02E-10 & 4.82E-10 & 3.10E-10  & 7.69E-10 & 3.28E-10 & 8.20E-10\\ \hline
\end{tabular}
\caption{Test 4.1. Initial condition computed by using \eqref{incondimpl}. Errors in $L^1$ norm for DWBM$i$.} \label{ex6_error_wb_newton_ec}

\end{table}

%\begin{table}[H]
%\centering
%\begin{tabular}{|c|cc|cc|cc|} \hline
%Cells & \multicolumn{2}{|c|}{Error ($1^{st}$)} &\multicolumn{2}{|c|}{Error ($2^{nd}$)}&\multicolumn{2}{|c|}{Error ($3^{rd}$)}\\
%  & $h$&$q$ &$h$& $q$ &$h$&$q$\\\hline
%100& 2.24E-8 & 4.80E-8 & 2.47E-8  & 6.15E-8 & 4.96E-8 & 1.30E-7\\
%200 & 5.99E-9 & 1.98E-8& 7.96E-10 & 1.86E-9 & 1.94E-9 & 4.91E-9\\
%400 & 3.30E-9 & 1.08E-8 & 2.39E-11 & 1.19E-11 & 1.05E-10 & 2.67E-10\\
%800 & 6.62E-9 & 2.34E-8 & 3.19E-12  & 6.51E-12 & 1.00E-11 & 2.66E-12\\ \hline
%\end{tabular}
%\caption{Errors in $L^1$ norm for the well-balanced schemes (Descent method).} \label{ex6_error_wb_ec}

%\end{table} 

%Conclusions concerning to the computational cost are similar to the ones included in table \ref{ex6_times}.

\subsubsection{Test 4.2}
\begin{figure}[H]
\begin{center}
  \subfloat[Free surface and bottom. ]{
   \includegraphics[width=0.5\textwidth]{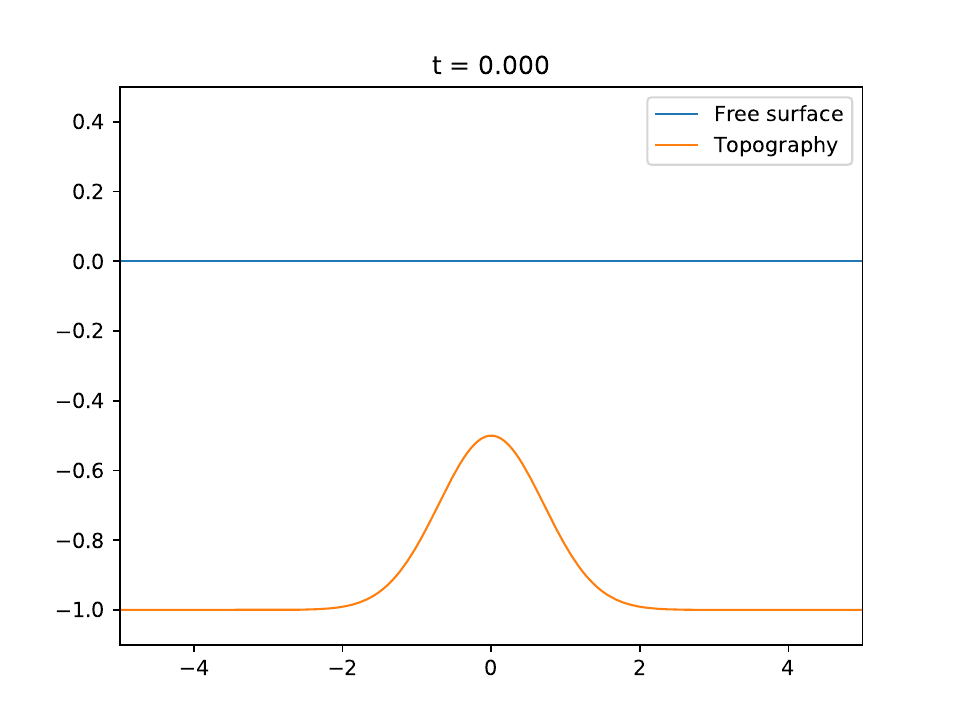}}
 \subfloat[Velocity.]{
   \includegraphics[width=0.5\textwidth]{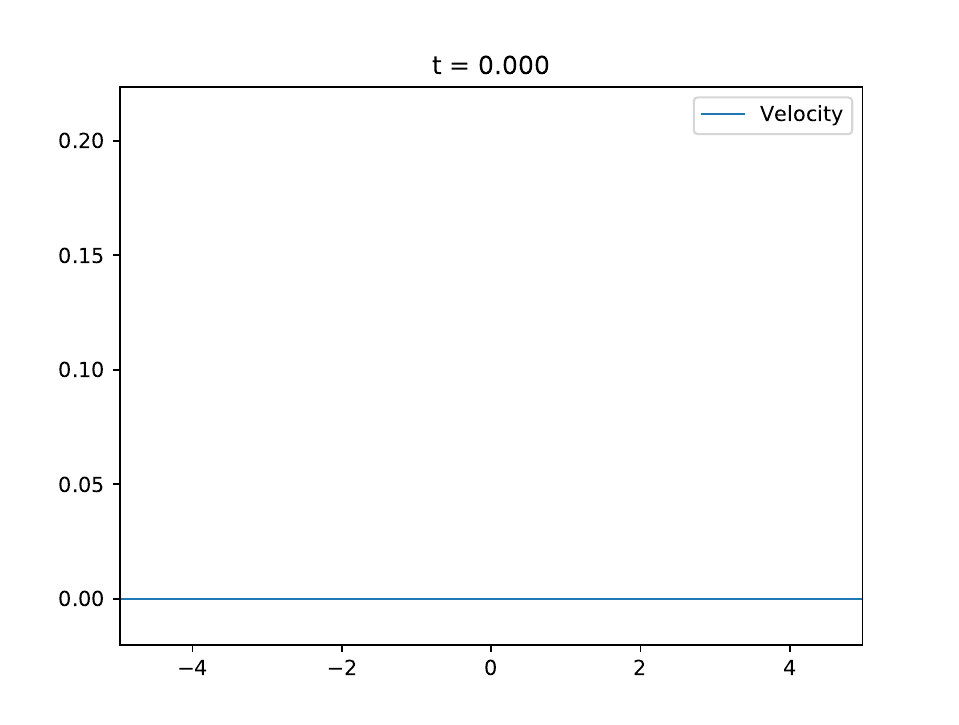}}
    \caption{Test4.2. Initial condition. } \label{test42_cini}
 \end{center}
\end{figure}

\begin{figure}[H]
\begin{center}
   \subfloat[$t = 100s$.]{
   \includegraphics[width=0.5\textwidth]{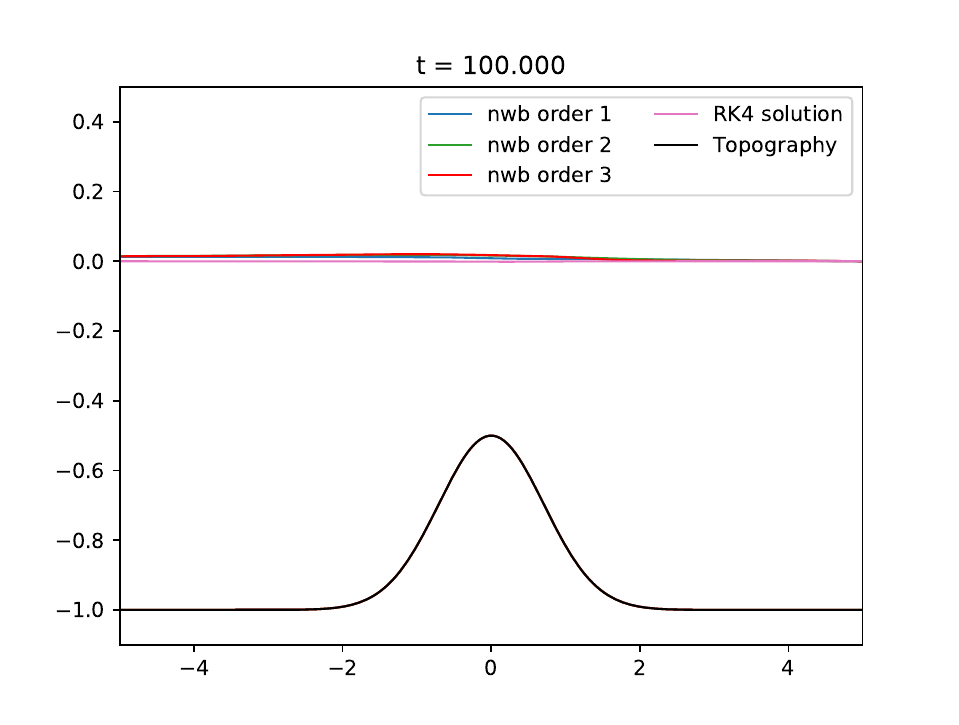}}
   \subfloat[$t = 100s$.]{
   \includegraphics[width=0.5\textwidth]{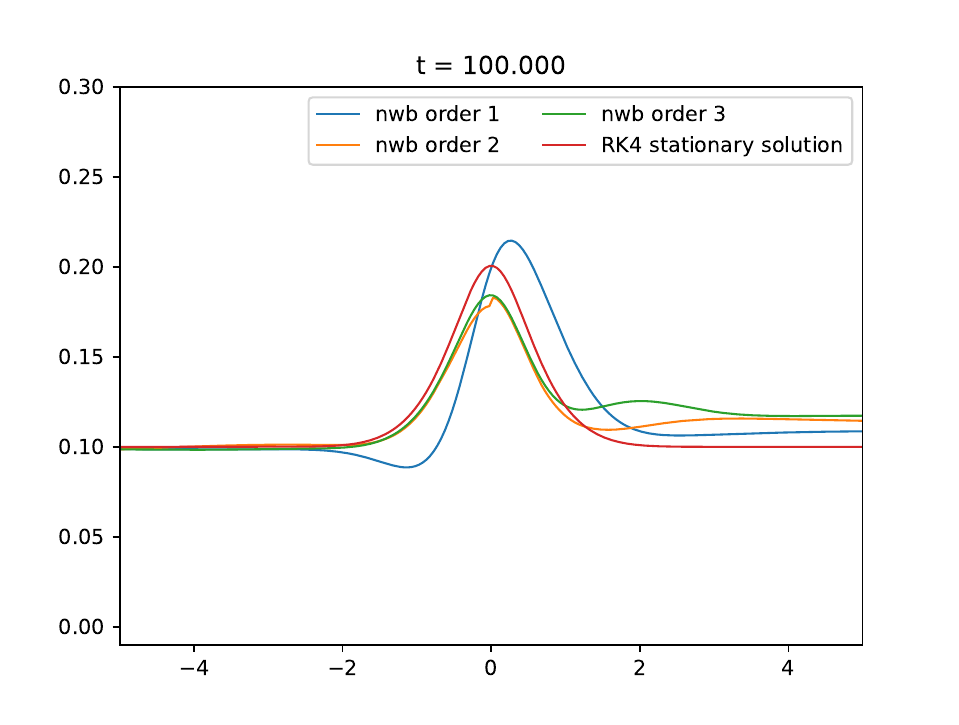}}\vspace{0.0001mm}
   \subfloat[$t = 1000s$.]{
   \includegraphics[width=0.5\textwidth]{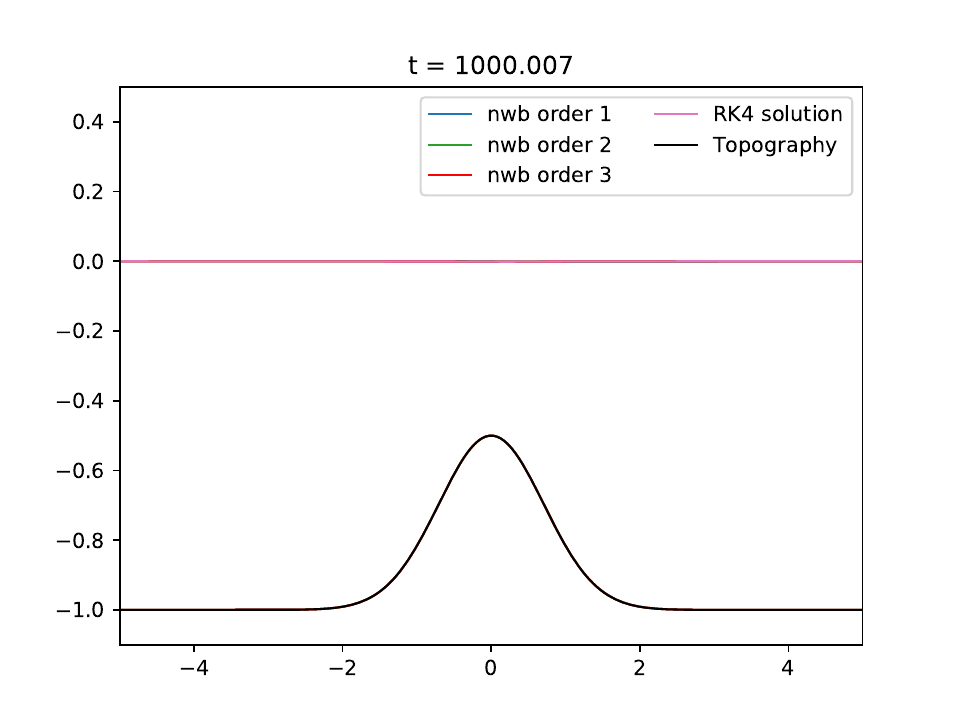}}
   \subfloat[$t = 1000s$.]{
   \includegraphics[width=0.5\textwidth]{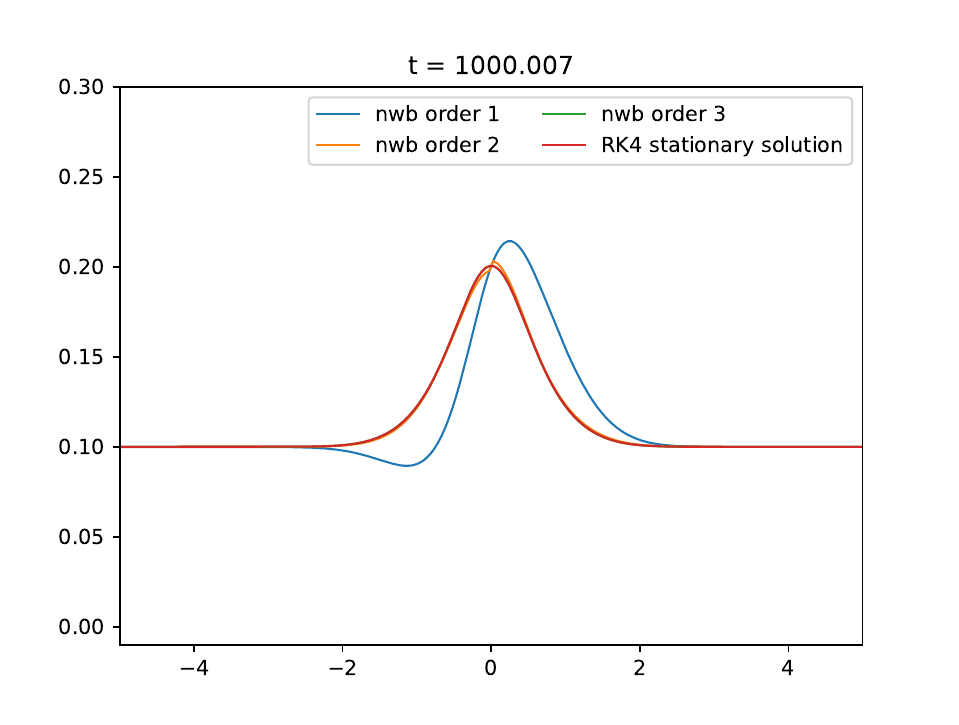}}\vspace{0.0001mm}
   \subfloat[$t = 5000s$.]{
   \includegraphics[width=0.5\textwidth]{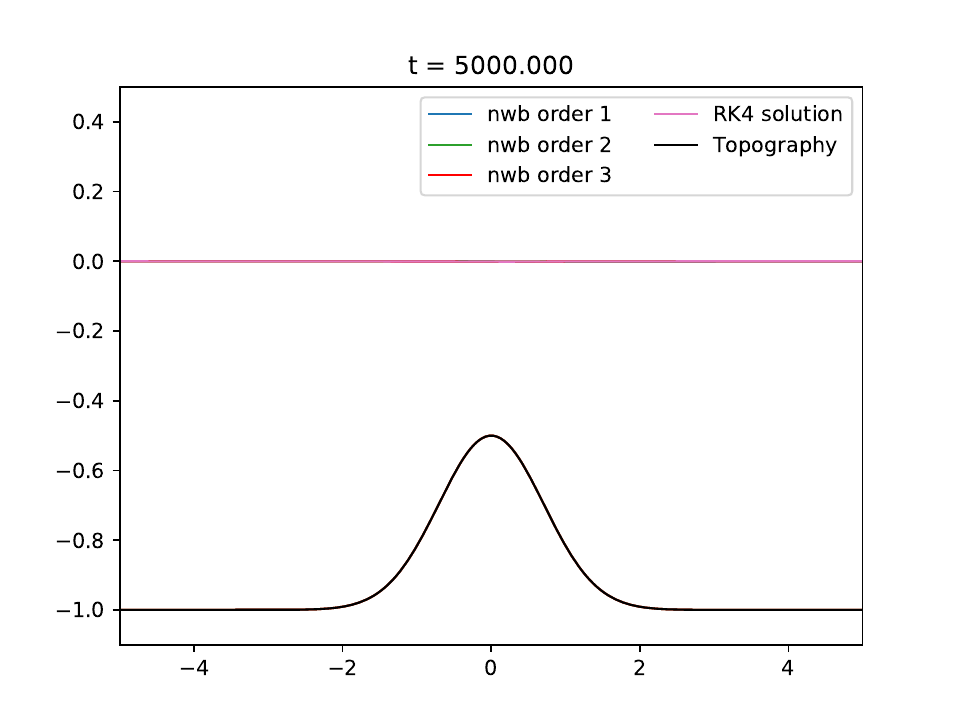}}
   \subfloat[$t = 5000s$.]{
   \includegraphics[width=0.5\textwidth]{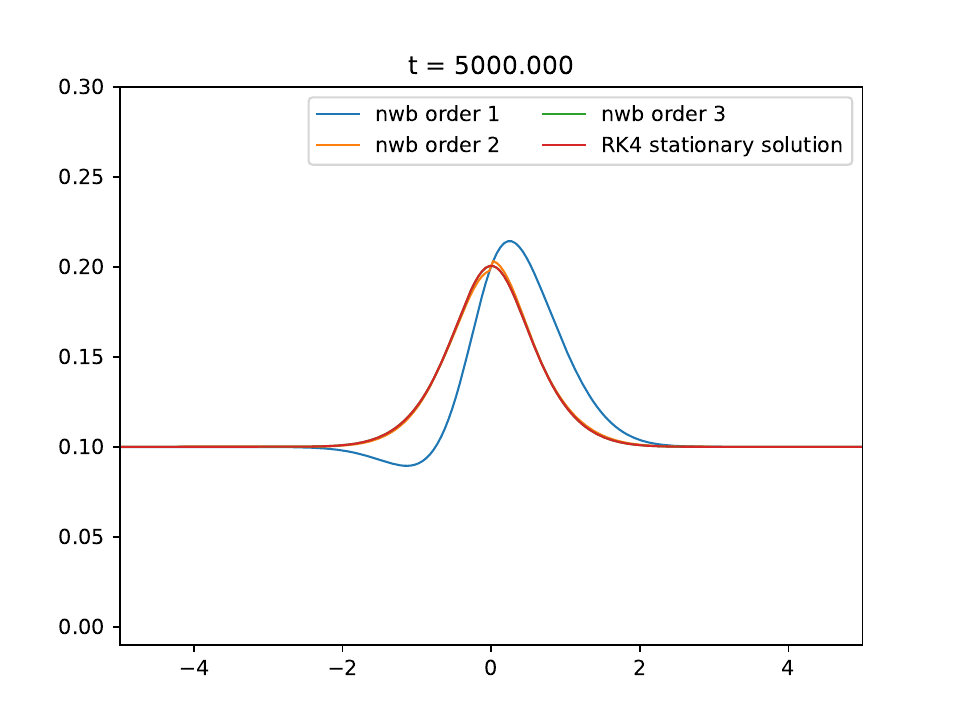}}
    \caption{Test 4.2. Numerical solutions and stationary solution computed with RK4 at times $t=100,1000, 5000s$: free surface and topography (left) and velocity (right).  SM$i$, $i=1,2,3$. Number of cells: 200.} \label{test42_nwb}
 \end{center}
\end{figure}

\begin{figure}[H]
\begin{center}
   \subfloat[$t = 100s$.]{
   \includegraphics[width=0.5\textwidth]{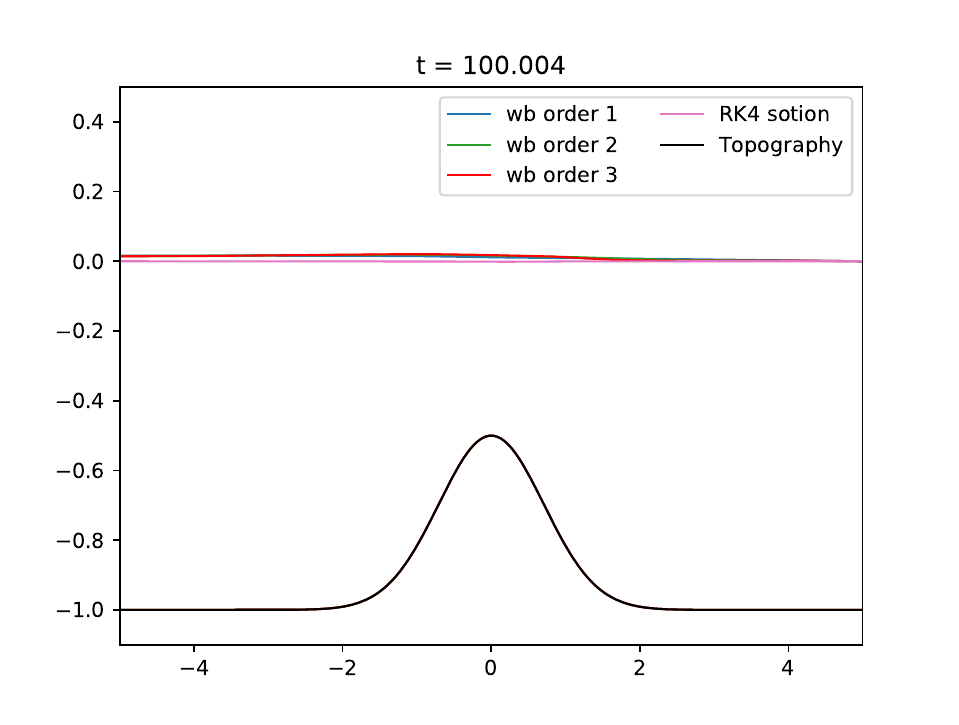}}
   \subfloat[$t = 100s$.]{
   \includegraphics[width=0.5\textwidth]{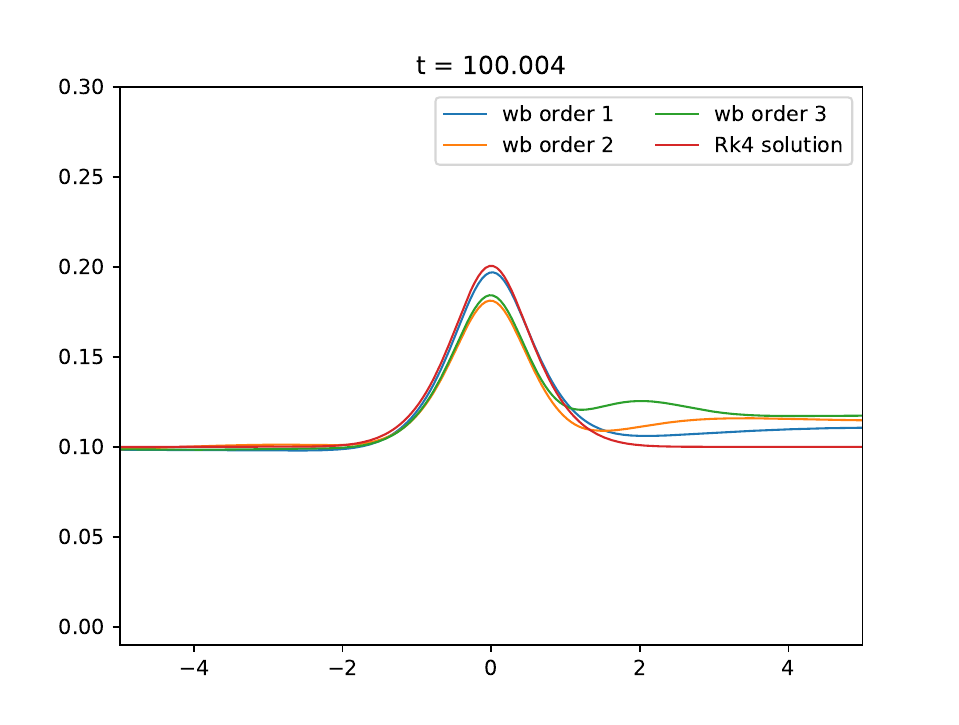}}\vspace{0.0001mm}
   \subfloat[$t = 1000s$.]{
   \includegraphics[width=0.5\textwidth]{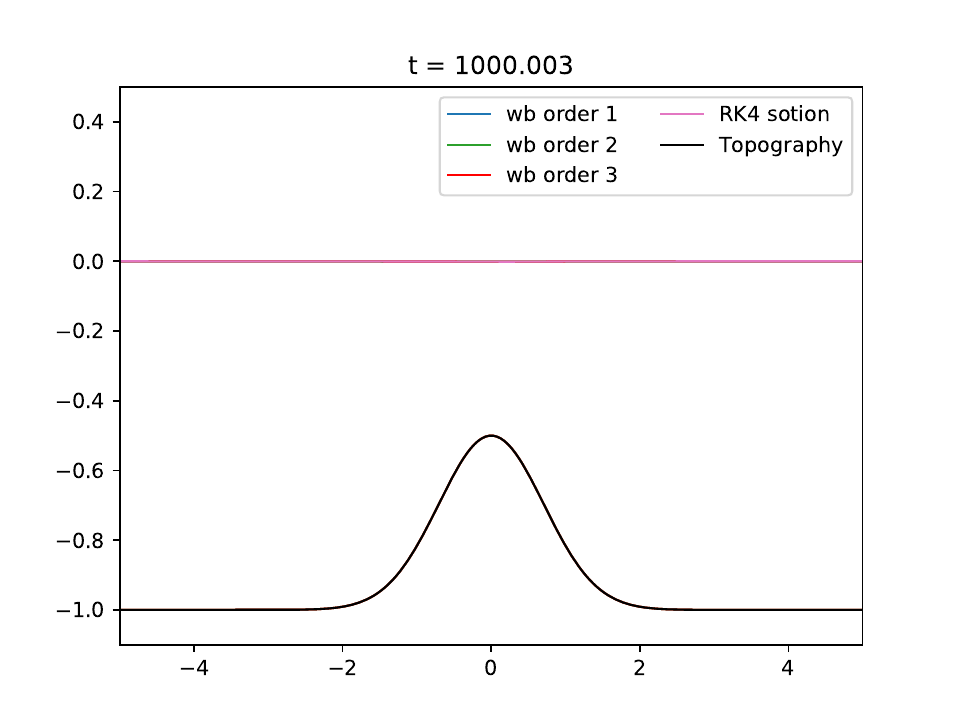}}
   \subfloat[$t = 1000s$.]{
   \includegraphics[width=0.5\textwidth]{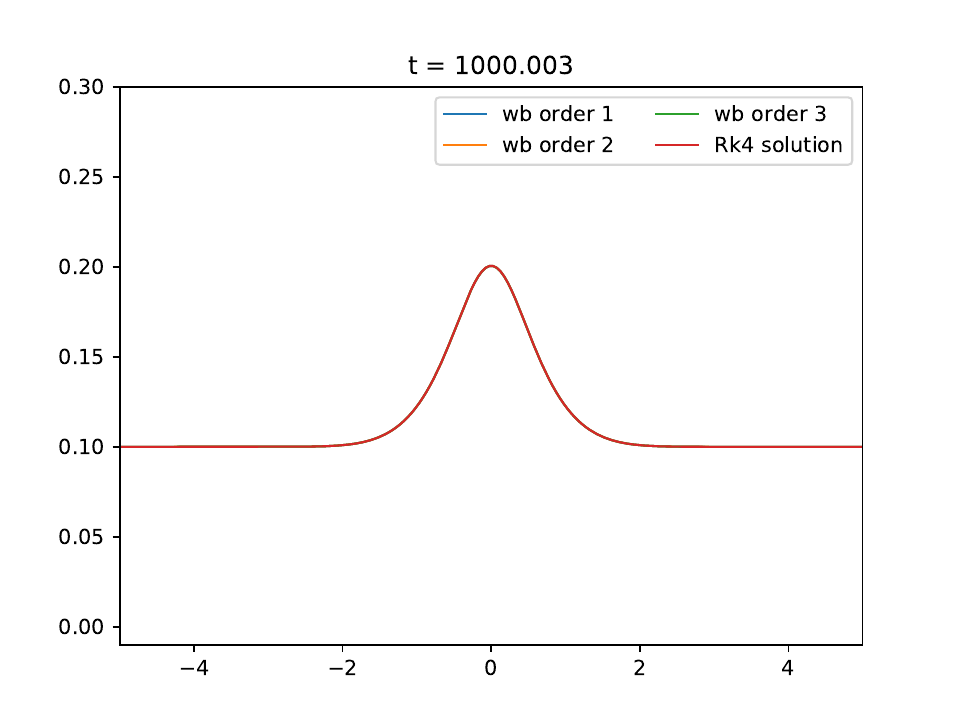}}\vspace{0.0001mm}
   \subfloat[$t = 5000s$.]{
   \includegraphics[width=0.5\textwidth]{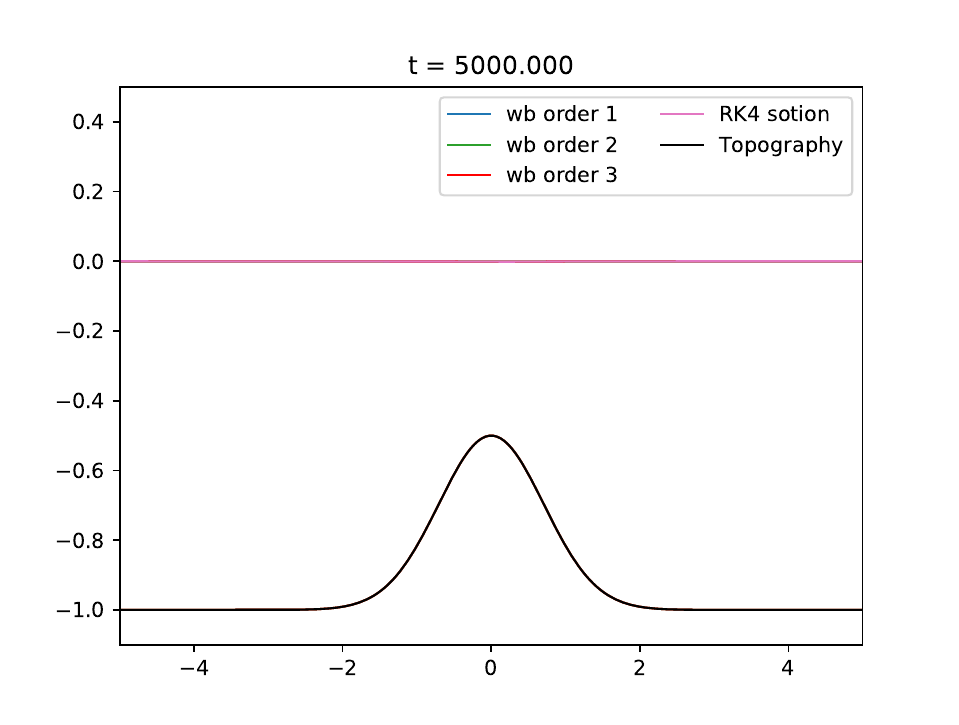}}
   \subfloat[$t = 5000s$.]{
   \includegraphics[width=0.5\textwidth]{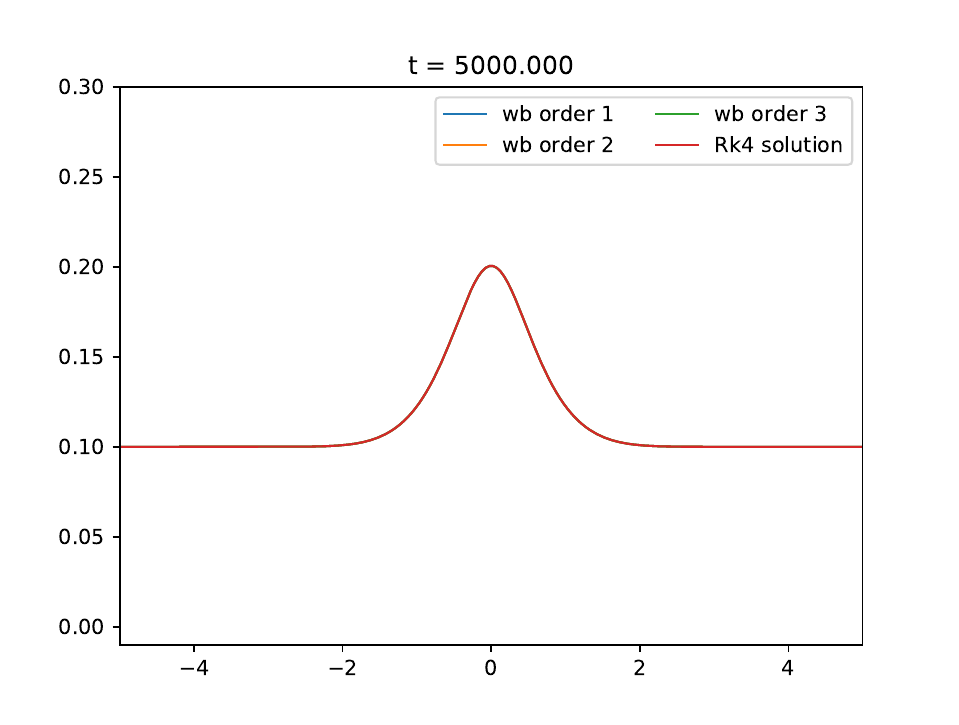}}
    \caption{{Test 4.2. Numerical solutions and stationary solution computed with RK4 at times $t=100,1000, 5000s$: free surface and topography (left) and velocity (right). DMWB$i$, $i=1,2,3$. Number of cells: 200.} } \label{test42_wb}
 \end{center}
\end{figure}

The goal of this test is to study the convergence in time of the numerical solutions to a steady state. 
We consider $x \in [-5,5]$ and $CFL=0.5$. The depth function is given by
\begin{equation}
H(x)= 1-\frac{e^{-x^2}}{2}, 
\end{equation}
and the initial condition is $h(x,0)=H(x)$ and $q(x,0)=0.0$ (see Figure \ref{test42_cini}). The imposed boundary conditions are the following
$$
q(-5,t) = 0.1, \quad h(5,t) = 1.0.
$$
The numerical solution is run in a time interval large enough so that a stationary state is reached: $t \in [0,5000]$. Figures \ref{test42_nwb} and \ref{test42_wb} show the evolution of the numerical solution for SM$i$ and NWB$i$, $i=1,2,3$. As it can be checked the well-balanced methods converge faster and better to the stationary solution. Tables \ref{test42_error_nwb} and \ref{test42_error_wb} show the $L^1$ errors at time $t = 5000s$.

%The main of this test is to study the convergence of a solution to the subcritical stationary solution corresponding to this system:
%\begin{equation}
%\begin{cases}
%q_x=0,\\
%h_x= \displaystyle \frac{ghH_x}{- \displaystyle \frac{q^2}{h^2} + gh },\\
%h(-5)=1.0,\, q(-5)=0.1
%\end{cases}
%\end{equation}

\begin{table}[H]
\centering
\begin{tabular}{|c|cc|cc|cc|} \hline
Cells & Error ($i = 1$)  &Order&Error ($i=2$)  & Order&Error ($i=3$) &Order\\
 & \multicolumn{2}{|c|}{$h$}  & \multicolumn{2}{|c|}{$h$} & \multicolumn{2}{|c|}{$h$} \\\hline
50& 1.92E-2 & - & 6.69E-3  & - & 7.59E-4 & -\\
100 & 8.83E-3 & 1.143 & 1.72E-3 & 1.960 & 7.46E-5 & 3.347\\
200 & 4.24E-3 & 1.058 & 4.21E-4 & 2.031 & 8.05E-6 & 3.212\\
400 & 2.08E-3 & 1.027 & 1.03E-4 & 2.031 & 9.51E-7 & 3.081\\  \hline
Cells & Error ($i = 1$)  &Order&Error ($i = 2$)  & Order&Error ($i = 3$) &Order\\
 & \multicolumn{2}{|c|}{$q$}  & \multicolumn{2}{|c|}{$q$} & \multicolumn{2}{|c|}{$q$} \\\hline
50& 2.73E-1 & - & 4.77E-2  & - & 6.89E-3 & -\\
100 & 1.36E-1 & 1.143 & 1.18E-2 & 2.015 & 8.97E-4 & 2.941\\
200 & 6.77E-2 & 1.058 & 2.89E-3 & 2.030 & 1.13E-4 & 2.989\\
400 & 3.38E-2 & 1.027 & 7.14E-4 & 2.017 & 1.42E-5 & 2.992\\  \hline
\end{tabular}
\caption{Test 4.2. Errors in $L^1$ norm and convergence rates for SM$i$, $i=1,2,3$. $t = 5000s$} \label{test42_error_nwb}

\end{table}

\begin{table}[H]
\centering
\begin{tabular}{|c|cc|cc|cc|} \hline
Cells & \multicolumn{2}{|c|}{Error ($i = 1$)} &\multicolumn{2}{|c|}{Error ($i = 2$)}&\multicolumn{2}{|c|}{Error ($i = 3$)}\\

  & $h$&$q$ &$h$& $q$ &$h$&$q$\\\hline
50& 1.41E-8 & 2.14E-14 & 1.37E-8  & 4.73E-11 & 3.29E-8 & 4.43E-10\\
100 & 4.39E-10 & 1.10E-14 & 4.32E-10 & 4.48E-13 & 1.03E-9 & 7.26E-12\\
200 & 1.37E-11 & 1.15E-14 & 1.36E-11 & 1.57E-14 & 3.23E-11 & 3.14E-13\\
400 & 4.32E-13 & 1.83E-14 & 4.32E-13 & 5.34E-14 & 1.00E-12 & 7.88E-13\\ \hline
\end{tabular}
\caption{Test 4.2. Errors in $L^1$ norm for DWBM$i$, $i=1,2,3$. $t = 5000s$} 
\label{test42_error_wb}

\end{table}

\subsection{Problem 5: Compressible Euler equations with gravitational force}
Let us consider now the Euler equations of gas dynamics with source term for the simulation of the flow of a gas in a gravitational field:
\begin{equation}\label{euler_equations}
\begin{cases}
\rho_t+(\rho u)_x=0,\\
(\rho u)_t +\left( \rho u^2 + p \right)_x= -\rho H_x,\\
(E)_t +\left( u (E + p) \right)_x= -\rho u H_x.\\
\end{cases}
\end{equation}
Here, $\rho \geq 0$ is the density, $u$ the velocity, $q=\rho u$ the momentum, $p\geq 0$ the pressure, $E$ the total energy per unit volume, and $H(x)$ the gravitational potential. Futhermore, the internal energy $e$ is given by $\rho e = E - \displaystyle \frac{1}{2} \rho u^2$. Pressure is determined from $e$ through the equation of state. Here we suppose for simplicity an ideal gas, therefore
$$p= (\gamma -1 ) \rho e,$$
where $\gamma >1$ is the adiabatic constant: here $\gamma = 1.5$. 

System \eqref{euler_equations} is a particular case of \eqref{sle} corresponding to the choices $N=3$,
$$ U =\begin{pmatrix}
\rho \\
\rho u \\
E \\
\end{pmatrix} , \quad f(U) =\begin{pmatrix}
\rho u \\
 \rho u^2 + p\\
 u(E+p)\\
\end{pmatrix}, \quad S(U) =\begin{pmatrix}
0\\
- \rho\\
- \rho u\\
\end{pmatrix}.$$

The system of ODE satisfied by the stationary solutions is:
\begin{equation} \label{euler_edo_estacionaria}
\begin{cases}
q_x=0,\\
\left( \displaystyle \frac{q^2}{\rho} + p \right)_x= -\rho H_x,\\
\left( \displaystyle \frac{q}{\rho} (E+p) \right)_x= -q  H_x.\\
\end{cases}
\end{equation}
It can be easily checked that \eqref{euler_edo_estacionaria} can be written in the following form, as we suppose that the system is strictly hyperbolic: 

\begin{equation}\label{u'=Geuler}
\begin{cases}
q_x=0,\\
\displaystyle \frac{d \hat{U}}{dx} = G(x, \hat{U}), \\
\end{cases}
\end{equation}
where
$$ \hat{U} =\begin{pmatrix}
\rho \\
E \\
\end{pmatrix} , \quad G(x, \hat{U}) = - \begin{pmatrix}
\displaystyle \frac{\rho}{c^2 - u^2} \\
 \displaystyle \frac{\rho}{\gamma-1} \left( 1 + \frac{3-\gamma}{2} \frac{u^2}{c^2 - u^2} \right)\\

\end{pmatrix} H_x,$$
where
$$c = \sqrt{ \gamma  \frac{p}{\rho}}$$
is the wave speed,
and thus $\nabla_{\hat{U}} G $ is given by
\begin{equation}\label{nablaGeuler}
  - \left[ \begin{array}{cc}
\displaystyle \frac{c^2-u^2 +(\gamma -1) \displaystyle \frac{\gamma E}{\rho} - (\gamma(\gamma -1)+2)u^2}{(c^2-u^2)^2 } & - \displaystyle \frac{\gamma}{(c^2-u^2)^2 } \\
& \\
\displaystyle \frac{1}{\gamma -1} \left[ 1 + \displaystyle \frac{(\gamma -3) u^2 }{2(c^2-u^2)}  + \displaystyle \frac{u^2(\gamma -3)}{(c^2-u^2)^2} \left( \displaystyle \frac{\left( \gamma (\gamma -1) +2 \right) u^2}{2 (\gamma-1)}  - \frac{\gamma E}{2 \rho} \right) \right]  &  \displaystyle \frac{\gamma (\gamma -3) u^2}{(c^2-u^2)^2}\\
\end{array}
\right] H_x
\end{equation}
%\begin{equation}\label{nablaGsw}
%\nabla_u G  = \left[ \begin{array}{cc}
%\displaystyle
%-\frac{3gH_xh^2q^2}{(- {q^2}+ gh^3)^2 } & \displaystyle \frac{2gH_xh^3q}{(- {q^2}+ gh^3)^2 } \\
%& \\
%0  & 0\\
%\end{array}
%\right].
%\end{equation}
%

Reasoning like in the shallow water case, it can be shown that only $2 \times 2$ ODE systems in the variables $\rho$, $E$
have to be solved to compute the state and the adjoint variables at every step of Newton's method. 

\subsubsection{Test 5.1}

\begin{figure}[H]
\begin{center}
  \subfloat[ Density. ]{
   \includegraphics[width=0.33\textwidth]{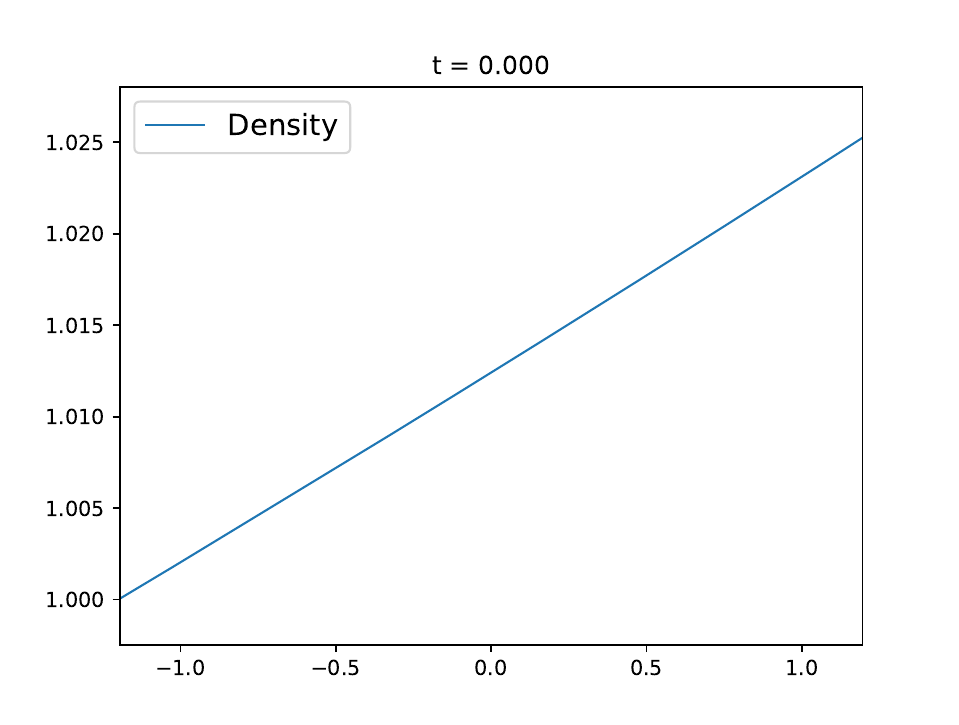}}
 \subfloat[Velocity.]{
   \includegraphics[width=0.33\textwidth]{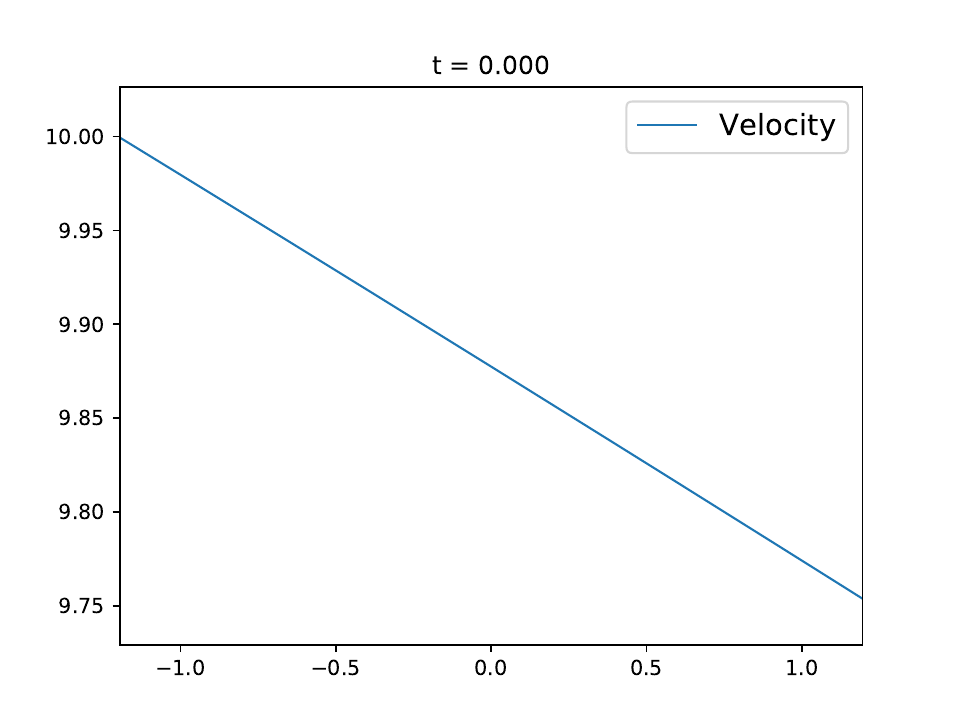}}
    \subfloat[ Energy.]{
   \includegraphics[width=0.33\textwidth]{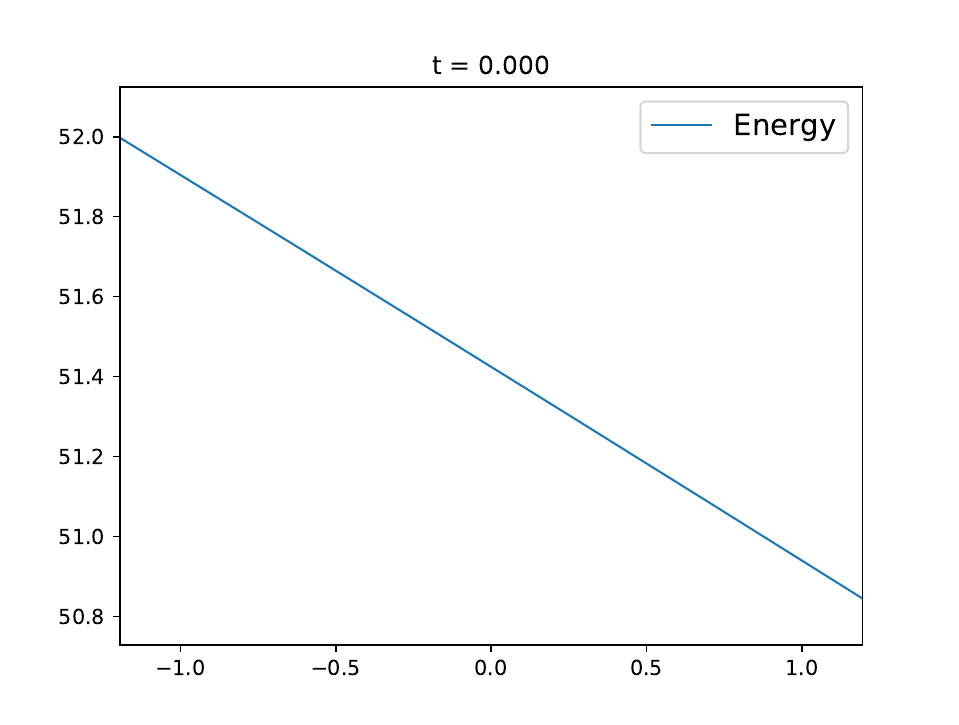}}
    \caption{Test 5.1. Initial condition: a supersonic stationary solution computed with the RK4 method.} \label{test51cini}
 \end{center}
\end{figure}

We consider $x \in [-1,1]$, $t \in [0,5]$, $CFL=0.9$, and the gravity potential is the identity function $H(x) = x$.
As initial condition, we compute using RK4 the supersonic stationary solution which solves the Cauchy problem:
\begin{equation} \label{test51_cini}
\begin{cases}
q_x=0,\\
\displaystyle \frac{d \hat{U}}{dx} = G(x, \hat{U}), \\
\rho(-1)=1, \, q(-1)=10, \, E(-1)=52.
\end{cases}
\end{equation}
See Figure \ref{test51cini}.

Boundary conditions
$$\rho(-1,t), \quad q(-1,t) = 10, \quad E(-1,t) = 52,$$
are imposed at $x=-1$ and free boundary conditions at $x = 1$. $N_p=1$ is considered, and the modified Newton's  method in which  $\Lambda(0)$ is only computed once is applied to this problem. 

Figure \eqref{test51_nwb} shows the numerical results obtained with SM$i$, $i=1,2,3$ (up) together with some zooms (down) where it can be clearly observed how the stationary solution is perturbed. 

We have made a zoom to the pictures of the non well-balanced results in order to show that these schemes perturb the stationary solution:

\begin{figure}[H]
\begin{center}
  \subfloat[Density. ]{
   \includegraphics[width=0.33\textwidth]{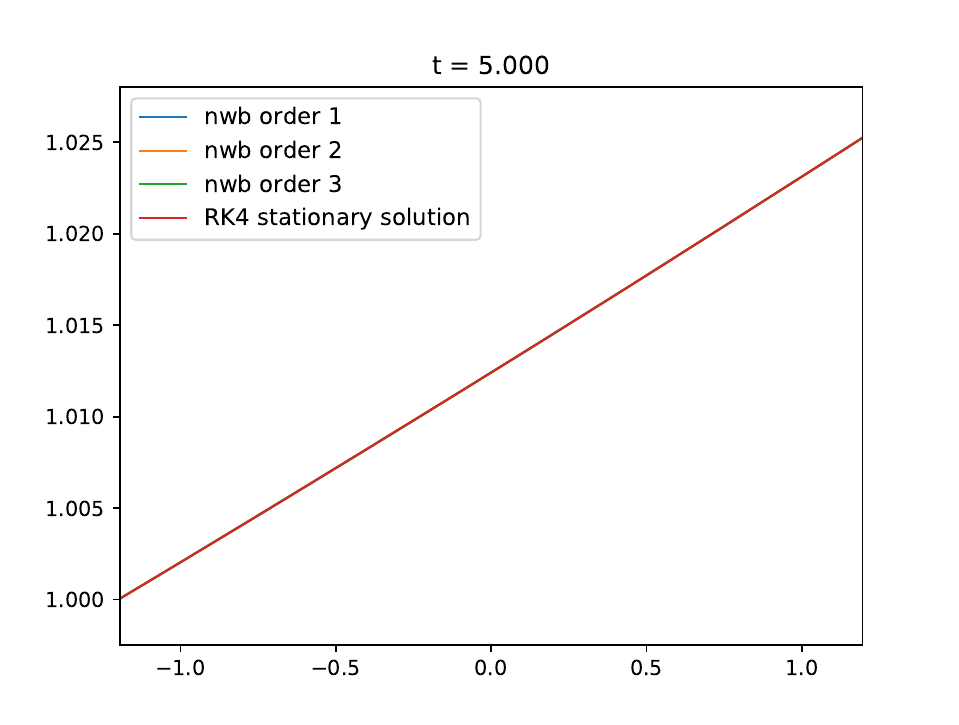}}
 \subfloat[Velocity.]{
   \includegraphics[width=0.33\textwidth]{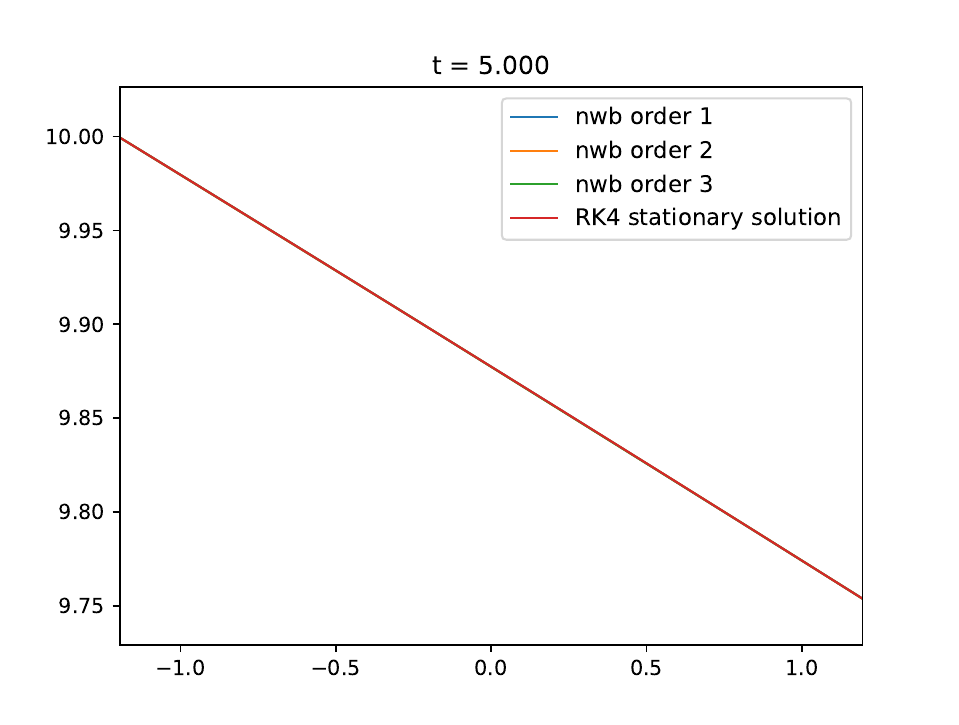}}
   \subfloat[ Energy.]{
   \includegraphics[width=0.33\textwidth]{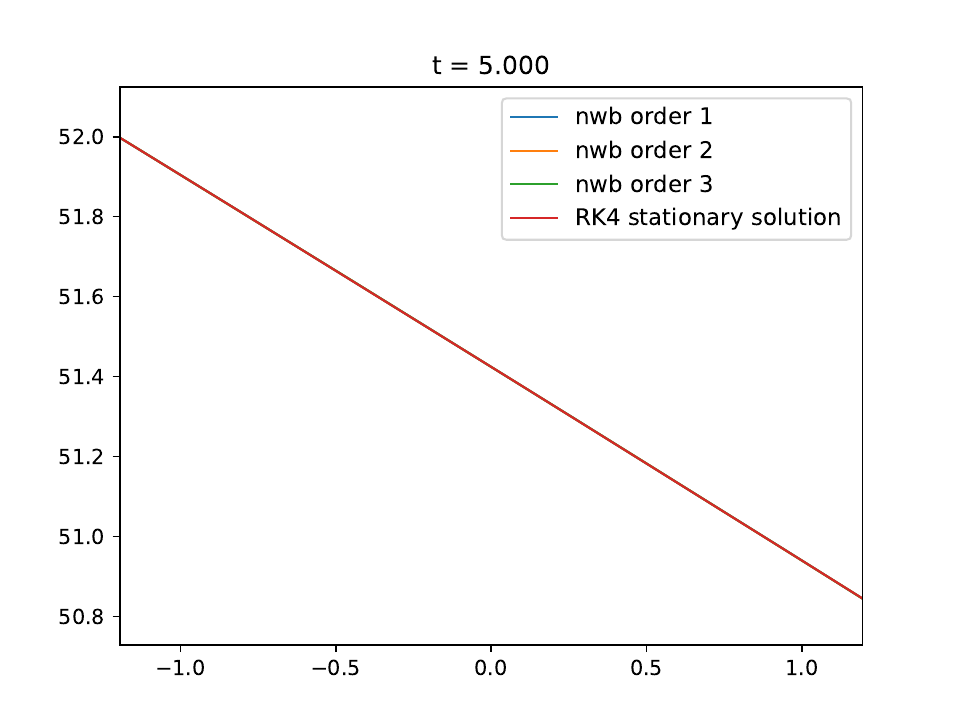}}\vspace{0.0001mm}
   \subfloat[Density. Zoom.]{
   \includegraphics[width=0.33\textwidth]{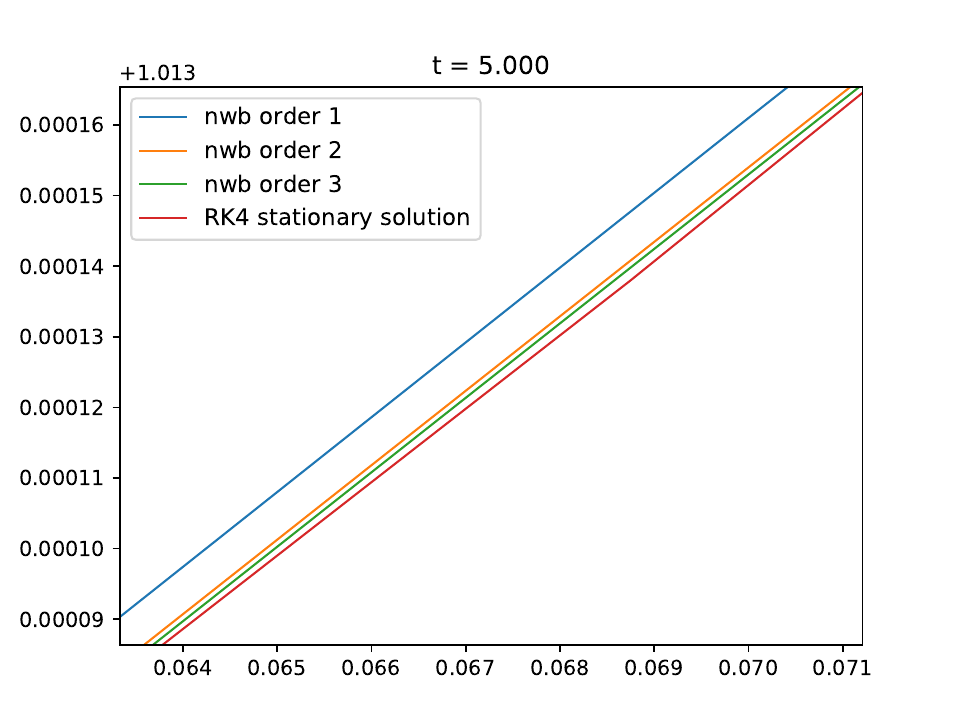}}
   \subfloat[ Velocity. Zoom.]{
   \includegraphics[width=0.33\textwidth]{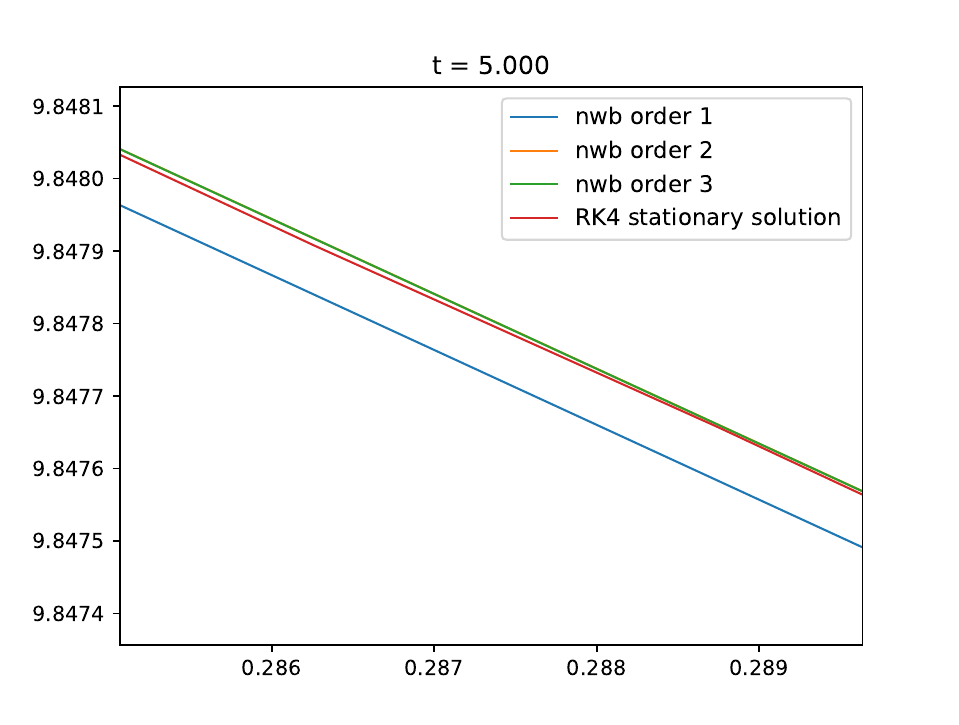}}
   \subfloat[ Energy. Zoom.]{
   \includegraphics[width=0.33\textwidth]{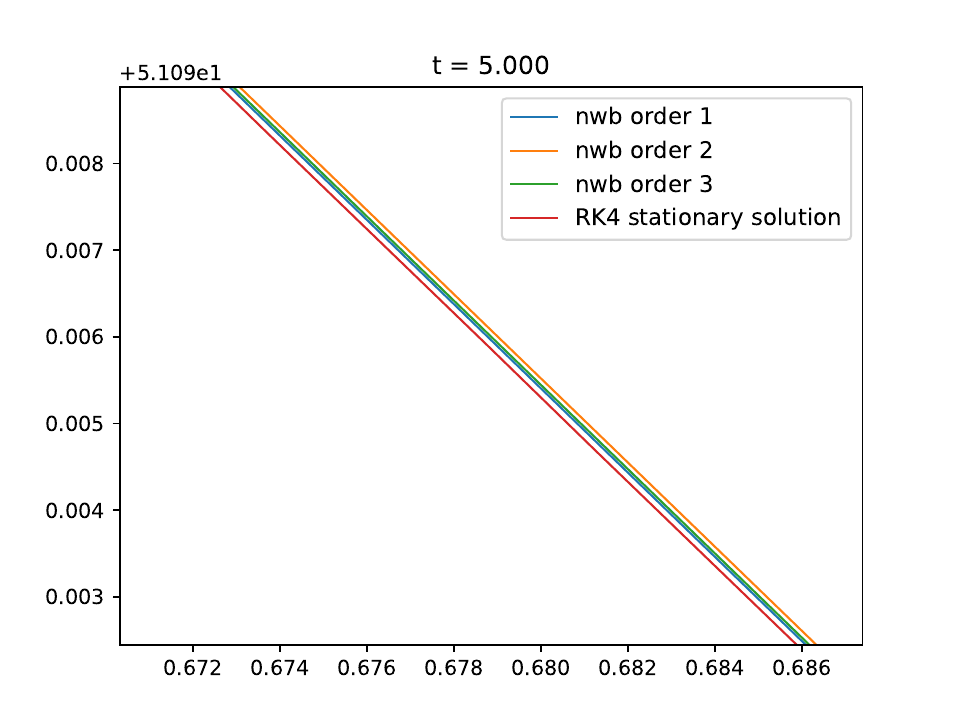}}
    \caption{Test 5.1. Numerical solutions  and stationary solution at time $t = 5s$: global view (up) and zoom (down). SM$i$, $i=1,2,3$. Number of cells: 100.} \label{test51_nwb}
 \end{center}
\end{figure}

\begin{figure}[H]
\begin{center}
  \subfloat[Density. ]{
   \includegraphics[width=0.33\textwidth]{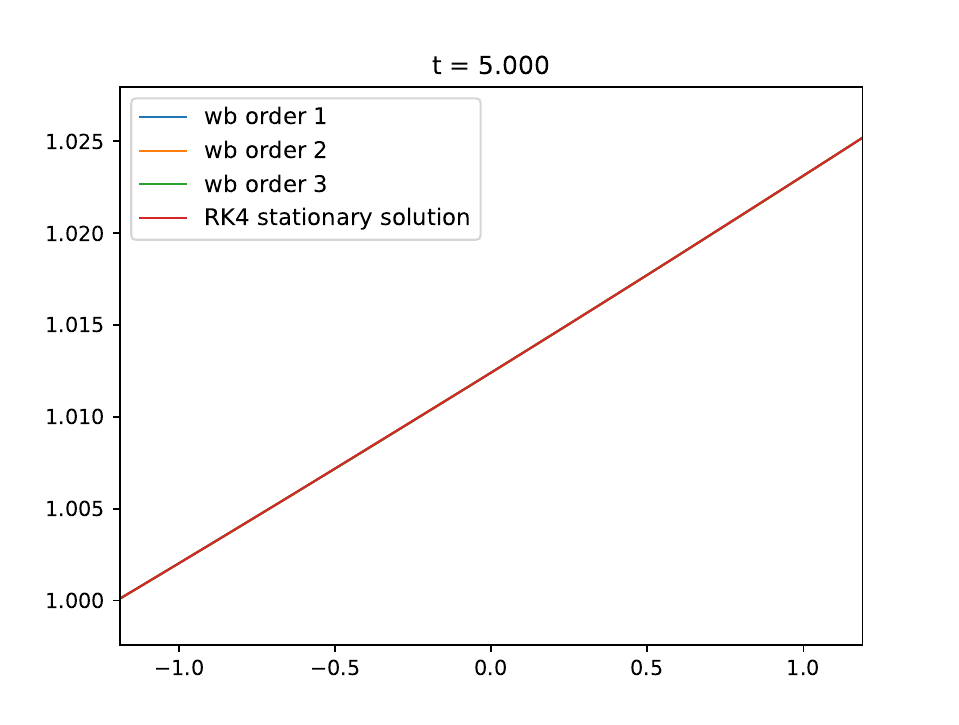}}
 \subfloat[Velocity.]{
   \includegraphics[width=0.33\textwidth]{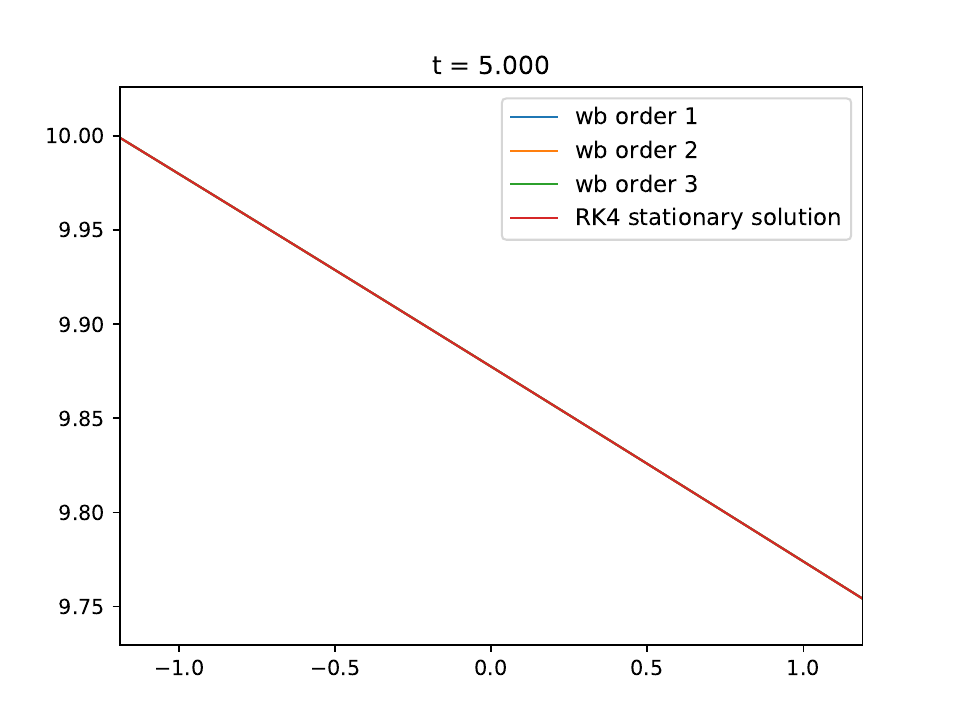}}
   \subfloat[ Energy.]{
   \includegraphics[width=0.33\textwidth]{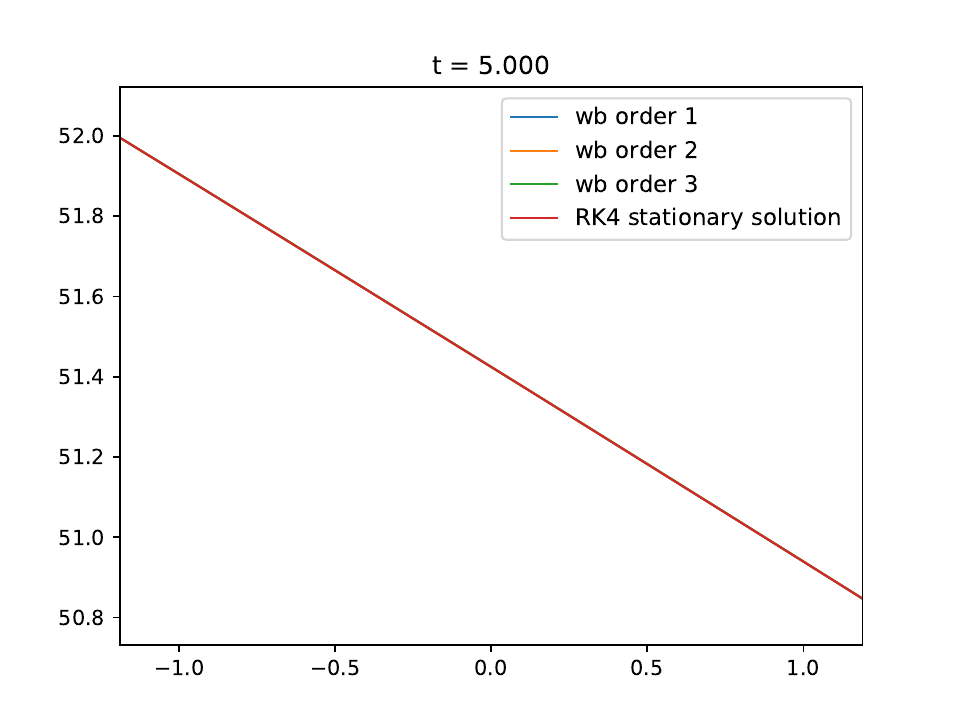}}\vspace{0.0001mm}
    \caption{Test 5.1. Numerical solutions  and stationary solution at time $t = 5s$. DWBM$i$, $i=1,2,3$. Number of cells: 100.} \label{test51_wb}
 \end{center}
\end{figure}

Tables \ref{ex51_error_nwb} and \ref{ex51_error_newton_wb} show the errors for SM$i$ and DWBM$i$, $i=1,2,3$, respectively. Computational times are shown in Table \ref{ex51_times}.

%\begin{table}[H]
%\centering
%\begin{tabular}{|c|ccc|ccc|ccc|} \hline
%Cells & \multicolumn{3}{|c|}{Error ($1^{st}$)} &\multicolumn{3}{|c|}{Error ($2^{nd}$)}&\multicolumn{3}{|c|}{Error ($3^{rd}$)}\\
%
%  & $\rho$&$q$ &$E$& $\rho$&$q$ &$E$& $\rho$&$q$ &$E$\\\hline
%100& 1.71E-5 & 5.91E-5 & 1.18E-4  & 4.30E-6 & 6.45E-5 & 4.40E-4 & 2.66E-6 & 4.10E-5 & 2.90E-4\\
%200 & 8.92E-6 & 3.90E-5 & 1.13E-4 & 7.24E-9 & 1.04E-7 & 6.90E-7 & 5.33E-10 & 1.38E-8 & 1.16E-7\\
%400 & 4.39E-6 & 1.93E-5 & 5.63E-5 & 6.38E-11 & 9.76E-12 & 3.13E-10 & 7.11E-14 & 7.13E-13 & 3.50E-12\\
%800 & 2.17E-6 & 9.55E-6 & 2.79E-5 & 1.58E-11 & 2.21E-12 & 7.89E-11 & 1.43E-13 & 1.41E-12 & 8.27E-12\\ \hline
%\end{tabular}
%\caption{Errors in $L^1$ norm for the non well-balanced schemes.} \label{ex51_error_nwb}
%
%\end{table}

\begin{table}[H]
\centering
\begin{tabular}{|c|c|c|c|} \hline
Cells & \multicolumn{1}{|c|}{Error ($i = 1$)} &\multicolumn{1}{|c|}{Error ($i = 2$)}&\multicolumn{1}{|c|}{Error ($i = 3$)}\\ \hline
  \multicolumn{4}{|c|}{$\rho$}\\\hline
100& 1.71E-5 &4.30E-6 & 2.66E-6 \\
200 & 8.92E-6 & 7.24E-9 & 5.33E-10 \\
400 & 4.39E-6& 6.38E-11 & 7.11E-14 \\
800 & 2.17E-6 & 1.58E-11 & 1.43E-13 \\ \hline
 \multicolumn{4}{|c|}{$q$}\\\hline
100& 5.91E-5 &6.45E-5 & 4.10E-5 \\
200 & 3.90E-5 & 1.04E-7 & 1.38E-8 \\
400 & 1.93E-5 & 9.76E-12 & 7.13E-13 \\
800 & 9.55E-6 & 2.21E-12 & 1.41E-12 \\ \hline 
 \multicolumn{4}{|c|}{$E$}\\\hline
100& 1.18E-4  &4.40E-4 & 2.90E-4 \\
200 & 1.13E-4 & 6.90E-7 & 1.16E-7 \\
400 & 5.63E-5 & 3.13E-10 & 3.50E-12 \\
800 & 2.79E-5 & 7.89E-11 & 8.27E-12 \\ \hline

\end{tabular}
\caption{Test 5.1. Errors in $L^1$ norm  for SM$i$, $i=1,2,3$.} \label{ex51_error_nwb}

\end{table}

\begin{table}[H]
\centering
\begin{tabular}{|c|c|c|c|} \hline
Cells & \multicolumn{1}{|c|}{Error ($i = 1$)} &\multicolumn{1}{|c|}{Error ($i = 2$)}&\multicolumn{1}{|c|}{Error ($i = 3$)}\\ \hline
  \multicolumn{4}{|c|}{$\rho$}\\\hline
100& 1.07E-14 & 2.15E-14 & 2.95E-14 \\
200 & 2.28E-14 & 4.45E-14 & 5.56E-14 \\
400 & 3.41E-14& 6.39E-14 & 8.79E-14 \\
800 & 4.27E-14 & 7.74E-14 & 1.71E-13 \\ \hline
 \multicolumn{4}{|c|}{$q$}\\\hline
100& 1.05E-13 & 2.16E-13 & 2.90E-13 \\
200 & 2.23E-13 & 4.29E-13 & 5.55E-13 \\
400 & 3.33E-13 & 6.15E-13 & 8.49E-13 \\
800 & 4.19E-13 & 7.53E-13 & 1.64E-12 \\ \hline 
 \multicolumn{4}{|c|}{$E$}\\\hline
100& 5.45E-12  & 1.10E-12 & 1.39E-12 \\
200 & 1.15E-13 & 2.15E-12 & 2.78E-12 \\
400 & 1.68E-12 & 3.04E-12 & 4.24E-12 \\
800 & 2.17E-12 & 3.72E-12 & 8.34E-12 \\ \hline

\end{tabular}
\caption{Test 5.1. Errors in $L^1$ norm  for DWBM$i$, $i=1,2,3$.} \label{ex51_error_newton_wb}

\end{table}

\begin{table}[H]
\centering
\begin{tabular}{|c|c|c|c|}
\hline
          Cells        & Order & Non well-balanced & Well-balanced  \\  \hline
\multirow{3}{*}{100} & $1^{st}$ O & 120 & 400  \\ %\cline{2-5} 
                  & $2^{nd}$ O  & 250 & 2320  \\ %\cline{2-5} 
                  & $3^{rd}$ O & 670 & 10290  \\ \hline
\multirow{3}{*}{200} & $1^{st}$ O & 430 & 1570  \\ %\cline{2-5} 
                  & $2^{nd}$ O  & 840 & 8370  \\ %\cline{2-5} 
                  & $3^{rd}$ O & 2790 & 39230  \\ \hline

\end{tabular}
\caption{Computational times (milliseconds).}\label{ex51_times}
\end{table}

\subsubsection{Test 5.2}
\begin{figure}[H]
\begin{center}
  \subfloat[Density. ]{
   \includegraphics[width=0.33\textwidth]{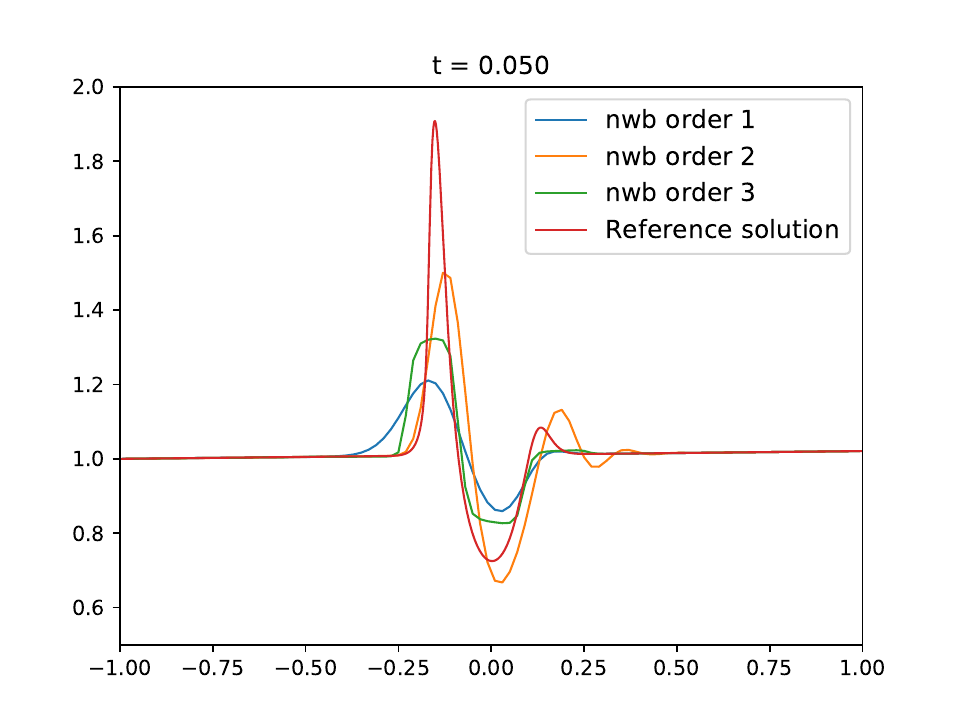}}
 \subfloat[Velocity.]{
   \includegraphics[width=0.33\textwidth]{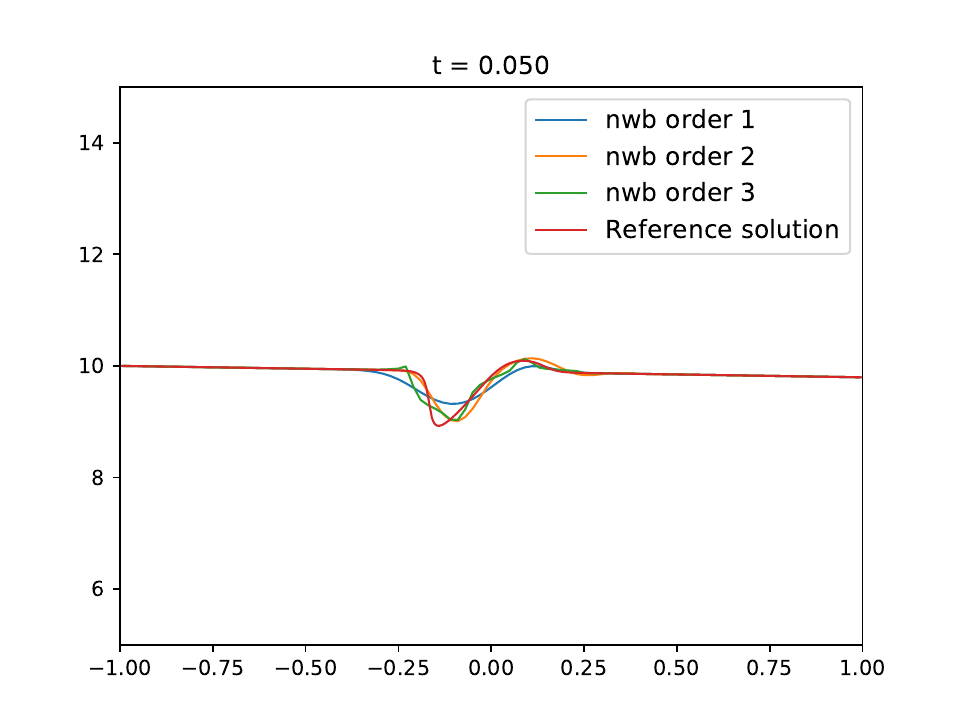}}
   \subfloat[ Energy.]{
   \includegraphics[width=0.33\textwidth]{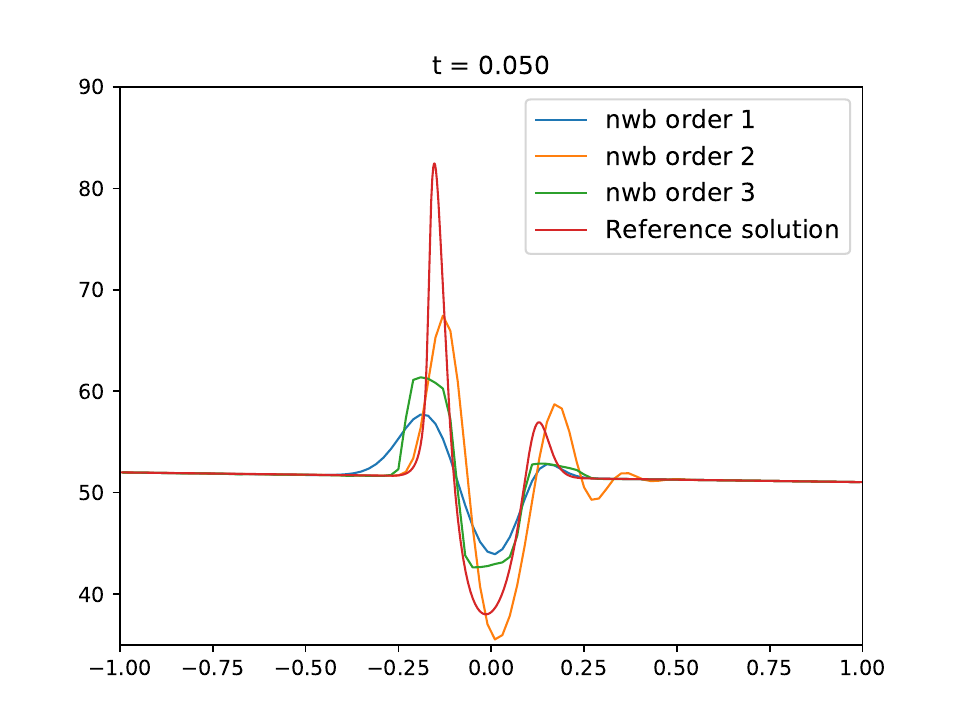}}\vspace{0.0001mm}
   \subfloat[Density. ]{
   \includegraphics[width=0.33\textwidth]{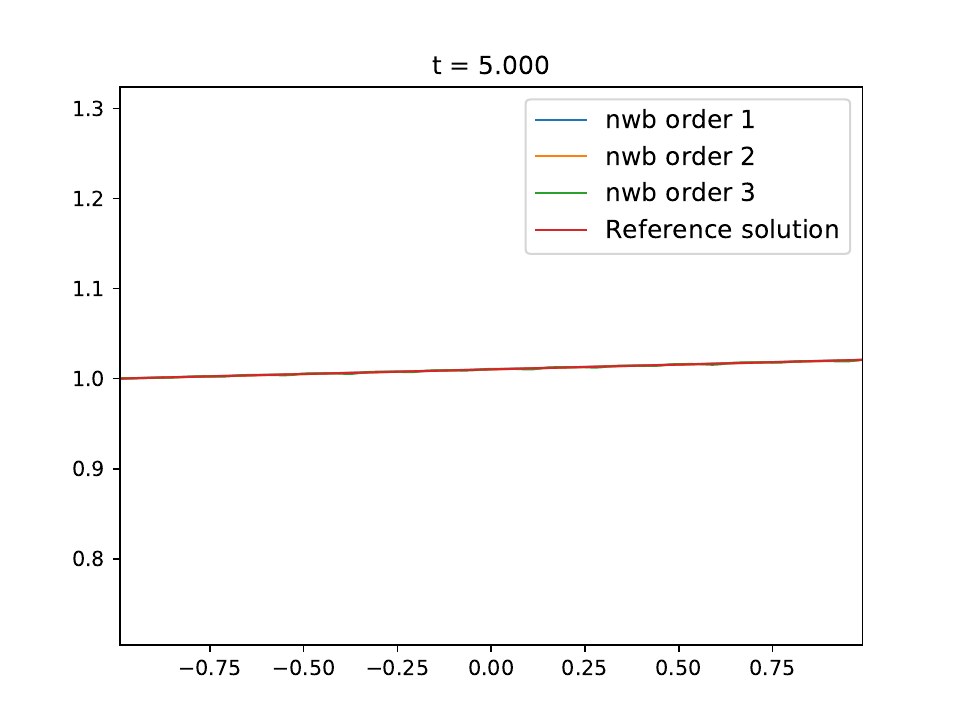}}
 \subfloat[Velocity.]{
   \includegraphics[width=0.33\textwidth]{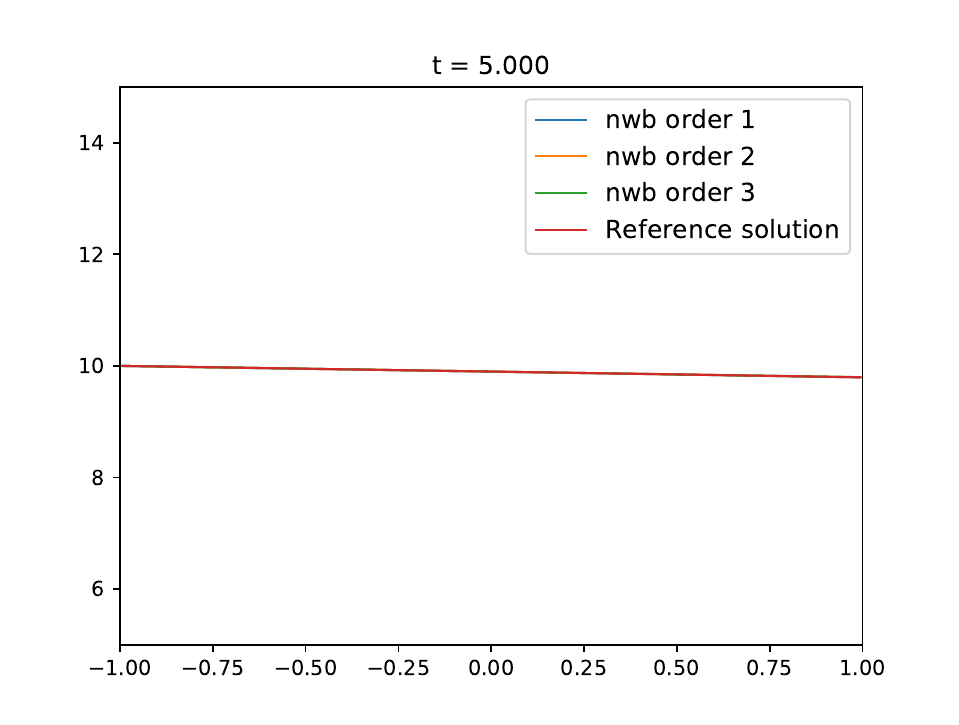}}
   \subfloat[ Energy.]{
   \includegraphics[width=0.33\textwidth]{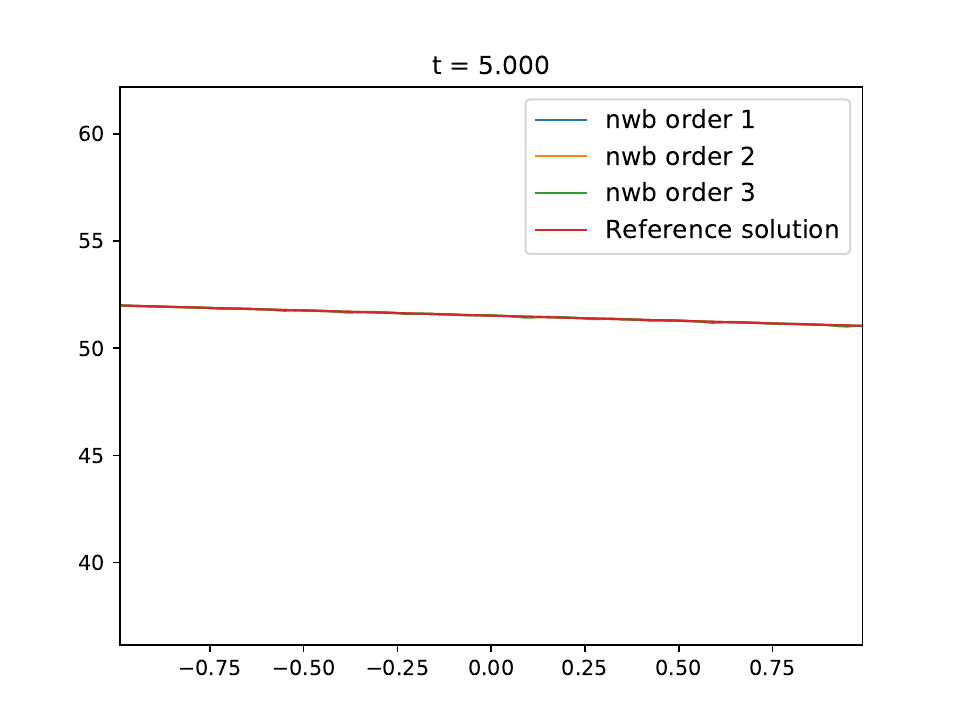}}\vspace{0.0001mm}
   \subfloat[Density. Zoom.]{
   \includegraphics[width=0.33\textwidth]{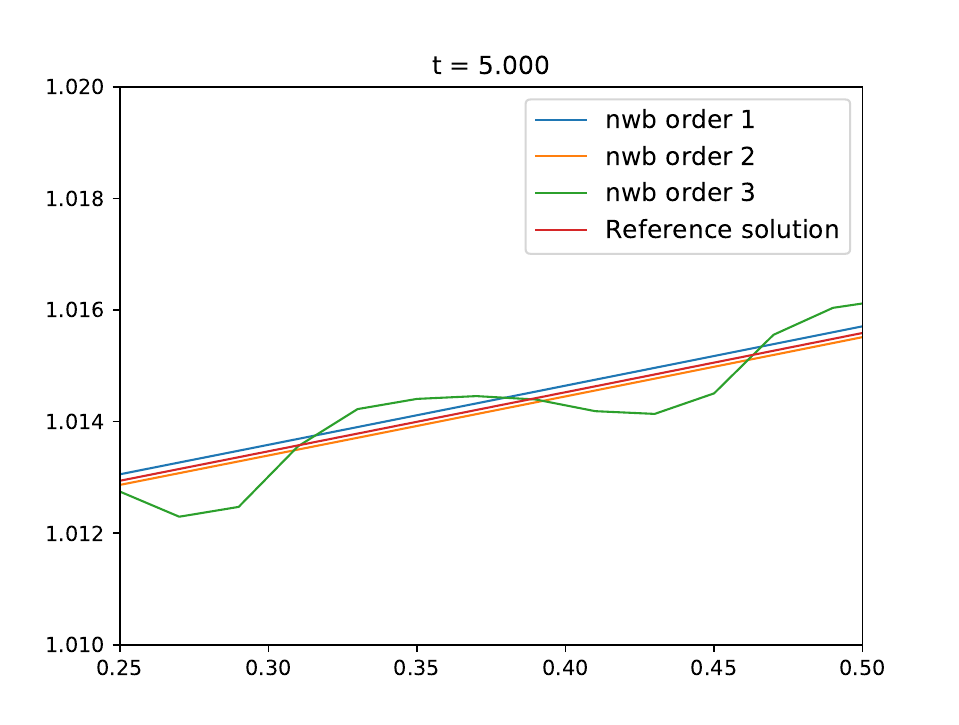}}
   \subfloat[ Velocity. Zoom.]{
   \includegraphics[width=0.33\textwidth]{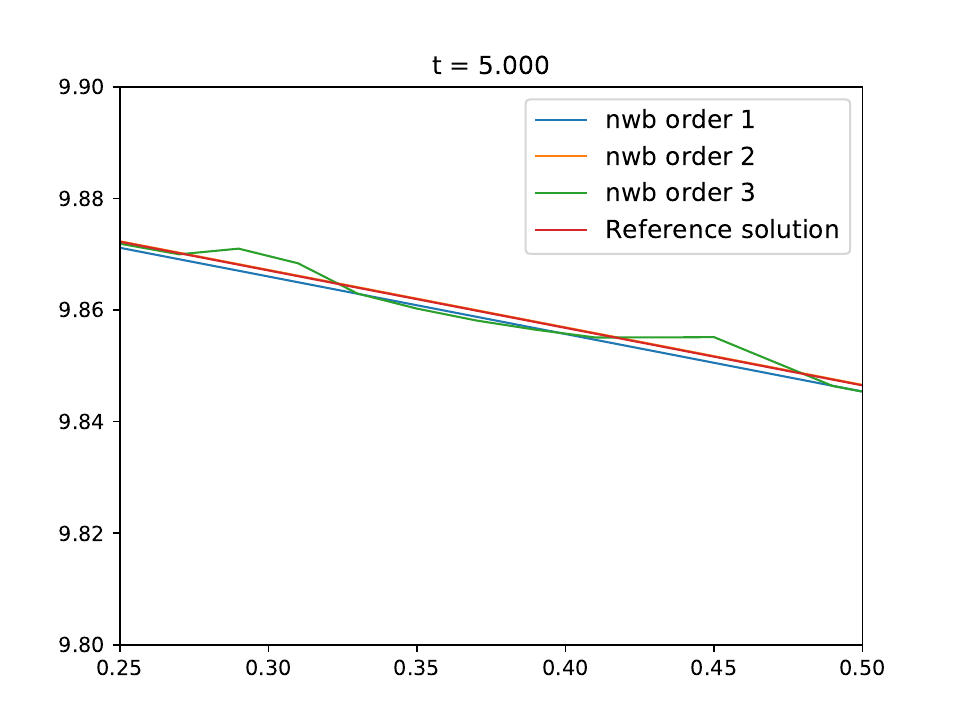}}
   \subfloat[ Energy. Zoom.]{
   \includegraphics[width=0.33\textwidth]{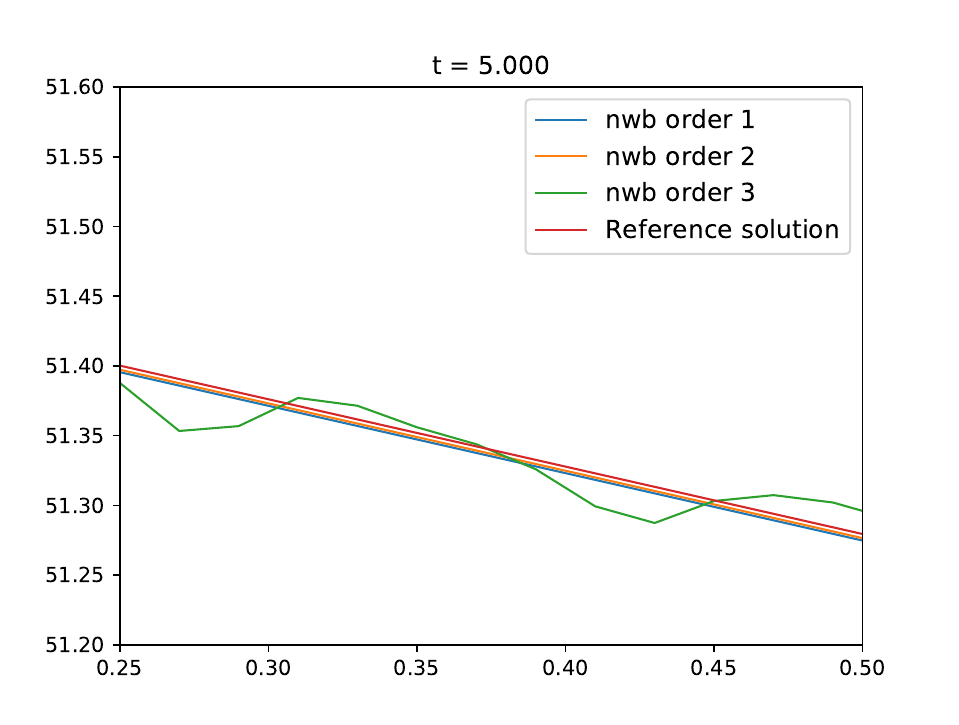}}
    \caption{Test 5.2. Numerical and reference solutions at times $t=0.05$ (up), $t=5s$ (global view (middle) and zoom (down). SM$i$, $i=1,2,3$. Number of cells: 100.} \label{test52_nwb}
 \end{center}
\end{figure}

\begin{figure}[H]
\begin{center}
  \subfloat[Density. ]{
   \includegraphics[width=0.33\textwidth]{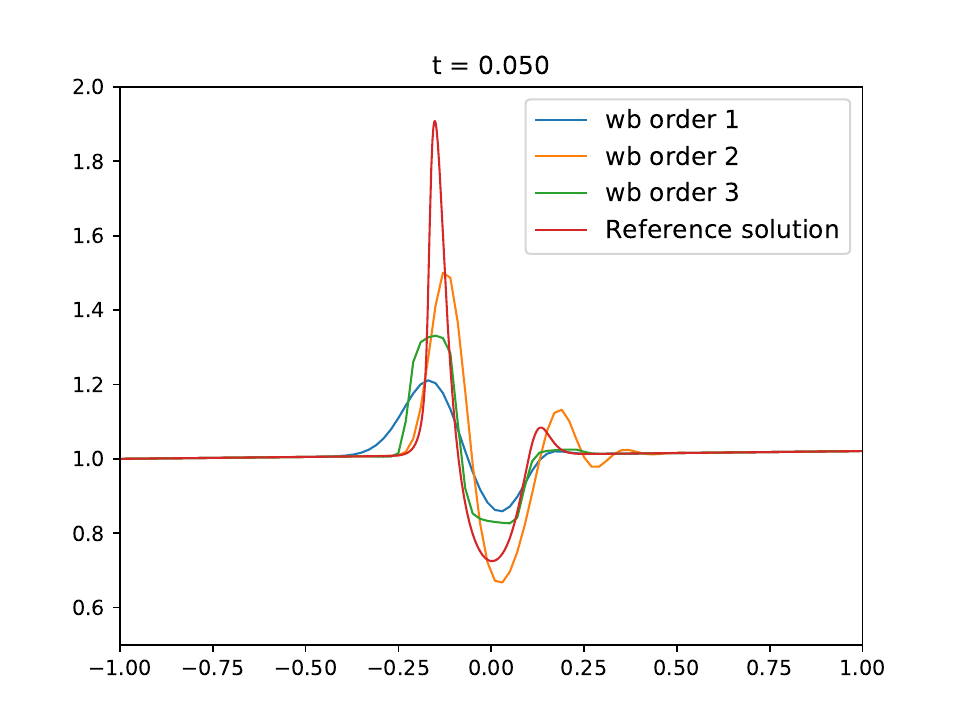}}
 \subfloat[Velocity.]{
   \includegraphics[width=0.33\textwidth]{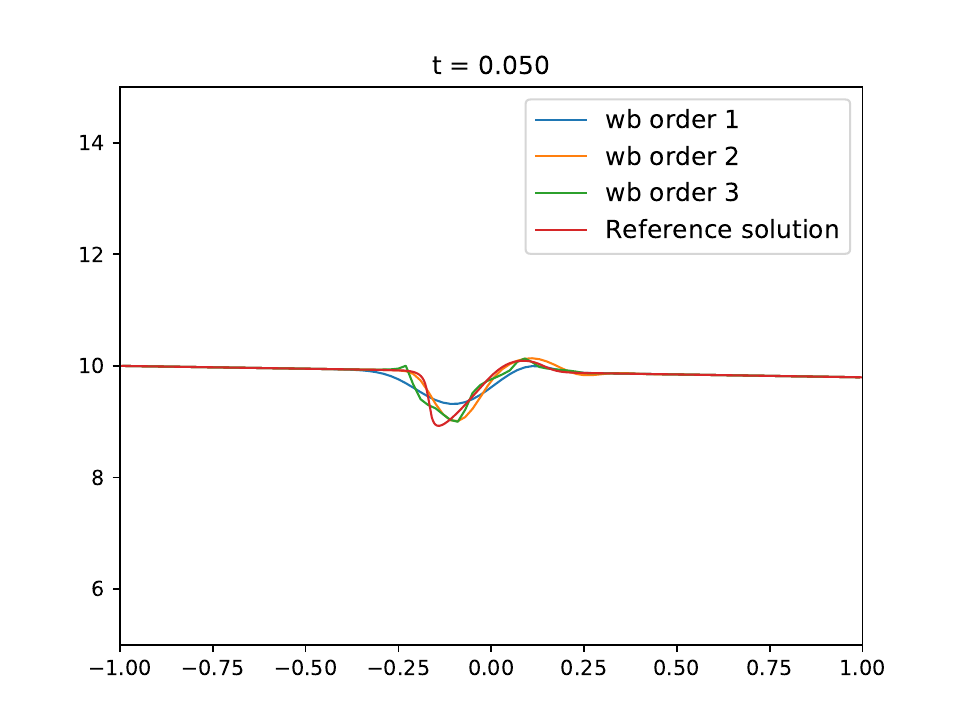}}
   \subfloat[ Energy.]{
   \includegraphics[width=0.33\textwidth]{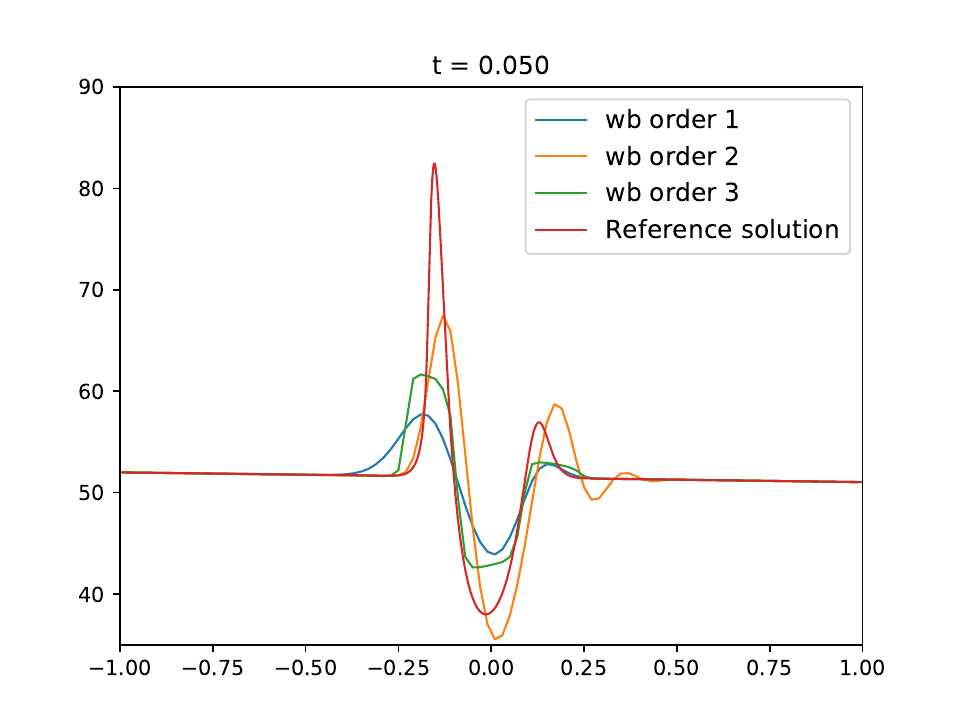}}\vspace{0.0001mm}
   \subfloat[Density. ]{
   \includegraphics[width=0.33\textwidth]{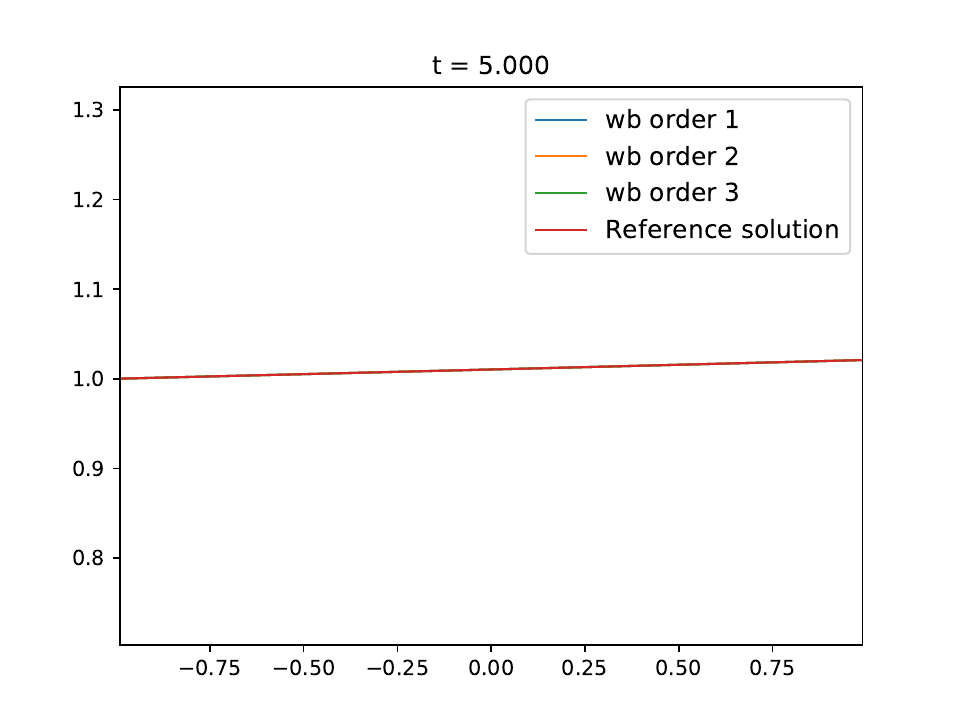}}
 \subfloat[Velocity.]{
   \includegraphics[width=0.33\textwidth]{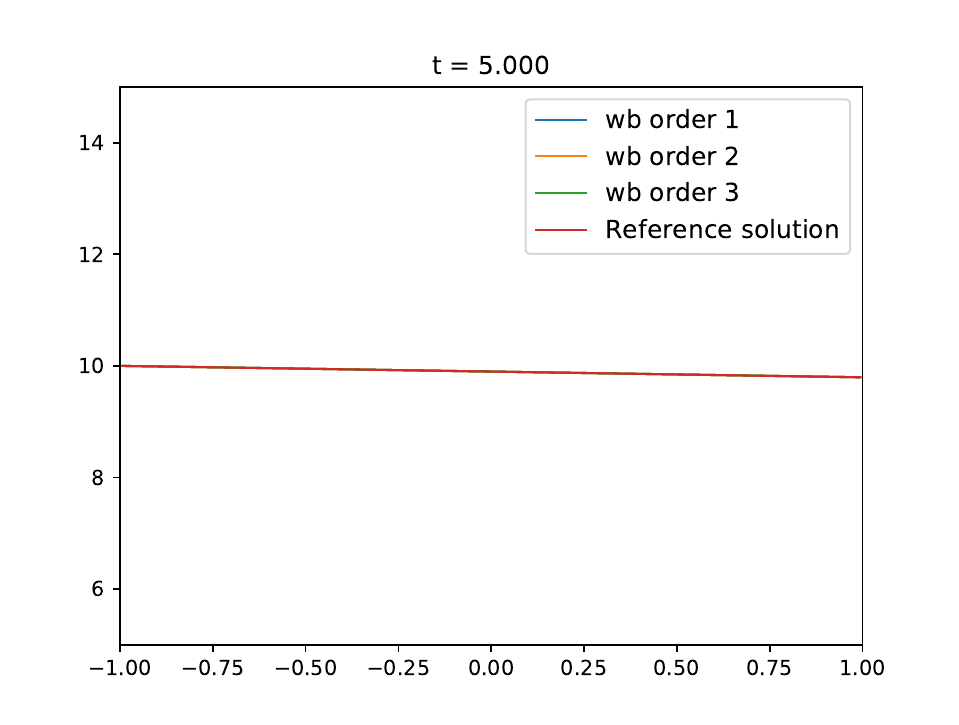}}
   \subfloat[ Energy.]{
   \includegraphics[width=0.33\textwidth]{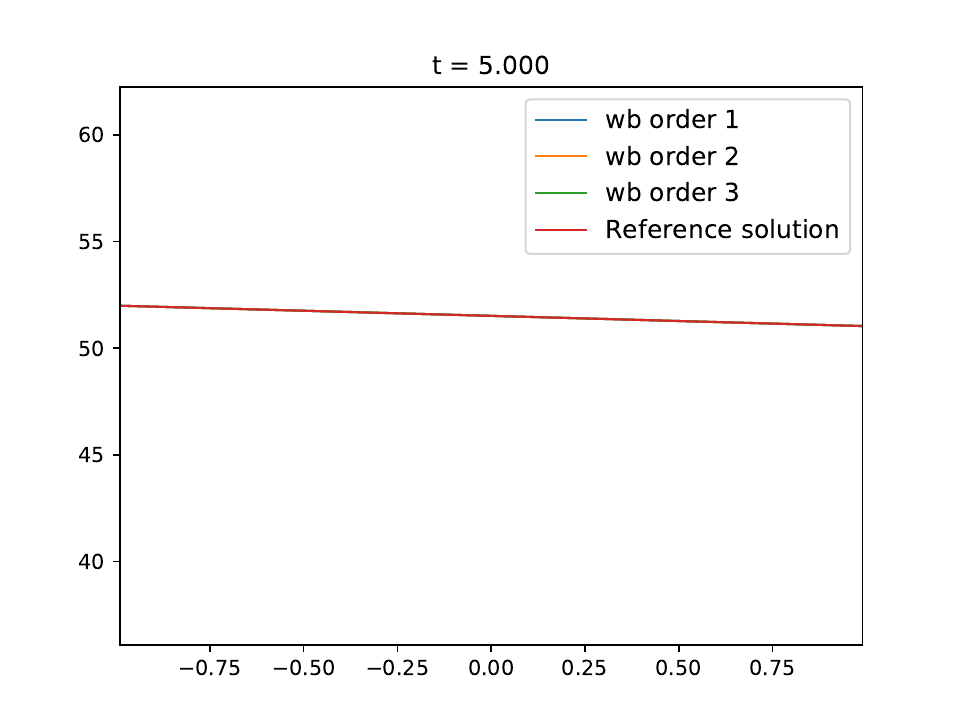}}\vspace{0.0001mm}
   \subfloat[Density. Zoom.]{
   \includegraphics[width=0.33\textwidth]{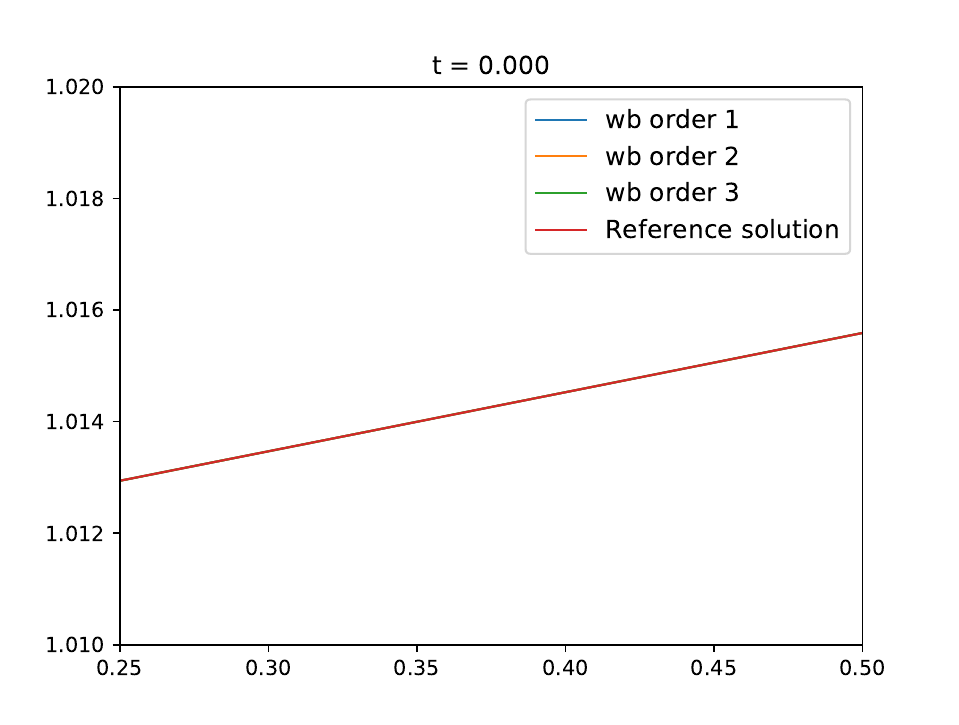}}
   \subfloat[ Velocity. Zoom.]{
   \includegraphics[width=0.33\textwidth]{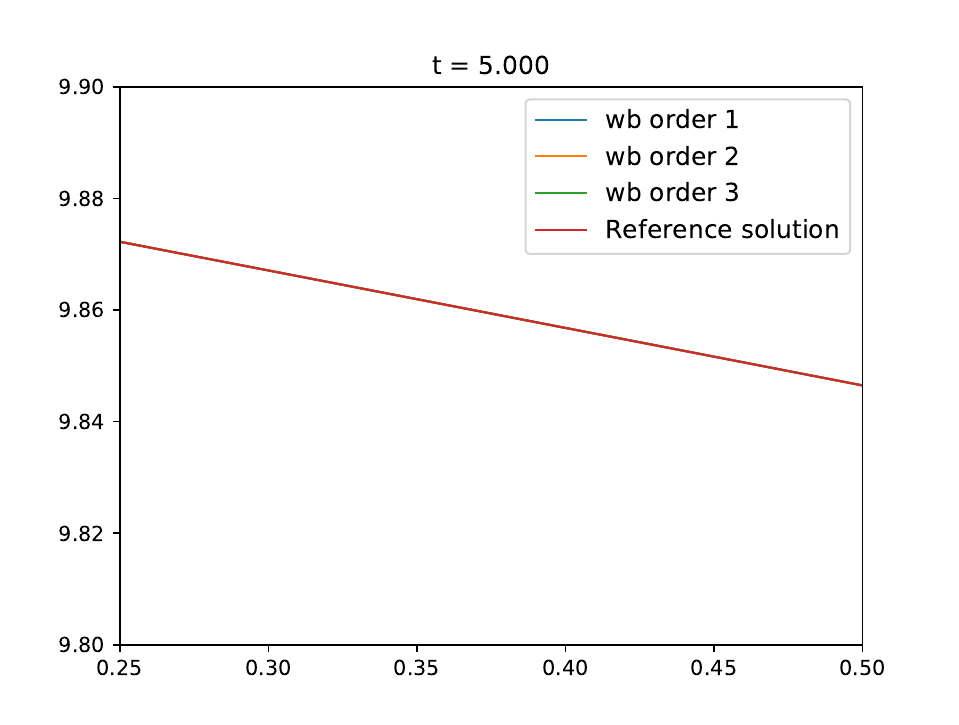}}
   \subfloat[ Energy. Zoom.]{
   \includegraphics[width=0.33\textwidth]{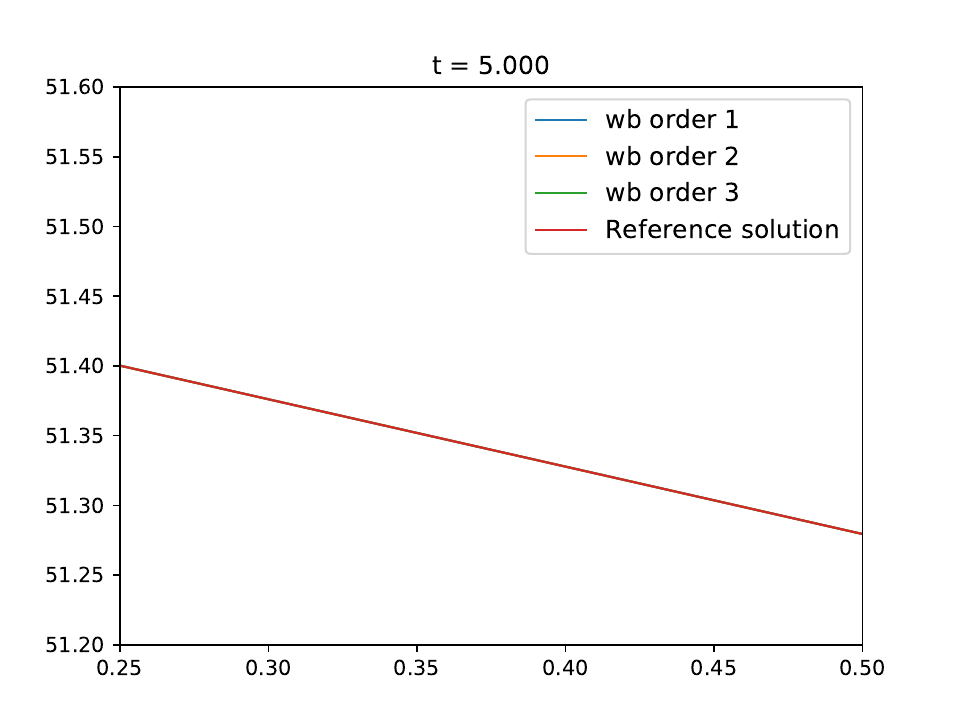}}
    \caption{Test 5.2. Numerical and reference solutions at times $t=0.05$ (up), $t=5s$ (global view (middle) and zoom (down). SM$i$, $i=1,2,3$. Number of cells: 100.} \label{test52_wb}
 \end{center}
\end{figure}

The evolution of a perturbation of the stationary solution considered in the previous test is now simulated. The only difference with Test 5.1. is that now the initial condition is given by:  
$$U_0(x)=U^*(x)+ \begin{pmatrix}
0.3e^{-200(x+0.5)^2} \\
0.0 \\
0.0\\
\end{pmatrix},$$
where $U^*(x)$ is again the stationary solution satisfying $U^*(-1) = (1,10,52)^t $: see Figure \ref{test51cini}.
\begin{figure}[H]
\begin{center}
  \subfloat[ Density. ]{
   \includegraphics[width=0.33\textwidth]{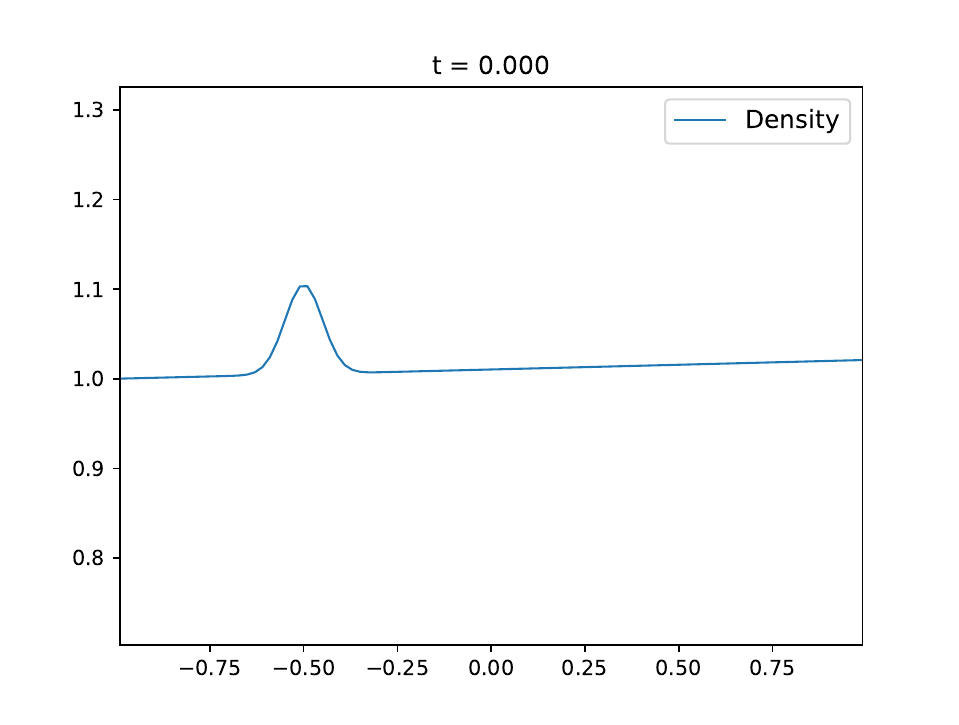}}
 \subfloat[Velocity.]{
   \includegraphics[width=0.33\textwidth]{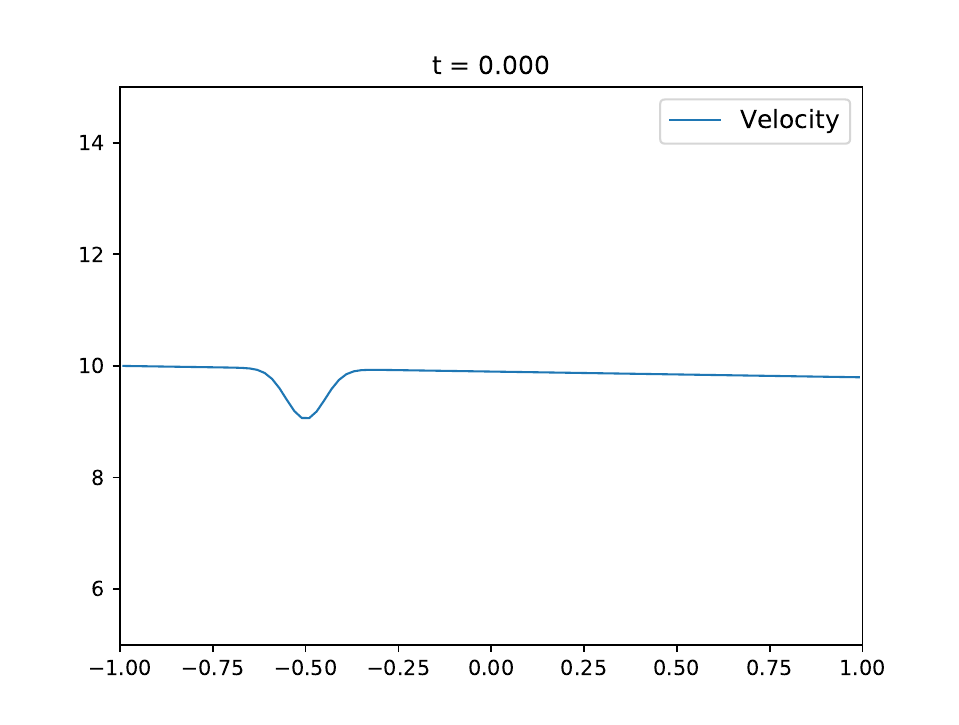}}
    \subfloat[ Energy.]{
   \includegraphics[width=0.33\textwidth]{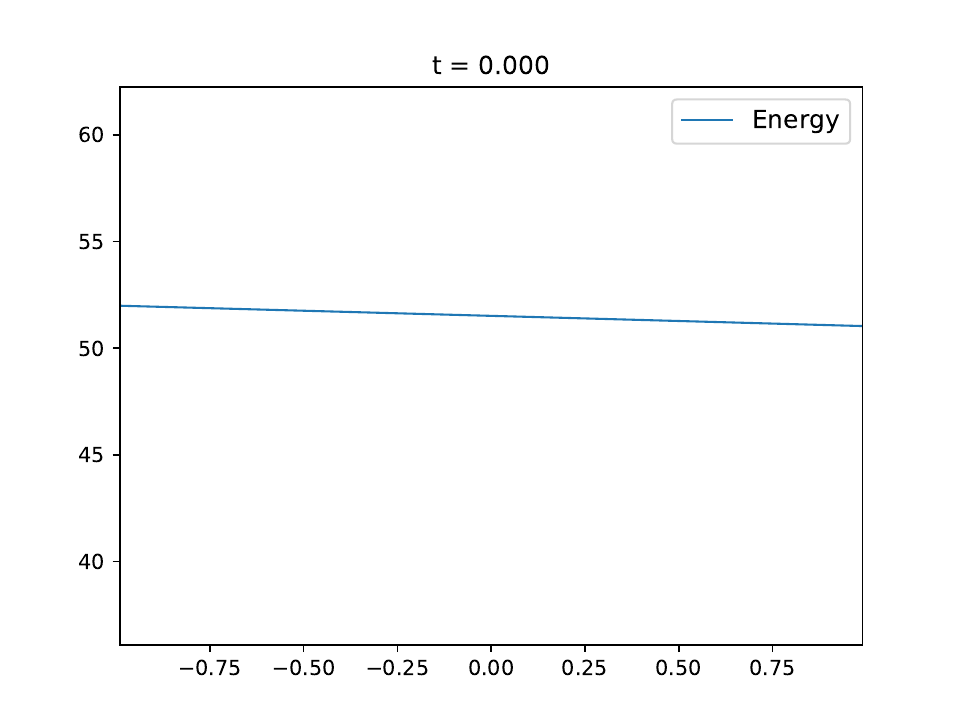}}
    \caption{Test 5.2. Initial condition} \label{test52cini}
 \end{center}
\end{figure}
Figures \ref{test52_nwb} and \ref{test52_wb} show the evolution of the perturbation at times $t=0.05, 5s$ obtained with SM$i$, $i=1,2,3$ and DWBM$i$, $i=1,2,3$.
A reference solution has been computed with a first order well-balanced scheme on a fine mesh (6400 cells). Like in previous cases, it can be observed how the
stationary solution is perturbed by non well-balanced methods. This is confirmed by Tables \ref{test52_error_nwb} and \ref{test52_error_wb}, where the errors at time  $t = 5s$ are shown.

\begin{table}[H]
\centering
\begin{tabular}{|c|c|c|c|} \hline
Cells & \multicolumn{1}{|c|}{Error ($i = 1$)} &\multicolumn{1}{|c|}{Error ($i = 2$)}&\multicolumn{1}{|c|}{Error ($i = 3$)}\\ \hline
  \multicolumn{4}{|c|}{$\rho$}\\\hline
100& 1.45E-4 &2.23E-4 & 7.95E-4 \\
200 & 7.35E-5 & 1.12E-4 & 4.67E-4 \\
400 & 3.70E-5& 5.63E-5 & 3.73E-4 \\
800 & 1.85E-5 & 2.82E-5 & 13.36E-4 \\ \hline
 \multicolumn{4}{|c|}{$q$}\\\hline
100& 1.28E-4 & 9.35E-5 & 6.31E-3 \\
200 & 6.41E-4 & 4.16E-5 & 3.90E-3 \\
400 & 3.21E-4 & 1.95E-5 & 3.45E-3 \\
800 & 1.61E-4 & 9.43E-6 & 3.01E-3 \\ \hline 
 \multicolumn{4}{|c|}{$E$}\\\hline
100& 5.96E-3  &9.49E-3 & 2.82E-2 \\
200 & 2.97E-3 & 4.73E-3 & 1.82E-2 \\
400 & 1.48E-3 & 32.36E-3 & 1.76E-2 \\
800 & 7.40E-4 & 1.18E-3 & 81.49E-2 \\ \hline

\end{tabular}
\caption{Test 5.2. Errors in $L^1$ norm for SM$i$, $i=1,2,3$. $t = 5s$} \label{test52_error_nwb}

\end{table}

\begin{table}[H]
\centering
\begin{tabular}{|c|c|c|c|} \hline
Cells & \multicolumn{1}{|c|}{Error ($i = 1$)} &\multicolumn{1}{|c|}{Error ($i = 2$)}&\multicolumn{1}{|c|}{Error ($i = 3$)}\\ \hline
  \multicolumn{4}{|c|}{$\rho$}\\\hline
100& 9.45E-15  & 2.52E-14 & 3.27E-14 \\
200 & 1.73E-14  & 3.53E-14 & 5.25E-14 \\
400 & 3.71E-14 & 5.04E-14 & 7.82E-14 \\
800 & 4.37E-14  & 7.43E-13 & 1.75E-13 \\ \hline
 \multicolumn{4}{|c|}{$q$}\\\hline
100& 9.03E-14 & 2.39E-13 & 3.15E-13 \\
200 & 1.71E-13 & 3.37E-13 & 5.54E-13 \\
400 & 3.64E-13 & 4.91E-13 & 8.53E-13 \\
800 & 4.30E-13 & 7.30E-13 & 1.69E-12 \\ \hline 
 \multicolumn{4}{|c|}{$E$}\\\hline
100& 5.44E-13  & 1.18E-12 & 1.50E-12 \\
200 & 8.93E-13 & 2.72E-12 & 2.47E-12 \\
400 & 1.84E-12 & 3.46E-12 & 4.79E-12 \\
800 & 2.27E-12 & 3.75E-12 & 8.52E-12 \\ \hline

\end{tabular}
\caption{Test 5.2. Errors in $L^1$ norm for DWBM$i$, $i=1,2,3$. $t = 5s$} \label{test52_error_wb}

\end{table}

\section{Conclusions}

The strategy introduced in \cite{sinum2008} has been followed to derive a family of high-order well-balanced numerical methods that can applied to general 1d systems of balance laws. The main difficulty in applying these methods comes from the first stage of the well-balanced reconstruction procedure: at every cell and at every time step a nonlinear problem has to be solved consisting in finding a stationary solution whose average is the given cell value. This problem has been interpreted as a control one related to an ODE system, in which the constraint is the given average and the control is the initial condition. The problem has been written in functional form, the gradient of the functional has been computed with the help of the adjoint system, and Newton's method can be then applied. The effects of the use of a quadrature formula to compute the cell averages and the integral source terms have been analyzed and numerical techniques have been introduced to preserve the well-balancedness of the methods. In particular, for first and second order methods, the use of the midpoint rule allows one to reduce the control problems to standard Cauchy problems.

In order to test the efficiency and the well-balancedness of the methods, they have been applied to a number of systems of balance laws ranging from academic tests systems consisting of Burgers equations with nonlinear source terms to flow models like the shallow water system or Euler equations of gas dynamics with gravity effects. In some cases the stationary solutions are known either in implicit or explicit form while in others the only information comes from the ODE that the stationary solutions solve: the former allow us to compare the efficiency of the new implementation while the latter allow us to show the generality of the methods. In particular, it is the first time, to the best of our knowledge, that a family of high-order methods that preserve moving stationary solutions for Euler equations with gravity have been designed. 

The tests put on evidence that the well-balanced modification increases the computational cost, specially for methods of order bigger than three. In any case, this extra computational cost is lower than the one that would require to lead the discretization errors to (close to) zero machine by refining the mesh or increasing the order of non-well-balanced methods. On the other hand, in cases in which the explicit form of the stationary solution is known, the computational cost of an implementation based on control techniques is of the same order --or even lower in some cases-- than an implementation based on the analytic expression of the solution of the non-linear problems related to the well-balanced reconstruction. 

Further developments include applications of the introduced technique to:

\begin{itemize}
    \item Systems of balance laws \eqref{sle} in which the function $H$ has jump discontinuities.
    \item Transcritical stationary solutions.
    \item Multidimensional problems.
    \end{itemize}

\end{document}